\documentclass[10pt]{amsart}
\usepackage{graphicx}
\usepackage{amsthm,amsfonts,amssymb,amsmath,amscd}
\usepackage[rotate,arrow,matrix]{xy}

\DeclareMathOperator{\id}{id}

\newtheorem{thm}{Theorem}[section]

\newtheorem{prop}[thm]{Proposition}
\newtheorem{claim}[thm]{Claim}

\theoremstyle{definition}
\newtheorem{df}[thm]{Definition}

\theoremstyle{remark}

\begin{document}

\title{On the Homology of Open-Closed String Field Theory}
\author{Eric Harrelson}
\date{\today}
\maketitle

\begin{abstract}
The homology of a 2-colored dioperad of decorated Riemann surfaces, relevant to open-closed string field theory, is computed. The structure it describes is realized in an open-closed setting of string topology via an action at the level of topological spaces.
\end{abstract}

\section*{Introduction}

In Zwiebach's study of oriented open-closed string theory \cite{zwie}, he considered a certain moduli space of Riemann Surfaces with boundary having closed punctures in the interior and open punctures on the boundary coming with parameterizations by the unit disk and upper half disk.  While he did not consider it as such, this moduli space forms a 2-colored PROP. In fact, an open-closed CFT (with one D-brane) can be defined as an algebra over this PROP, and an open-closed Topological CFT (or string background) can be defined as an algebra over its chains.  The first purpose of this paper is to describe completely the homology of the biggest genus 0 structure inside this PROP.  The operad inside of this PROP formed by spheres with no boundary is well known to be homotopy equivalent to the framed little disks operad.   It is shown by Getzler that its homology  defines a BV-algebra \cite{getz1}.  This extends the result by F. Cohen \cite{fcohen} showing that the  homology of the non-framed little disks operad describes a Gerstenhaber algebra.  In \cite{vor4}, Voronov invented the Swiss-cheese operad which is a (non framed) finite dimensional model of the operad inside this PROP formed by Riemann spheres with one or no boundary components.  He computed its homology and calls the algebra that it defines a Swiss-cheese algebra.  The algebra is defined on a pair of graded vector spaces $(V_C,V_O)$ and consists of a Gestenhaber structure on $V_C$, an associative multiplication on $V_O$ , and an algebra action of $V_C$ on $V_C$.  The framed version of Swiss, h.e. to the subspace of spheres with one or no boundary components, has the same result except that $V_C$  is a BV-algebra.

The biggest genus 0 operad inside this PROP is formed by all spheres with boundary having exactly one puncture labeled as an output.  Its homology contains the structure of (framed) Swiss-cheese.  However, there is a bigger genus 0 structure in this PROP containing this operad.  The subspace of all genus 0 surfaces, with an arbitrary number of inputs and outputs, forms what's called a {\em dioperad}, invented by Gan in \cite{gan}.  A dioperad only considers compositions which attach one input to one output so as to create no genus.  The first four sections of this paper together give a complete description of the homology of this dioperad.

Restricting the description of the dioperad to the generators with one output and the relations only involving them gives a description of the biggest genus 0 operad in the PROP.  Sec 5 considers what extra structure is given by the cyclic structure of this operad and the semi-modular structure given by self sewing open punctures on the same boundary component.

The homology of the free loop space of an oriented manifold $M$, $LM$ , was shown to form a BV-algebra in Chas and Sullivan's String topology \cite{chas1}.  In \cite{vor1}, Voronov invented the $Cacti$ operad and announced that it is h.e to the framed little disks operad.  He showed how we can obtain $H_*(LM)$ as an algebra over $H_*(Cacti)$ via a geometric action of $Cacti$ on $LM$.  In \cite{sull}, Sullivan considers the algebraic structure of open-closed string topology.  In particular, he studies the homology of $PM_K$, the space of paths in a manifold starting and ending in a fixed submanifold $K$.  The purpose of the last section is to extend the Cacti result to this open-closed setting by defining an open-closed version of Cacti, {\em OC Cacti}, having the same homology as the biggest genus 0 operad inside the PROP, and showing how to obtain the pair $(H_*(LM),H_*(PM_K))$ as algebra over $H_*(OC \, Cacti)$.

In my thesis, I have extended this open-closed Cacti operad to an open-closed 2-colored graph PROP modeling the entire moduli space and acting in string topology.  This extends the Sullivan chord diagrams used by R. Cohen and V. Godin in \cite{cohen} which act in closed string topology.  I have also extended to the setting where we consider a set of submanifolds of $M$, rather than just one, and all the spaces of paths starting in one submanifold and ending in another (see the last remark sec. 6).  I have recently become aware that some similar work involving an open-closed PROP and string topology is being done by A. Ramirez for his thesis \cite{rami}.

For good sources discussing the topics of operads and PROPs and their use in physics, see \cite{cohen2}, \cite{mark}, \cite{vor1}, \cite{vor2}, and \cite{vor3}.

\section*{Acknowledgements}
I would like the thank Sasha Voronov and Jim Stasheff for many helpful conversations and suggestions.

\tableofcontents

\section{Description of 2-colored dioperad}

Consider the moduli space of genus 0 Riemann surfaces with boundary.  RS with boundary means for us a complex surface based on the closed upper half plane.  Add to this punctures which can be in the interior or on a boundary component.  Each puncture is designated as an input or output.  The inputs are labeled $(1_i,...,n_i)$  and the outputs are labeled $(1_o,...,m_o)$.  It is required that there be at least one output  A puncture in the interior of a surface, called a {\em closed puncture}, comes with an analytic parametrization given by a biholomorphic mapping of the standard disk   into the surface sending  to the puncture.  A puncture in the boundary, called an {\em open puncture}, comes with a biholomorphic mapping of the half disk   into the surface sending the real axis into the bndry and  to the puncture.  These parameterizations may only overlap in their boundaries.

\begin{df}

The `open-closed' moduli space described above, denoted $OC={OC(n,m)}$, forms a 2-colored dioperad (see introduction) via sewing closed punctures to closed punctures and open to open using their parameterizations.  This is done in the standard way, using $w=1/z$ for closed punctures and $w=-1/z$ for open.

\end{df}

When sewing closed to closed, the number of boundary components in the resulting surface is the sum of the boundary components.  In the open to open case, its the sum minus one.

% \newpage

\underline{\bf Sewing closed to closed}

\begin{center}
\includegraphics[width=3in]{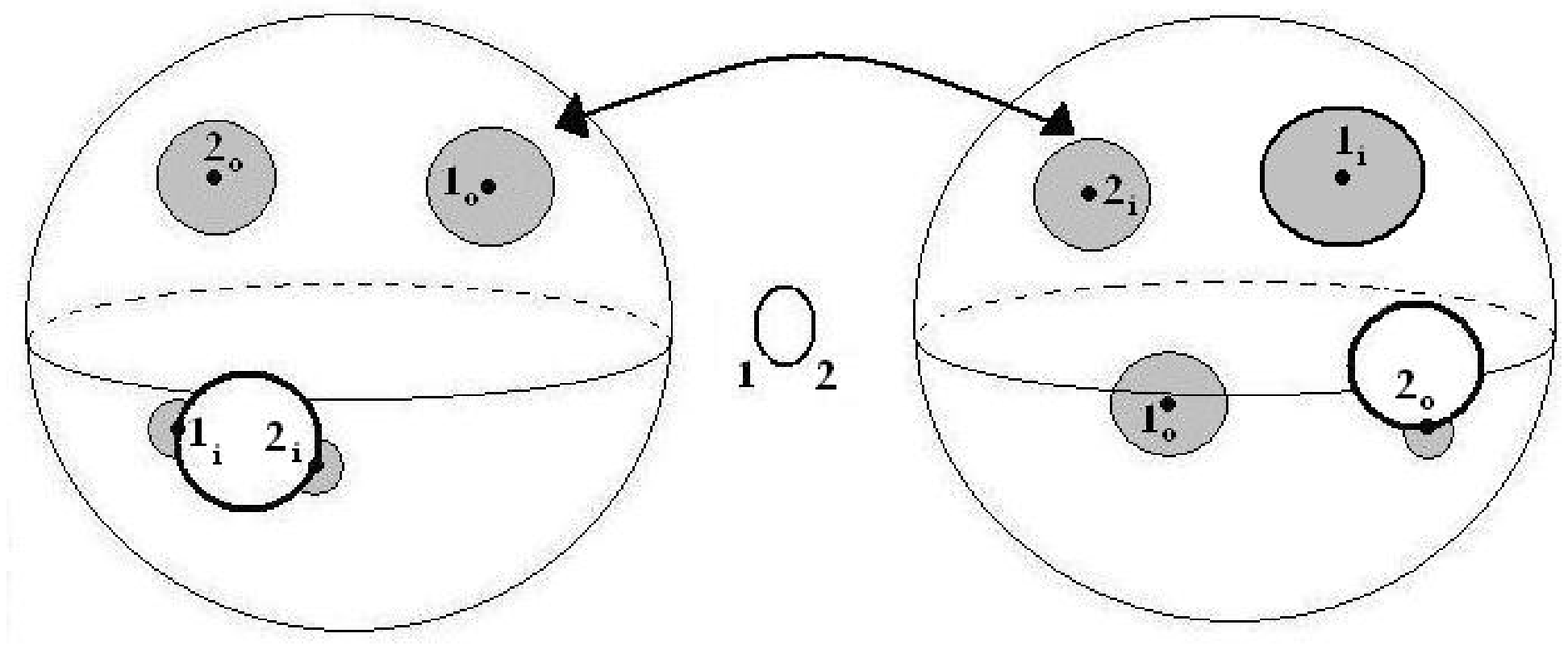} \par
\includegraphics[width=3in]{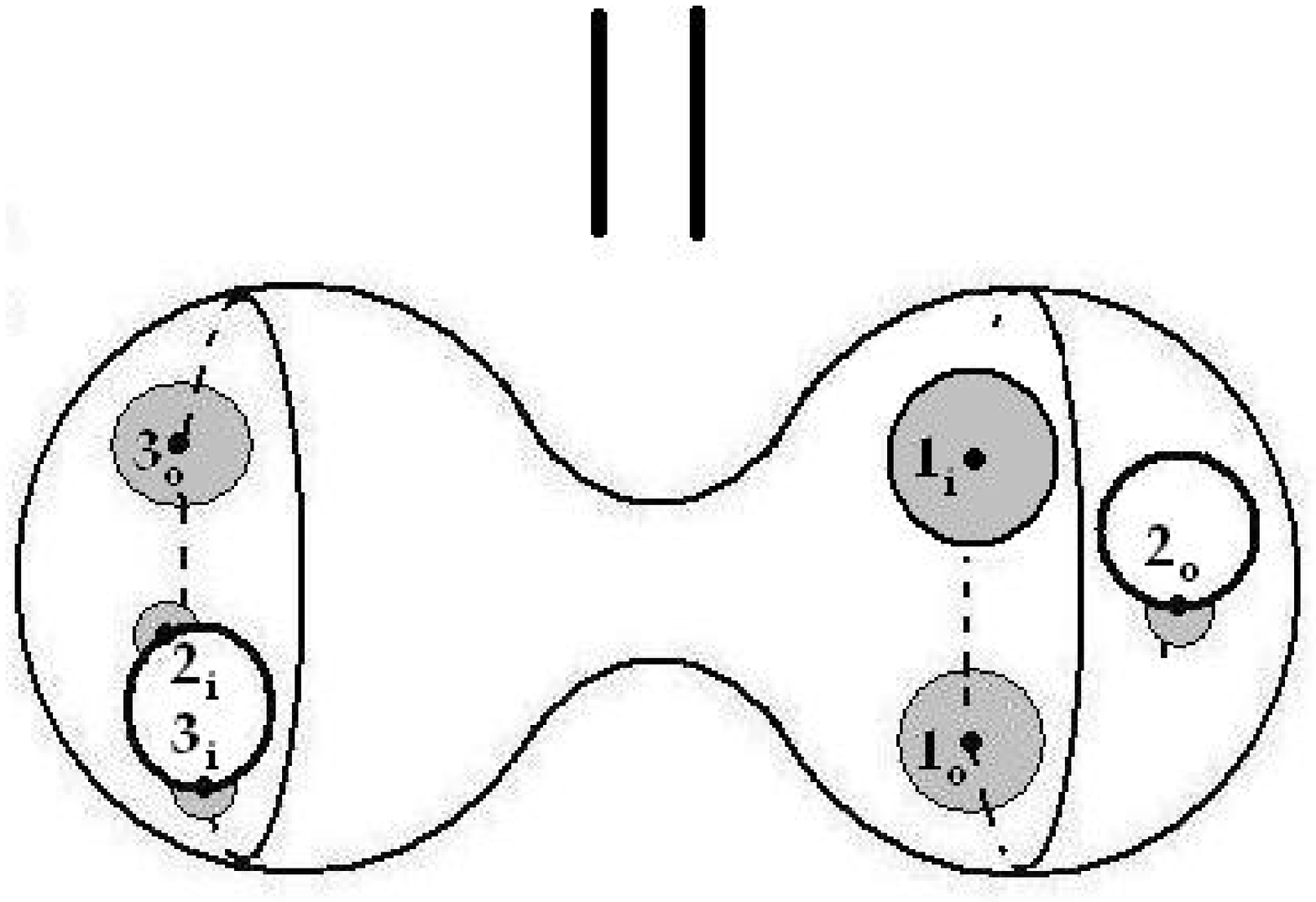}
\end{center}

\underline{\bf Sewing open to open}

\begin{center}
\includegraphics[width=3in]{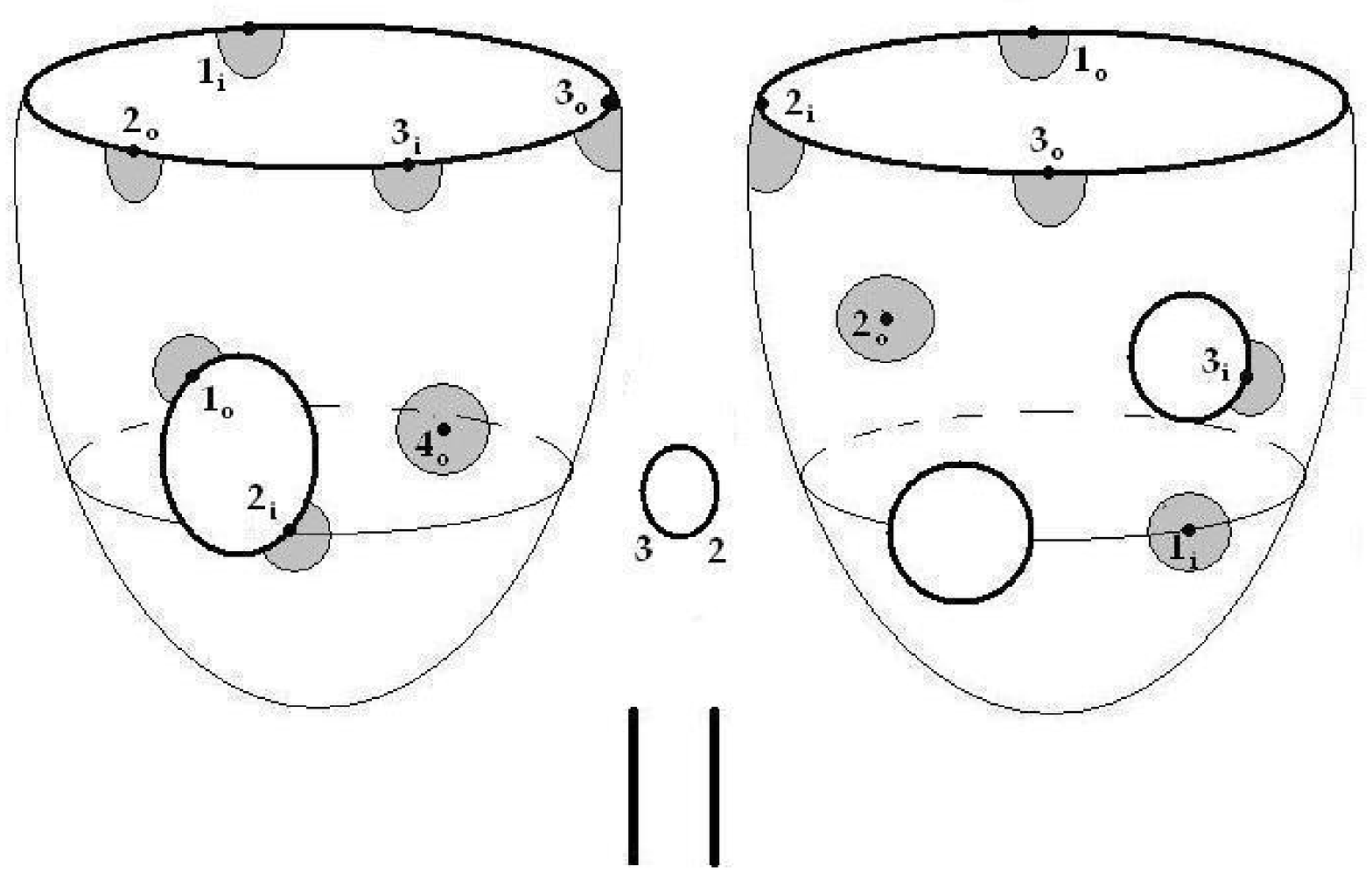} \par
\includegraphics[width=3in]{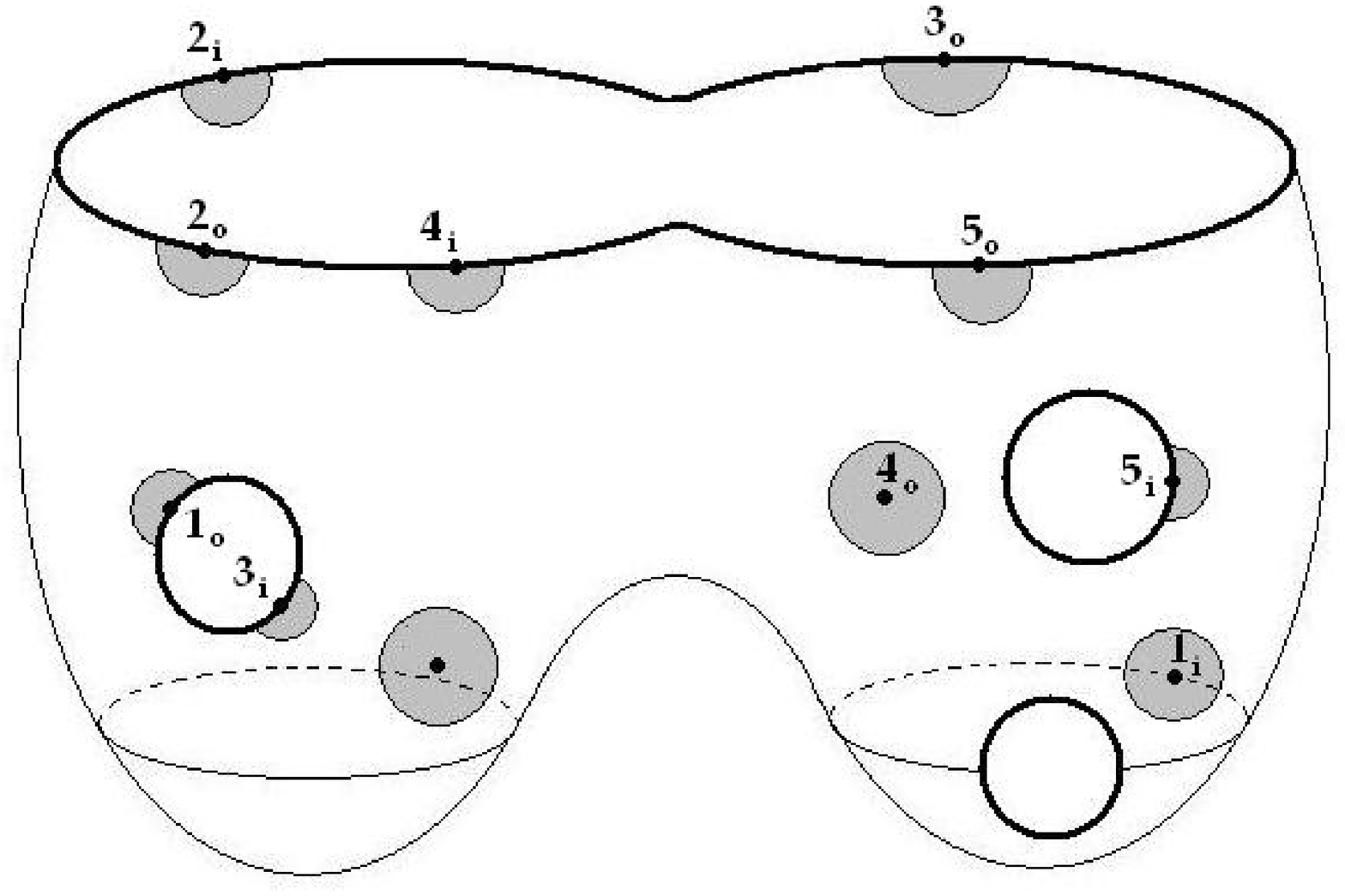}
\end{center}

\section{Description of $H_0(OC)$}

Lets start with a description of the path components of the moduli space and how composition (sewing) acts on them.  For each point in $OC$ we have the following data:
\begin{enumerate}
\item  A subset $C \subseteq {1_i,...,n_i,1_o,...,m_o}$ consisting of the labels of all closed punctures.
\item An unordered list of cyclically ordered subsets $(a_1,...,a_{k_1}),(b_1,...,b_{k_{2}})$, etc grouping together labels of open punctures lying on the same boundary component and giving them the cyclic order induced by the orientation of the boundary component (which is induced by the canonical orientation of the Riemann surface).  An empty parentheses $()$ is used for boundary components with no open punctures.
\end{enumerate}

\begin{df}
Its clear that two points are in the same path component iff they have the same data.  Call this data the {\em type} of the path component.
\end{df}

\underline{\bf Example} \par

\includegraphics[width=2in]{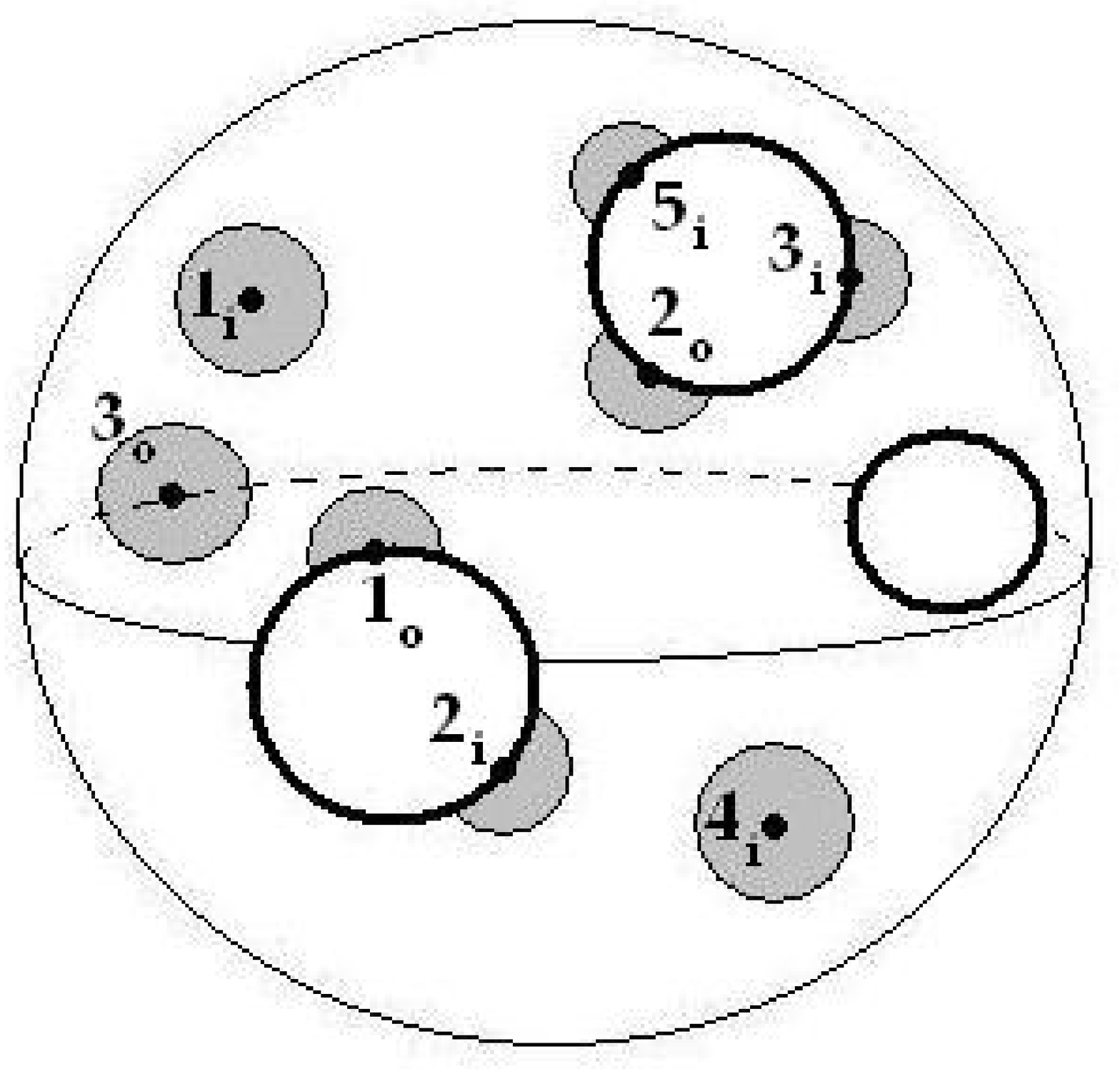} $\begin{array}[b]{c} Type: {1_i,4_i,3_o},(1_o,2_i),(2_o,5_i,3_i),() \\
                                                                \vspace{.75in} \end{array}$

It is also clear from the pictures for sewing how composition acts on the path components.  Note that the path components of type ${1_i,1_o}$  and $(1_i,1_o)$ are the identities for closed and open composition.

\begin{prop}
The following path components generate $H_0(OC)$.  (Listed with them are the degree 0 trees that we'll use to represent them and the corresponding degree 0 operations in an algebra $(V_c,V_o)$  over $H_0(OC)$  ).
\end{prop}

% $\begin{array}{lcl}
% \begin{array}{l} Type : {1_i,2_i,1_o} \\ \vspace{.4in} \end{array} & \hspace{.4in} \includegraphics[height=.8in]
% {fig4-1.ps} \hspace{.4in} & \begin{array}{l} m_c:V_c \otimes V_c \mapsto V_c \\ \vspace{.4in} \end{array}
% \end{array}$

\begin{description}

\item[g1] {\large\bf closed multiplication} $$ \begin{array}{lll} Type : {1_i,2_i,1_o} &  \hspace{1in} \put(0,-20){\includegraphics[height=.8in]{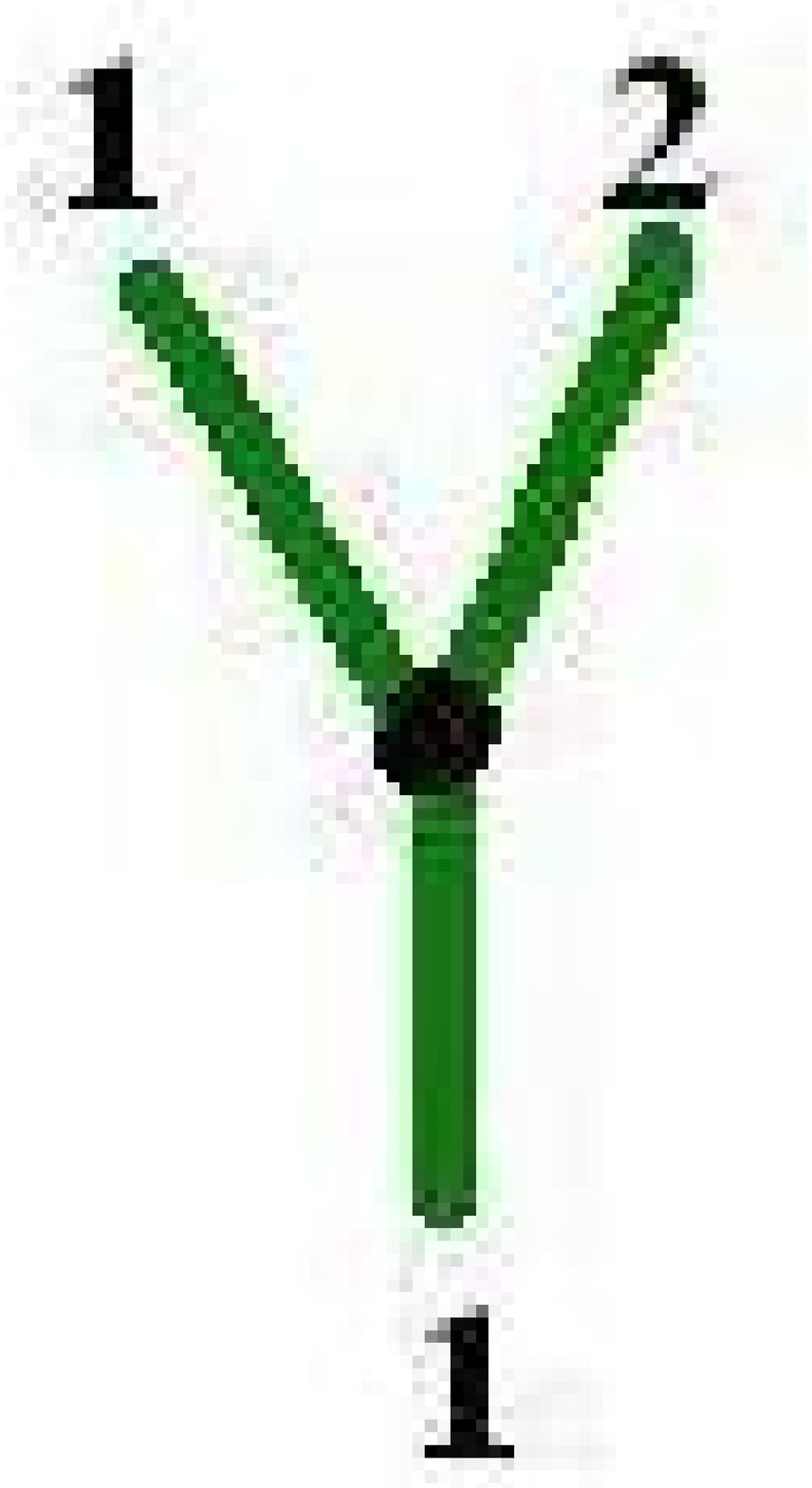}} & \hspace{1in} m_c:V_c \otimes V_c \mapsto V_c \end{array} $$

\item[g2] {\large\bf open multiplication} $$ \begin{array}{lll} Type : (1_i,1_i,1_o) & \hspace{1in} \put(0,-20){\includegraphics[height=.8in]{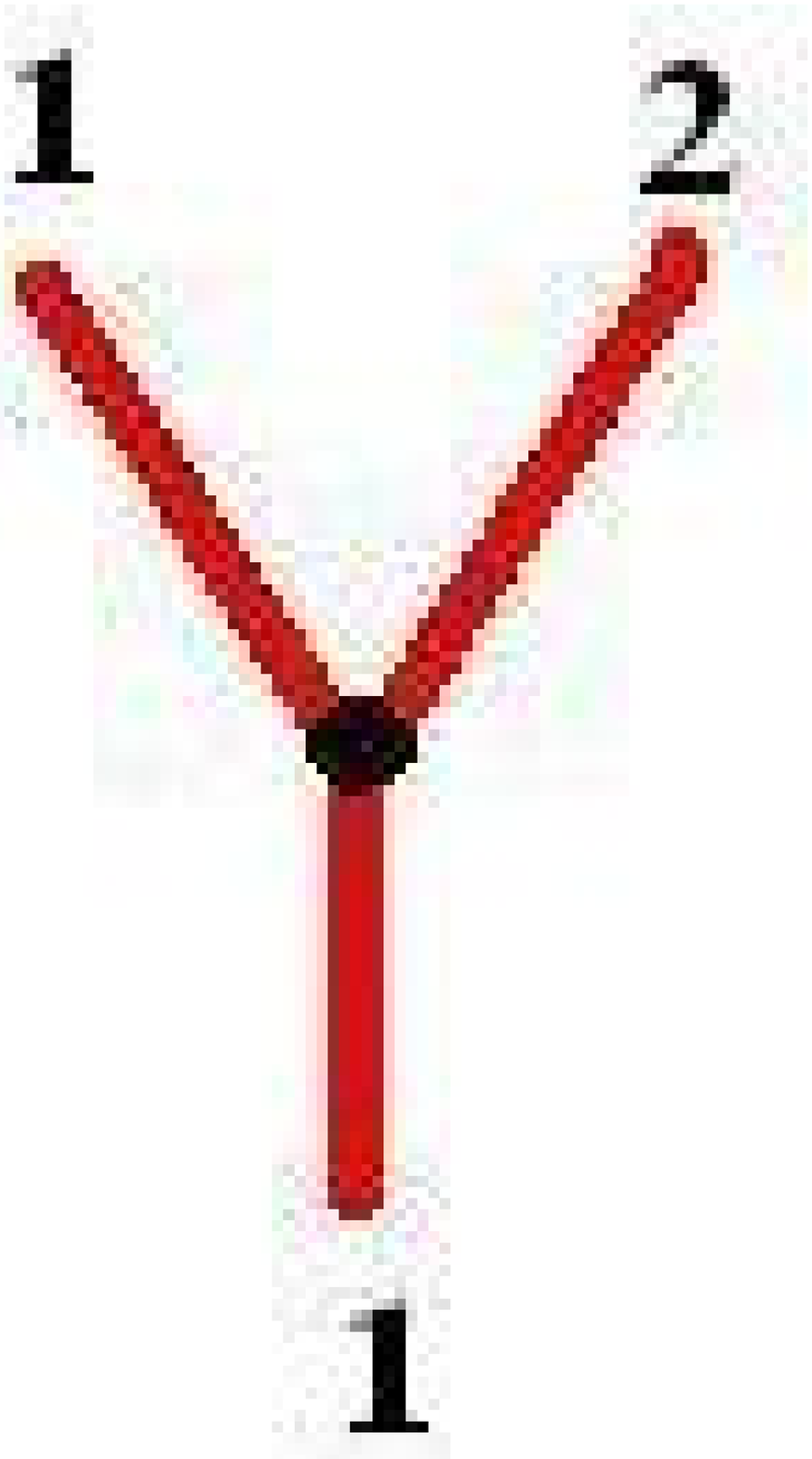}} & \hspace{1in} m_o:V_o \otimes V_o \mapsto V_o \end{array} $$

\item[g3]{\large\bf closed unit} $$ \begin{array}{lll} Type : {1_o} & \hspace{1in} \put(0,-10){\includegraphics[height=.4in]{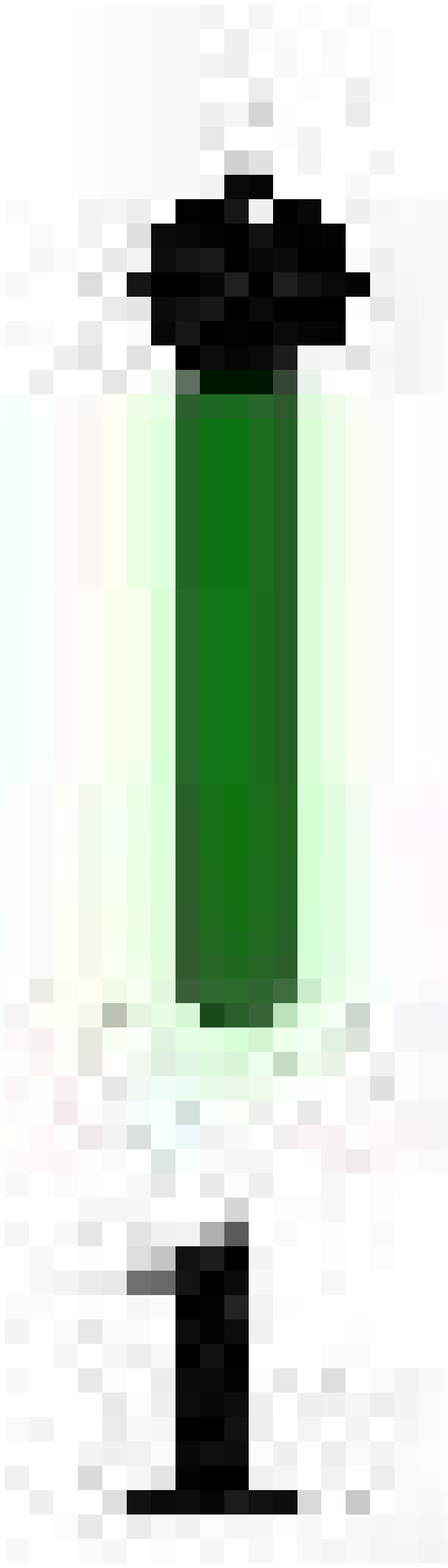}} & \hspace{1in} K \mapsto V_c \; 1 \mapsto e_c \end{array} $$

\item[g4]{\large\bf open unit} $$ \begin{array}{lll} Type : (1_o) & \hspace{1in} \put(0,-10){\includegraphics[height=.4in]{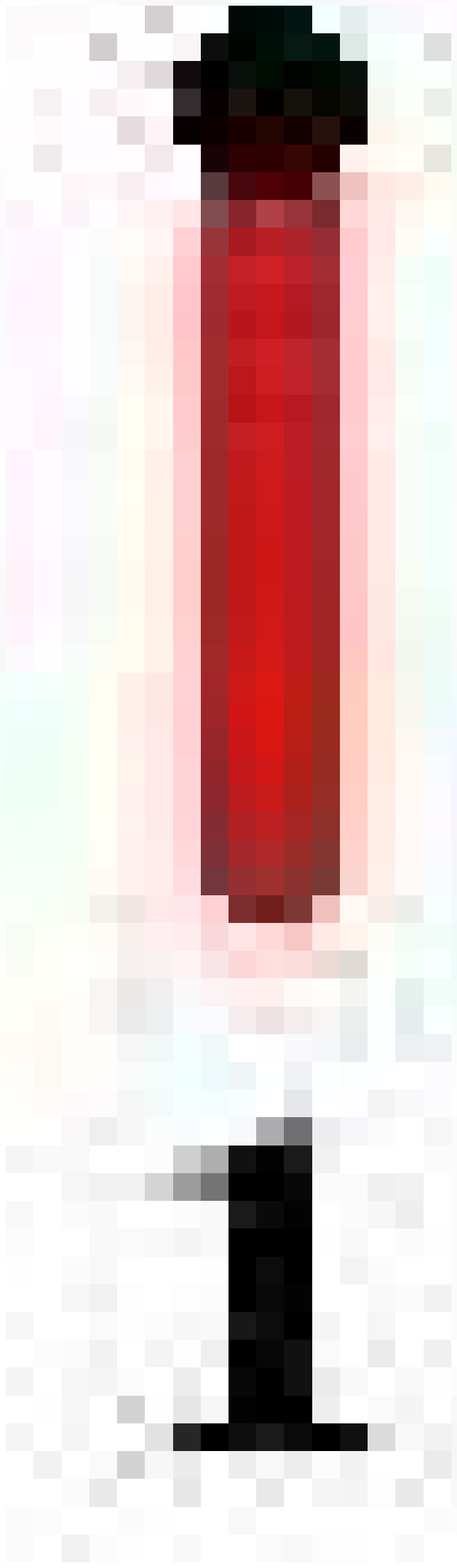}} & \hspace{1in} K \mapsto V_o \; 1 \mapsto e_o \end{array} $$

\item[g5]{\large\bf closed to open} $$ \begin{array}{lll} Type : {1_i}(1_o) & \hspace{1in} \put(0,-20){\includegraphics[height=.8in]{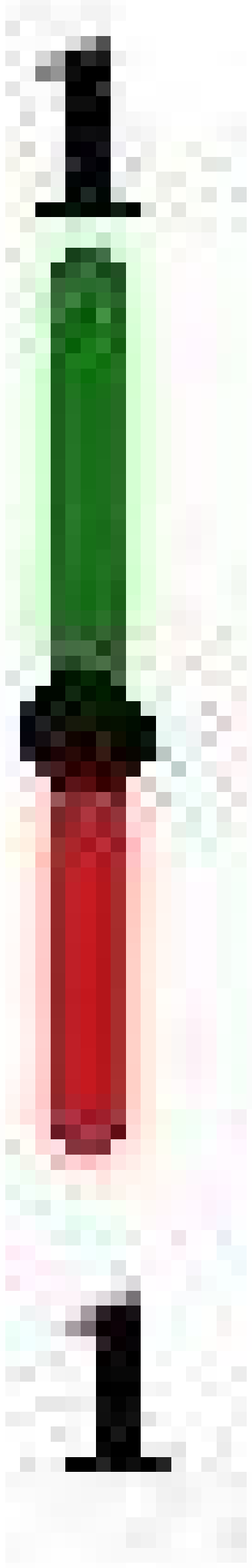}} & \hspace{1in} \phi_{c \mapsto o} : V_c \mapsto V_o \end{array} $$

\item[g6]{\large\bf open to closed} $$ \begin{array}{lll} Type : {1_o}(1_i) & \hspace{1in} \put(0,-20){\includegraphics[height=.8in]{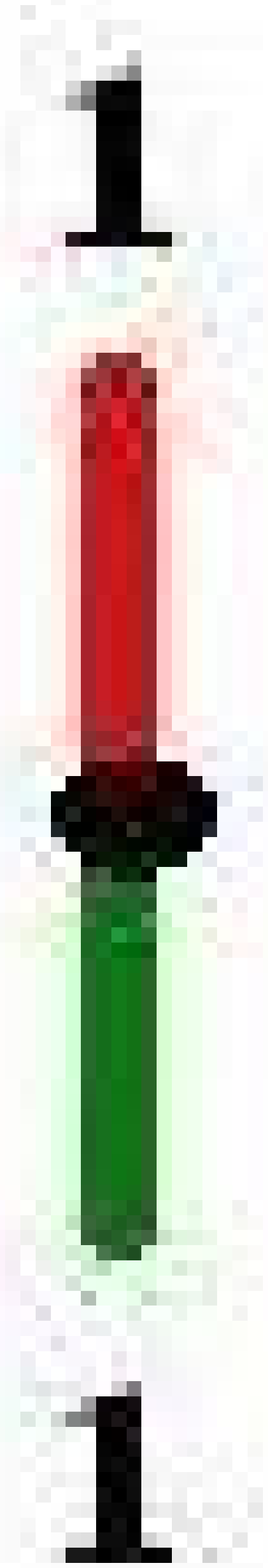}} & \hspace{1in} \phi_{o \mapsto c} : V_o \mapsto V_c \end{array} $$

\item[g7]{\large\bf closed comultiplication} $$ \begin{array}{lll} Type : {1_i,1_o,2_o}& \hspace{1in} \put(0,-20){\includegraphics[height=.8in]{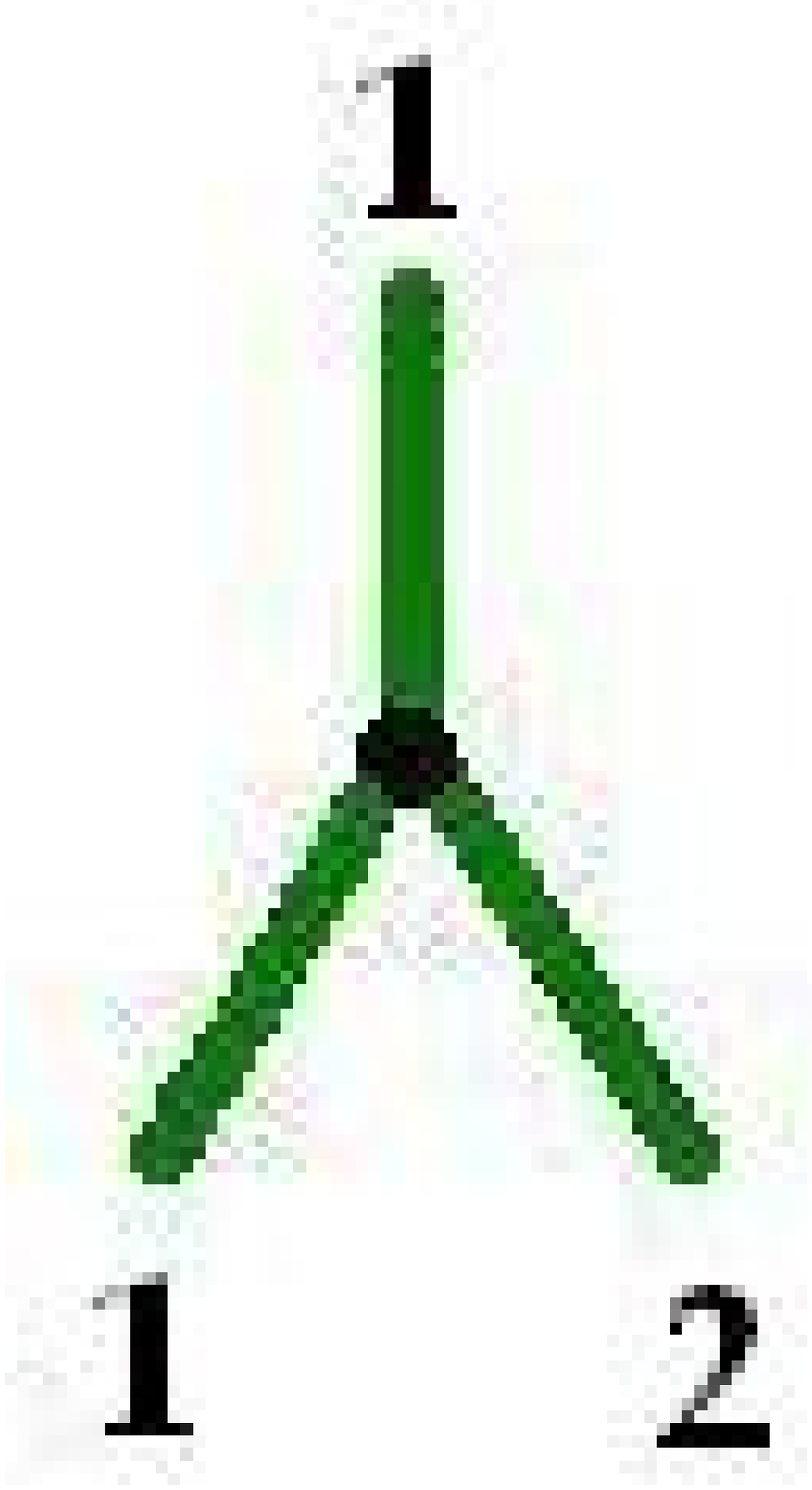}} & \hspace{1in} \bigtriangledown_c : V_c \mapsto V_c \otimes V_c \end{array} $$

\end{description}
\begin{proof}
In the proof of the main theorem of this section it will be defined what it means for a tree to be in {\em normal form}.  It is easy to see that all path components can be given by some tree in normal form.         \end{proof}

Next, we see that the degree 0 operations satisfy the following relations:

\begin{description}

\item[r1] Closed multiplication is associative and commutative with $e_c$ as a unit.
\item[r2] Open multiplication is associative with $e_o$ as a unit.
\item[r3] $V_o$ is an algebra over $V_c$ via $\phi_{c \mapsto o}$.  That is $\phi_{c \mapsto o}$ is an algebra homomorphism into the (graded) center of $V_o$ with $\phi_{c \mapsto o}(e_c)=e_o$ :

\begin{center}
\includegraphics[width=3in]{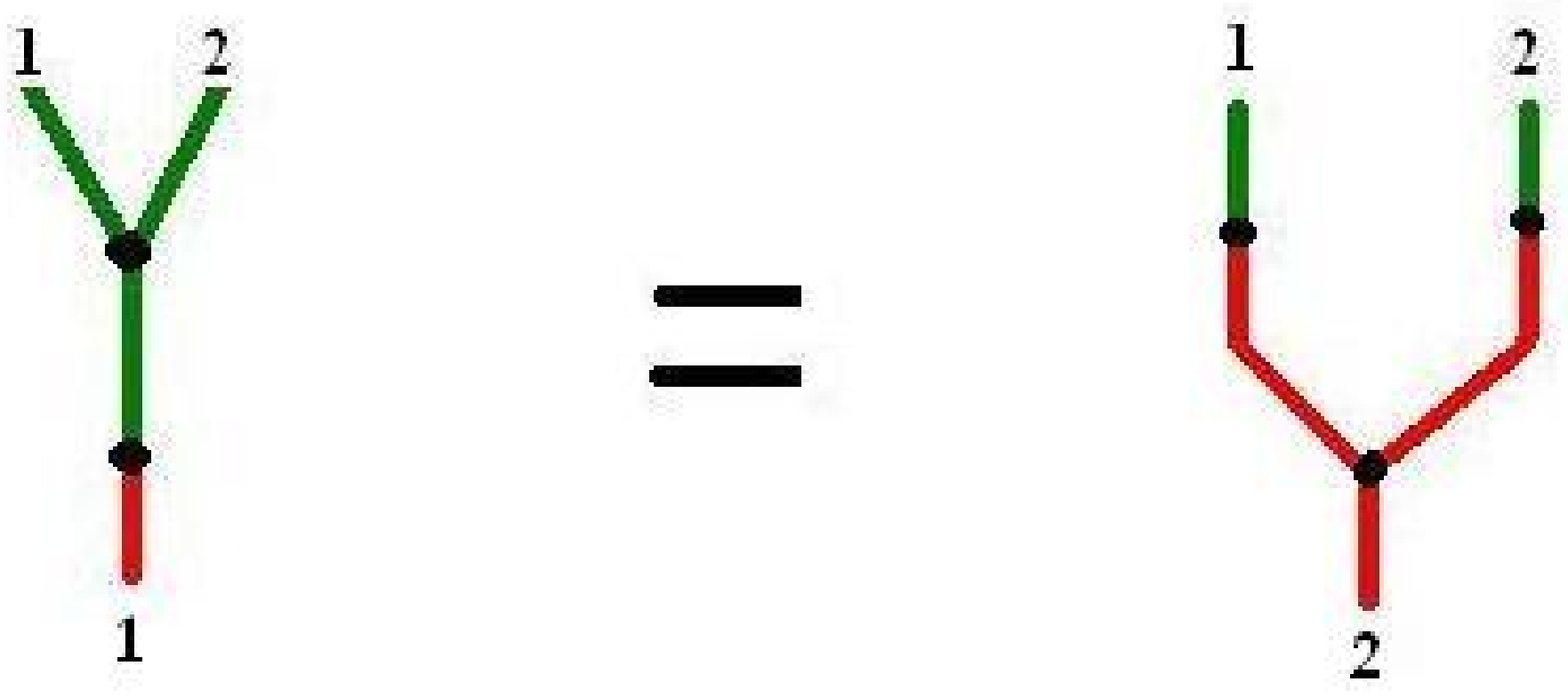}
\end{center}

\end{description}

\begin{center}
$\begin{array}{lr}
\includegraphics[width=3in]{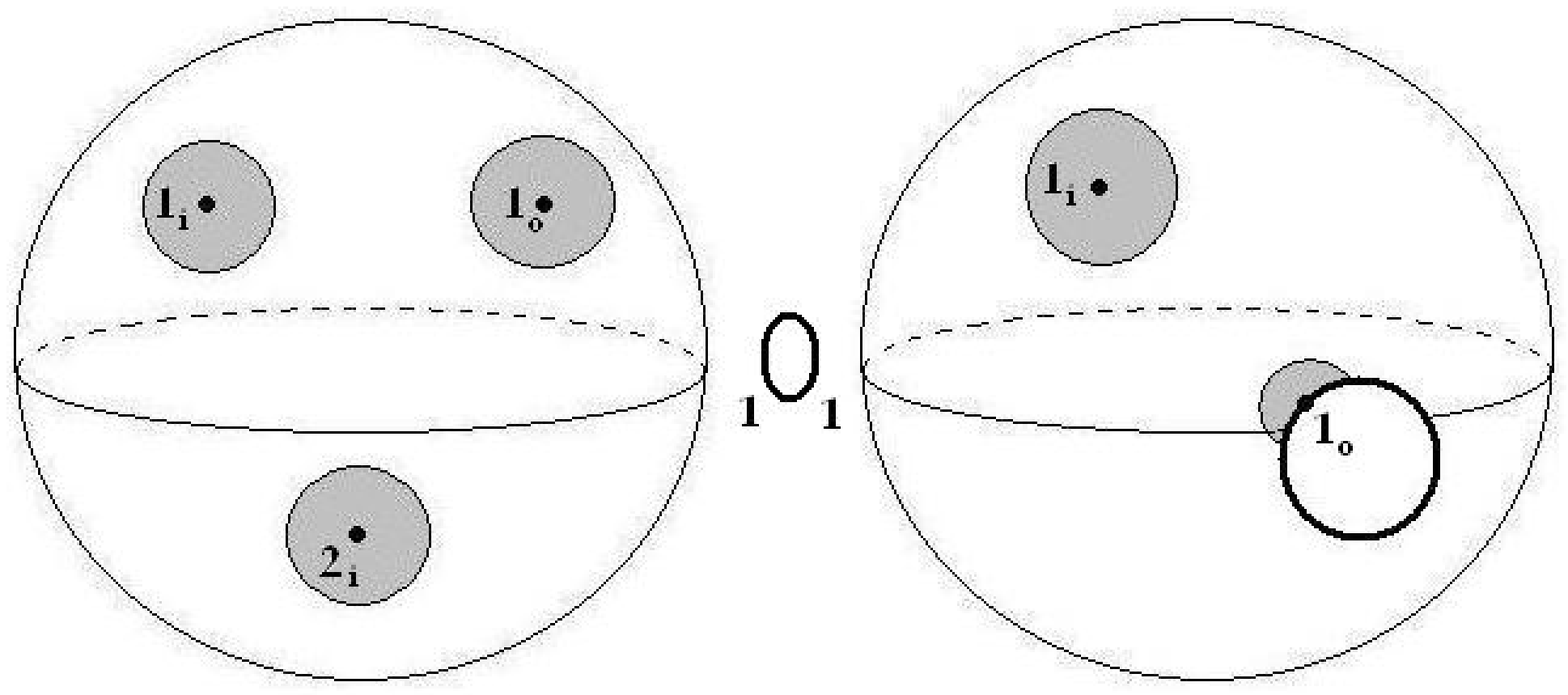} & \includegraphics[width=3in]{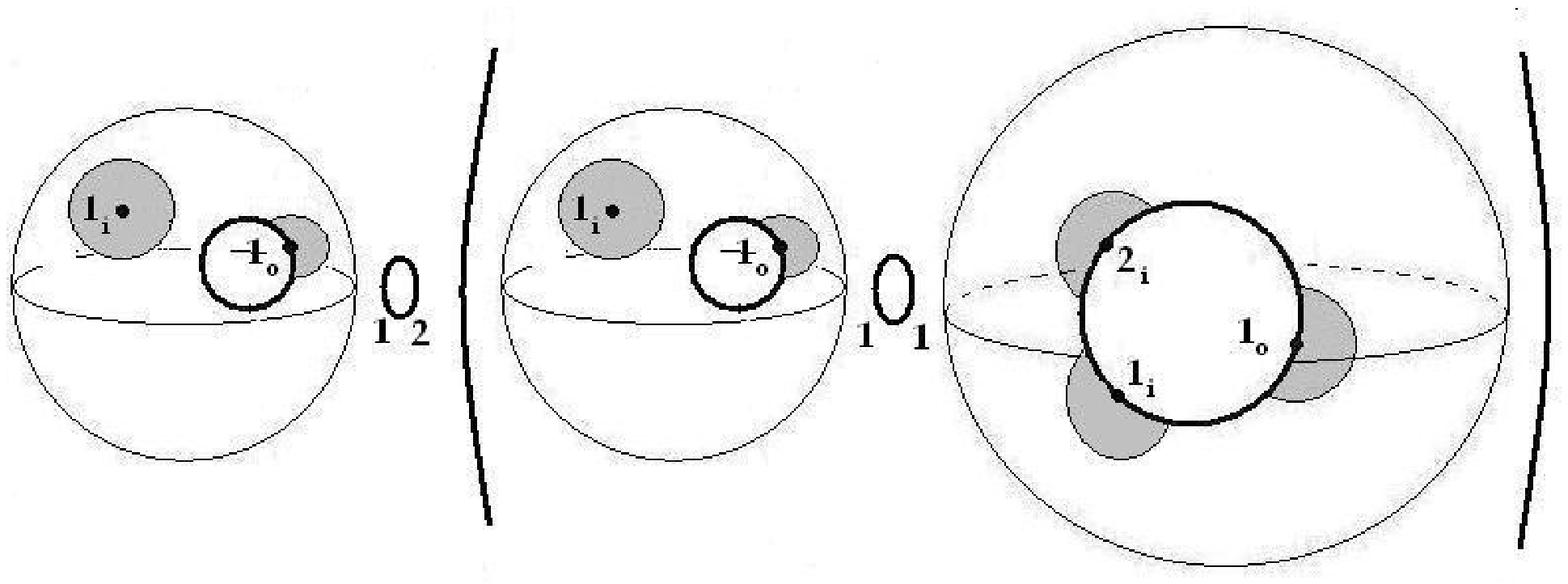} \\
\includegraphics[width=3in]{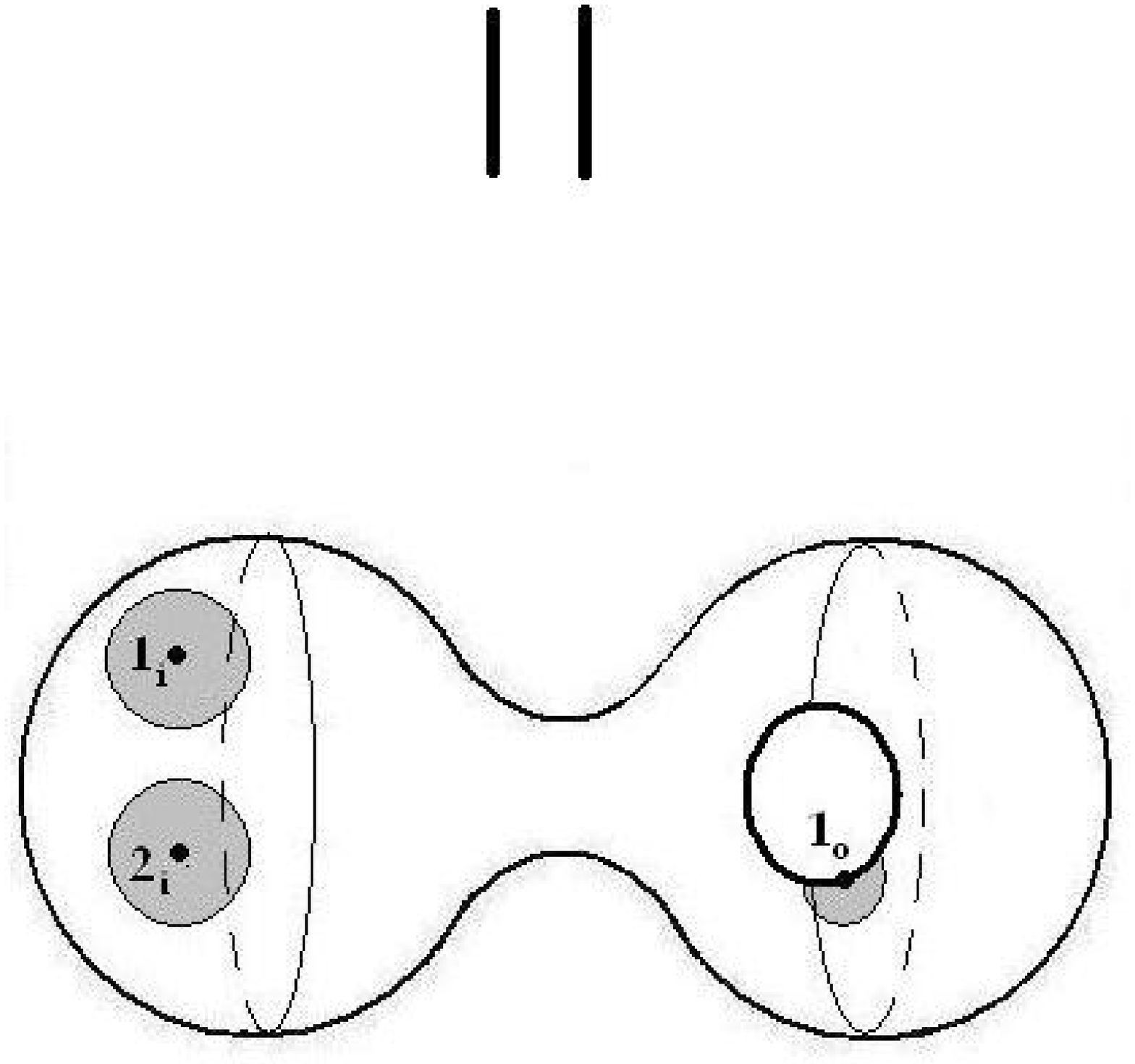} & \includegraphics[width=3in]{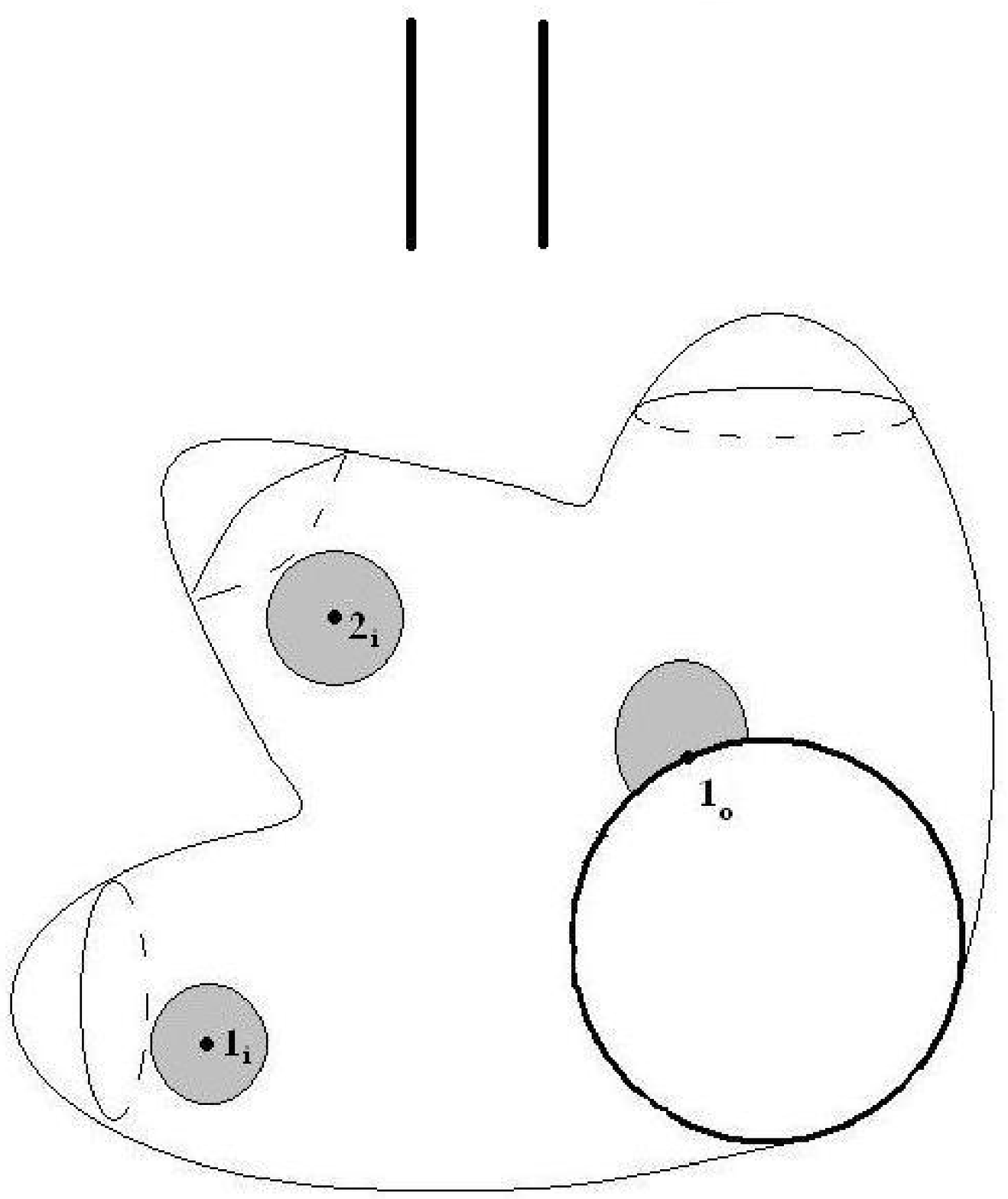}
\end{array}$
\end{center}

\begin{center}
{\huge $Type : \{1_i,2_i\},(1_o)$}
\end{center}
{\large and the picture giving that it is into the center:}

\begin{center}
\includegraphics[width=3in]{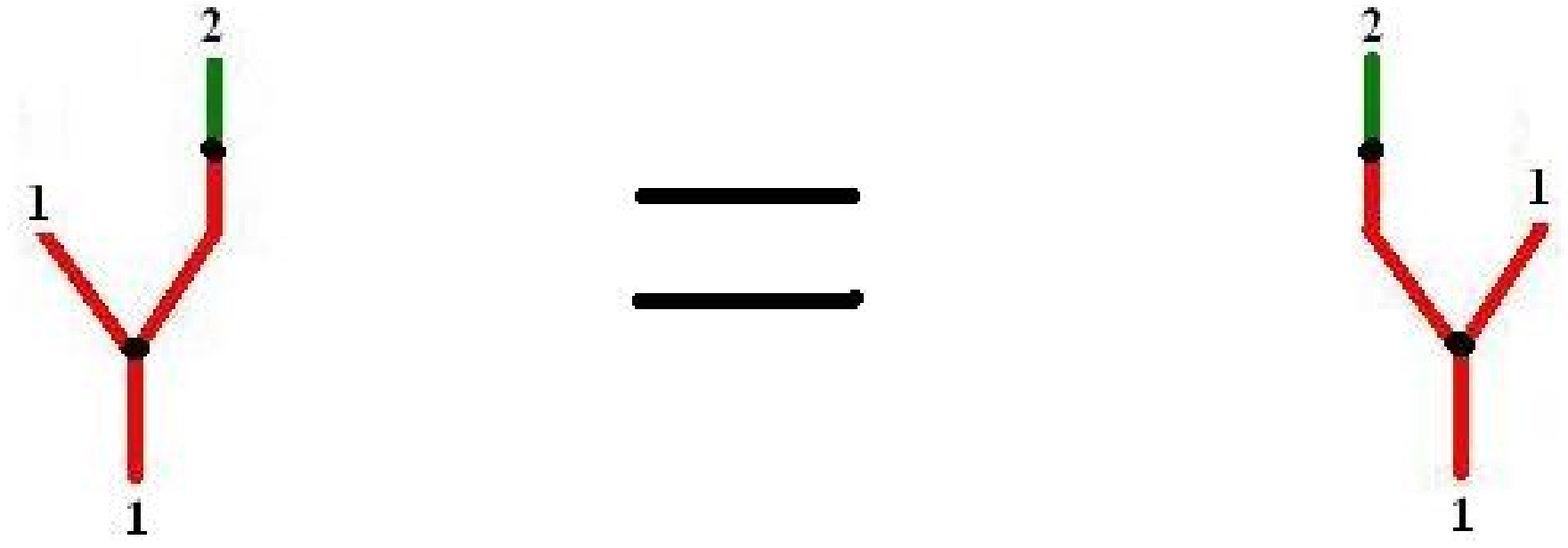}
\end{center}

\begin{center}
\includegraphics[width=6in]{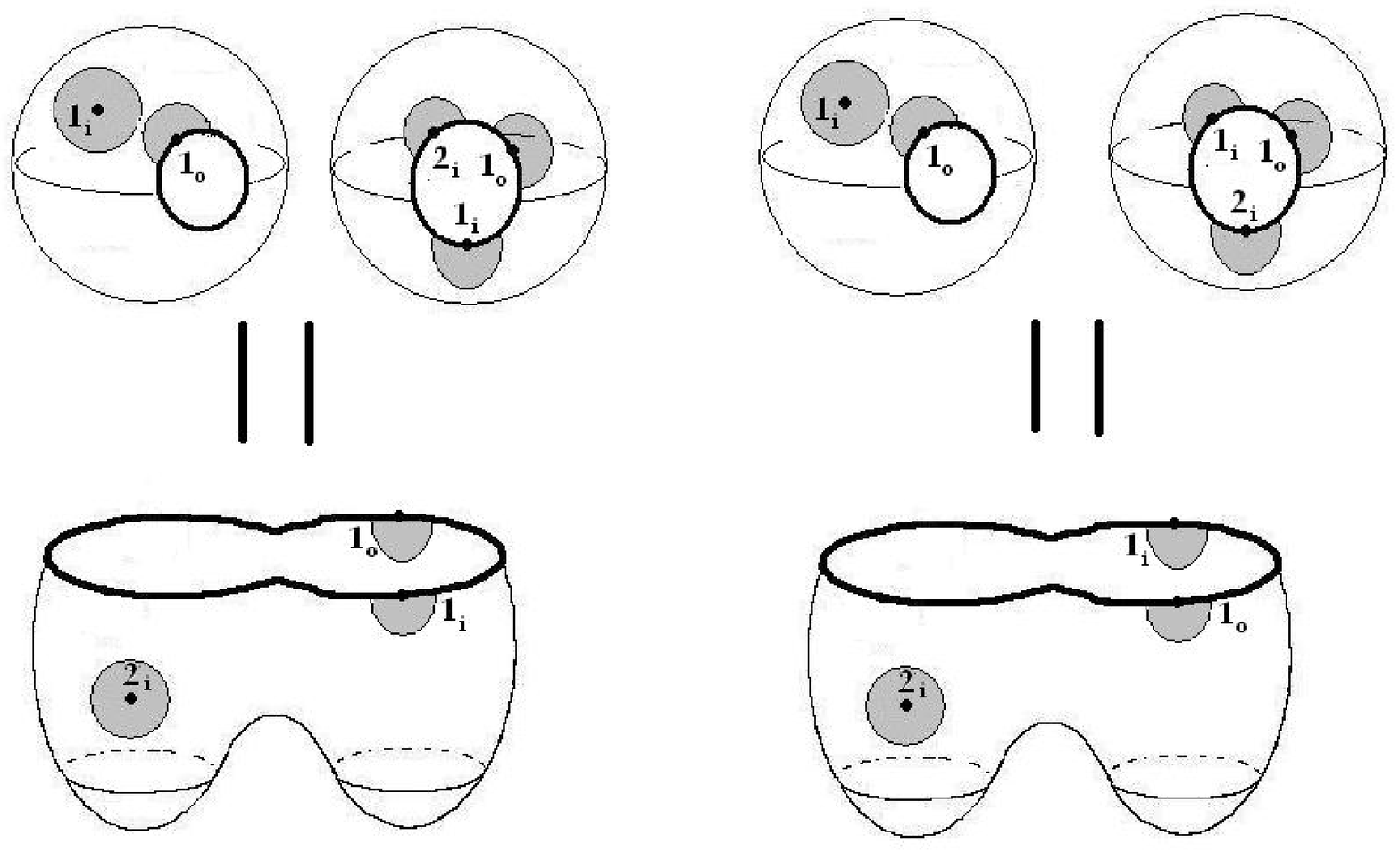}
\end{center}

\begin{center}
{\huge $Type : \{2_i\},(1_i,1_o)$}
\end{center}

\begin{description}
\item[r4] $\phi_{o \mapsto c}(\phi_{c \mapsto o}(a)b)=a \phi_{o \mapsto c}(b)$   (left side is open mult. and right side is closed):

\begin{center}
\includegraphics[width=3in]{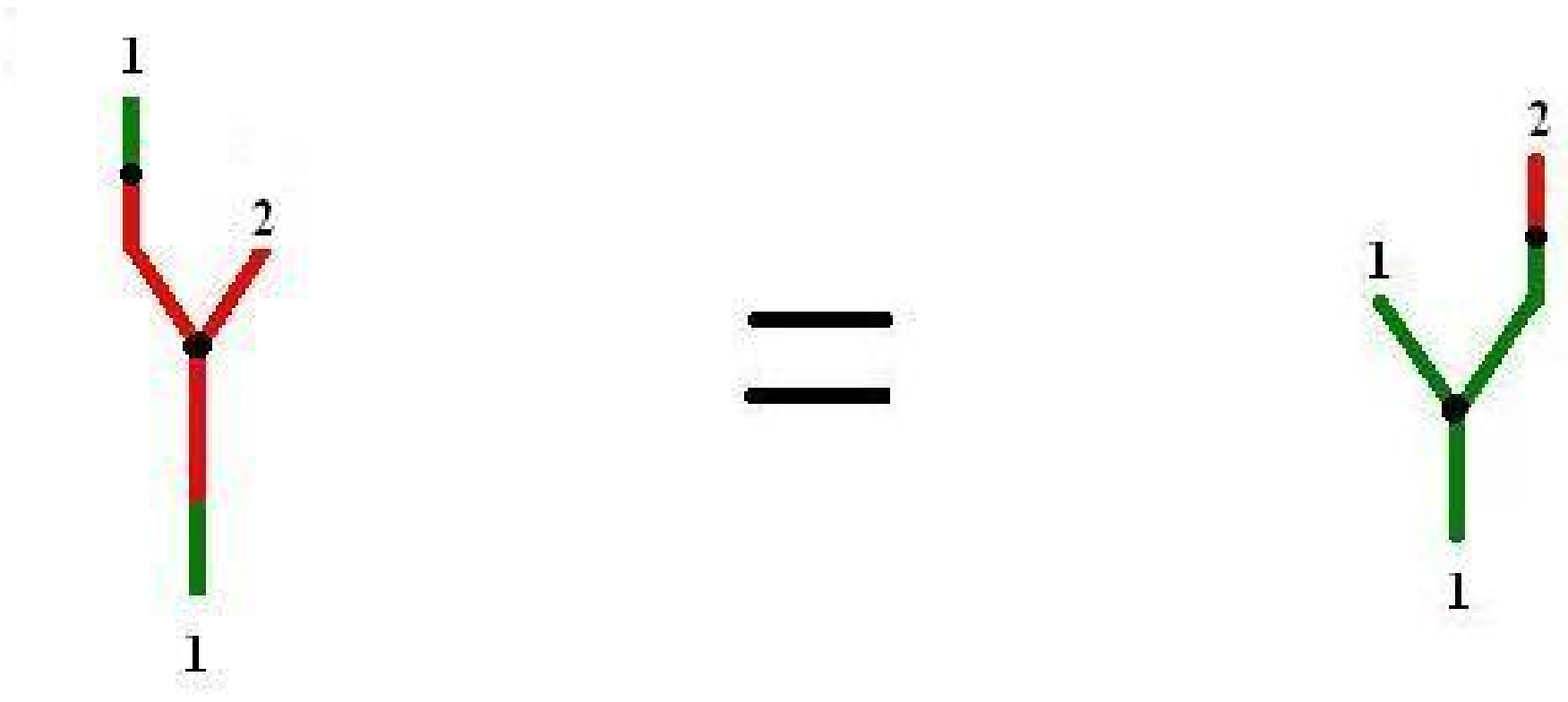}
\end{center}

\begin{center}
{\LARGE $Type : \{1_i,1_o\},(2_i)$}
\end{center}

\item[r5] $\phi_{o \mapsto c}(ab)= (-1)^{|a||b|} \phi_{o \mapsto c}(ba)$

\begin{center}
\includegraphics[width=3in]{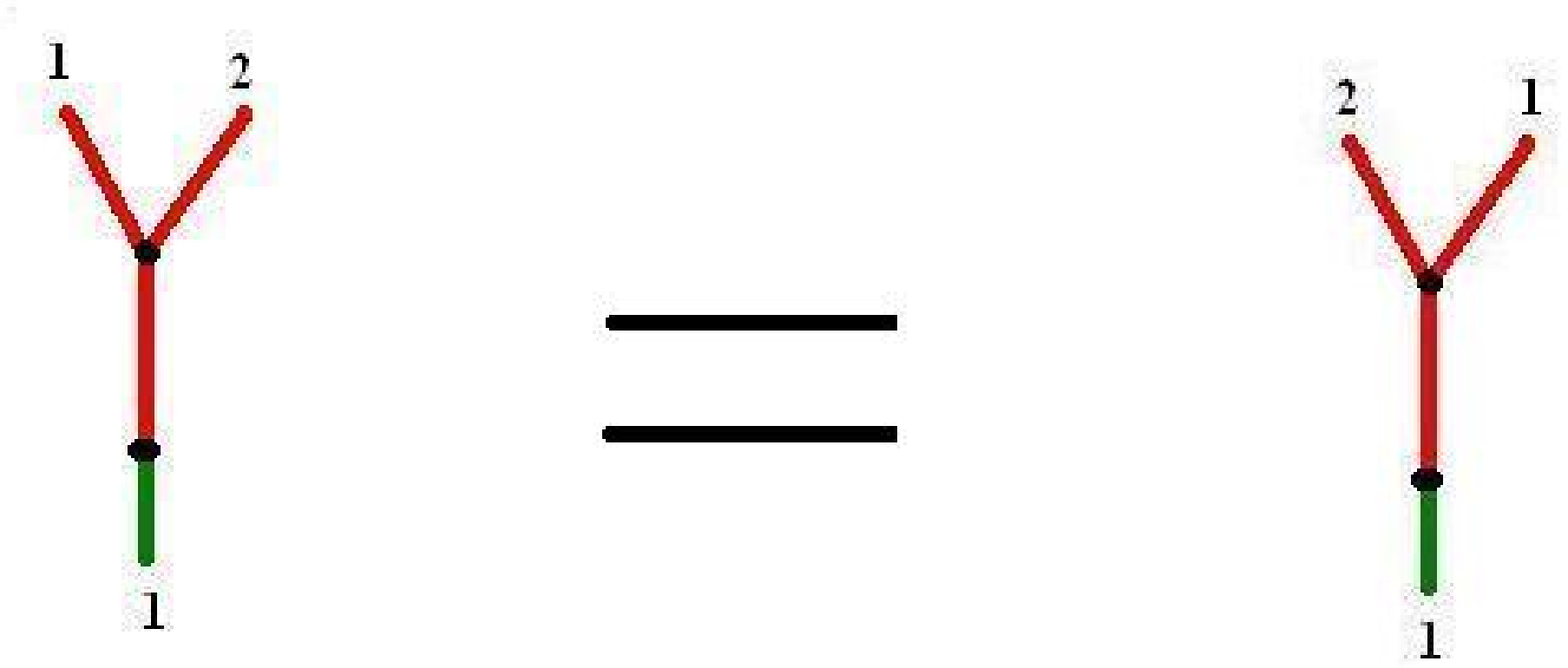}
\end{center}

\begin{center}
{\LARGE $Type : \{1_o\},(1_i,2_i)$}
\end{center}

\item[r6] We have the dual relations except there is no counit (actually, in the course of the proof I found that this relation is not needed, it can be deduced from the others).

\item[r7] $a \bullet \bigtriangledown_c (b)=\bigtriangledown_c (a) \bullet b=\bigtriangledown_c(ab) \mbox{  where  } a \bullet (b_1 \otimes b_2) = (ab_1) \otimes b_2 \mbox{ and } (a_1 \otimes a_2) \bullet b=a_1 \otimes (a_2b)$

\begin{center}
\includegraphics[width=3in]{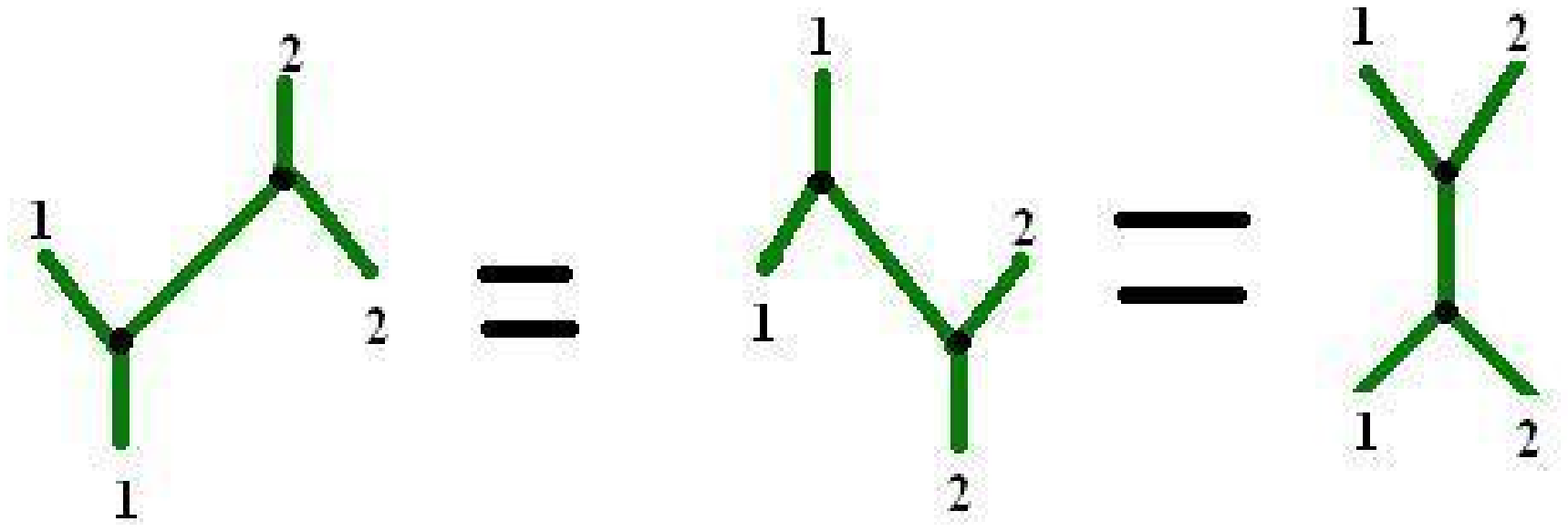}
\end{center}

\begin{center}
{\LARGE $Type : \{1_i,2_i,1_o,2_o\}$}
\end{center}

\item [r8] The same relation holds for the open case as well.
\end{description}

\vspace{.2in}

In the closed case, there is only one path component with 2 closed inputs and 2 closed outputs.  We can derive that any two green trees with two inputs and two outputs are equivalent using relations 7 and 1.  In the open case, there are many path components with one boundary component, 2 open inputs, and 2 open outputs.  Relations 2 and 8 are not enough to show that any two red trees with 2 inputs and 2 outputs going to the same path component are equivalent.  Thus we add:

\begin{description}
\item[r9] $b \bullet_R \bigtriangledown_o (a)=(\bigtriangledown_o (b) \cdot (1,2)) \bullet_L a=\bigtriangledown_o(a)\bullet_L b \mbox{  where  } b \bullet_R (a_1 \otimes a_2) = (a_1b) \otimes a_2,  \,\,  b_1 \otimes b_2 \cdot (1,2) = b_2 \otimes b_1 \mbox{ and } (a_1 \otimes a_2) \bullet_L b=a_1 \otimes (ba_2)$

\begin{center}
\includegraphics[width=3in]{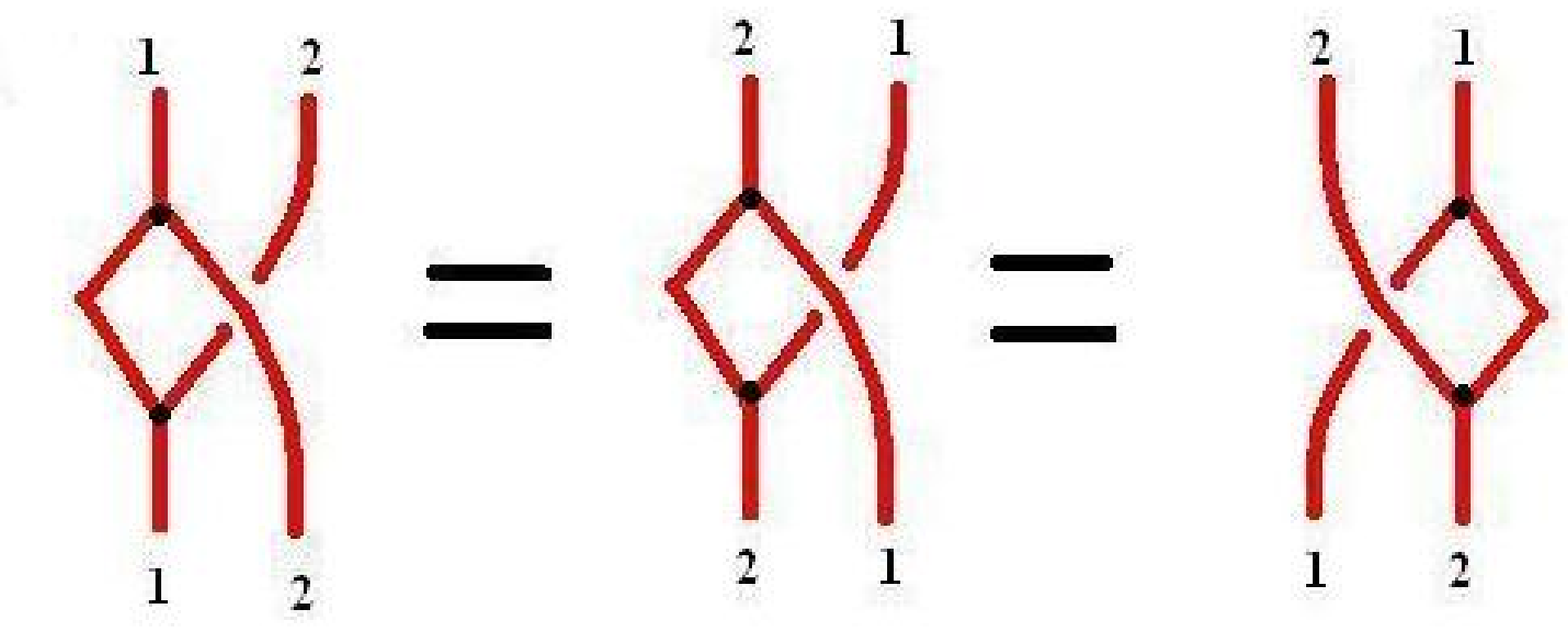}
\end{center}

\begin{center}
{\LARGE $Type : (1_i,2_o,2_i,1_o)$}
\end{center}
\item[r10] $(\bigtriangledown_c (a) \bullet \phi_{c \mapsto o} ) \bullet b = a \bullet (\phi_{o \mapsto c} \bullet \bigtriangledown_o (b)) \mbox{  where  } (a_1 \otimes a_2)\bullet \phi= \\ a_1 \otimes \phi (a_2) \mbox{ and } \phi \bullet (b_1 \otimes b_2)=\phi(b_1) \otimes b_2)$

\begin{center}
\includegraphics[width=3in]{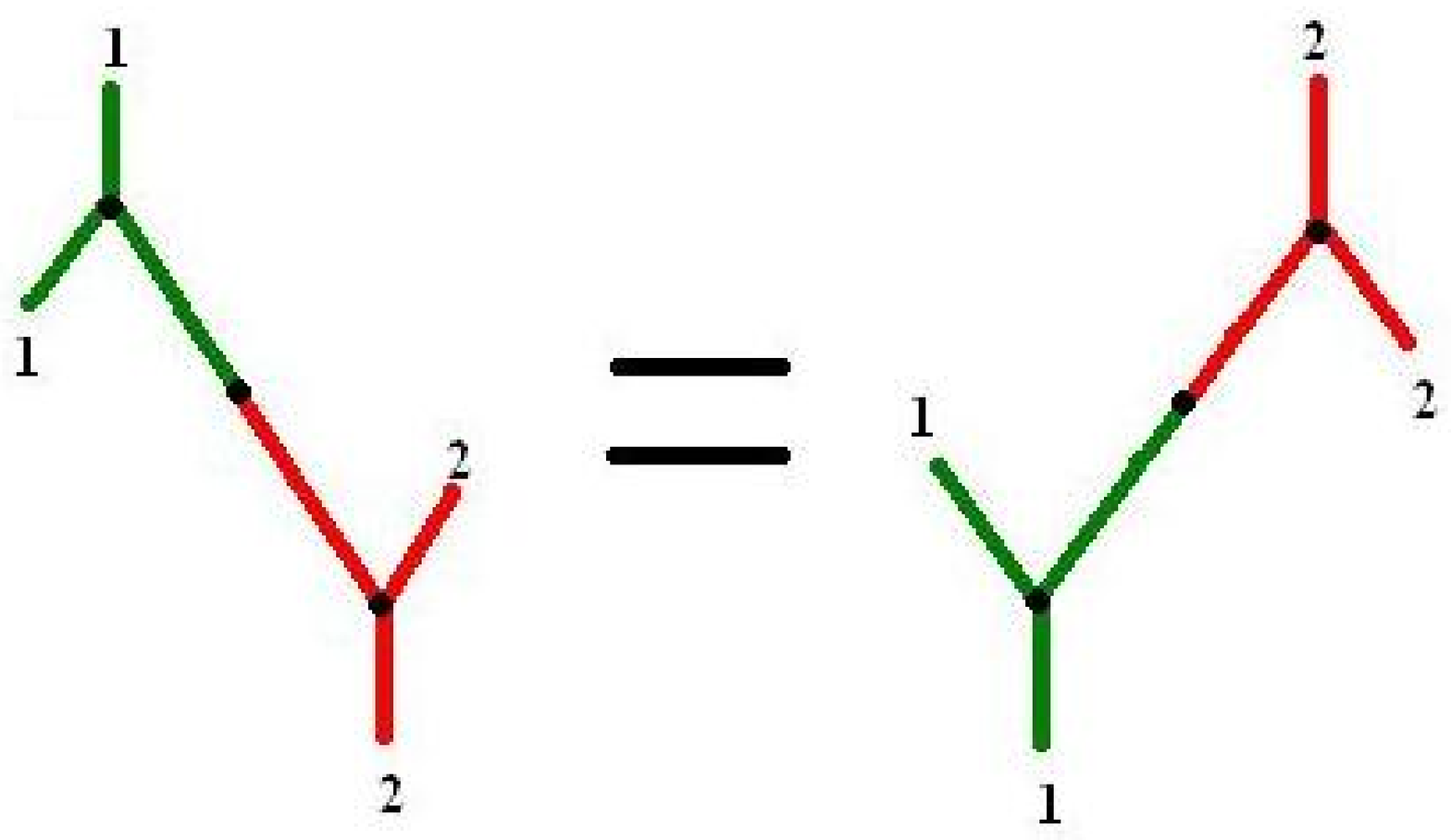}
\end{center}

\begin{center}
{\LARGE $Type : \{1_i,1_o\},(2_i,2_o)$}
\end{center}

\end{description}

\begin{thm}
The 8 degree 0 generators g1,...,g8 and the 10 relations r1,...,r10 completely describe $H_0(OC)$ or, equivalently, an algebra $(V_c,V_o)$ over it
\end{thm}

\begin{proof}
Let $F(G)/S$ be the free 2-colored dioperad generated by the degree 0 generators $G=\{g1,..,g8\}$ and relations $S=\{r1,..,r8\}$.  Then right now we have an onto dioperad morphism $F(G)/S \mapsto H_0(OC)$ .  We need to see that it is in fact 1-1.  I.e. we need to check that there is a 1-1 correspondence between path components of $OC$ and equivalence classes of labeled trees in $F(G)/S$.\vspace{.15in}

\noindent\underline{\bf Notation}\\
Since open/closed (co) multiplication is (co) associative, we'll just use a green or red tree like the following to denote (co) multiplication with n inputs (outputs):
\begin{center}
\includegraphics[height=1in]{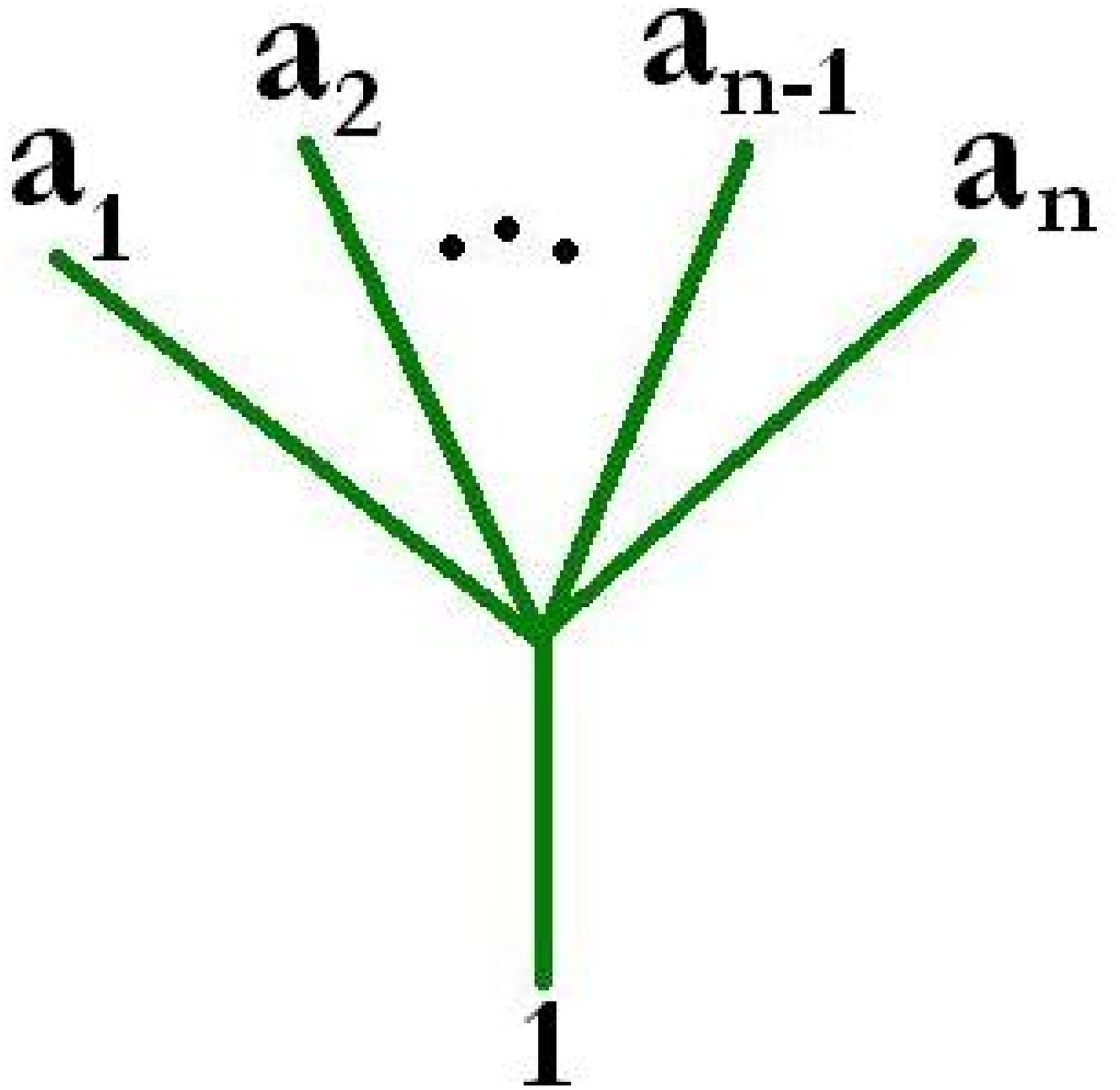}
\end{center}
\begin{df}
We'll say a tree is in {\em normal form} if it is in one of the following three forms:
\end{df}

\noindent \underline{\bf {\Large Form 1}}\\
A completely red tree is in normal form if it looks like:\\
$$\begin{array}{ll}
\put(0,-50){\includegraphics[width=2in]{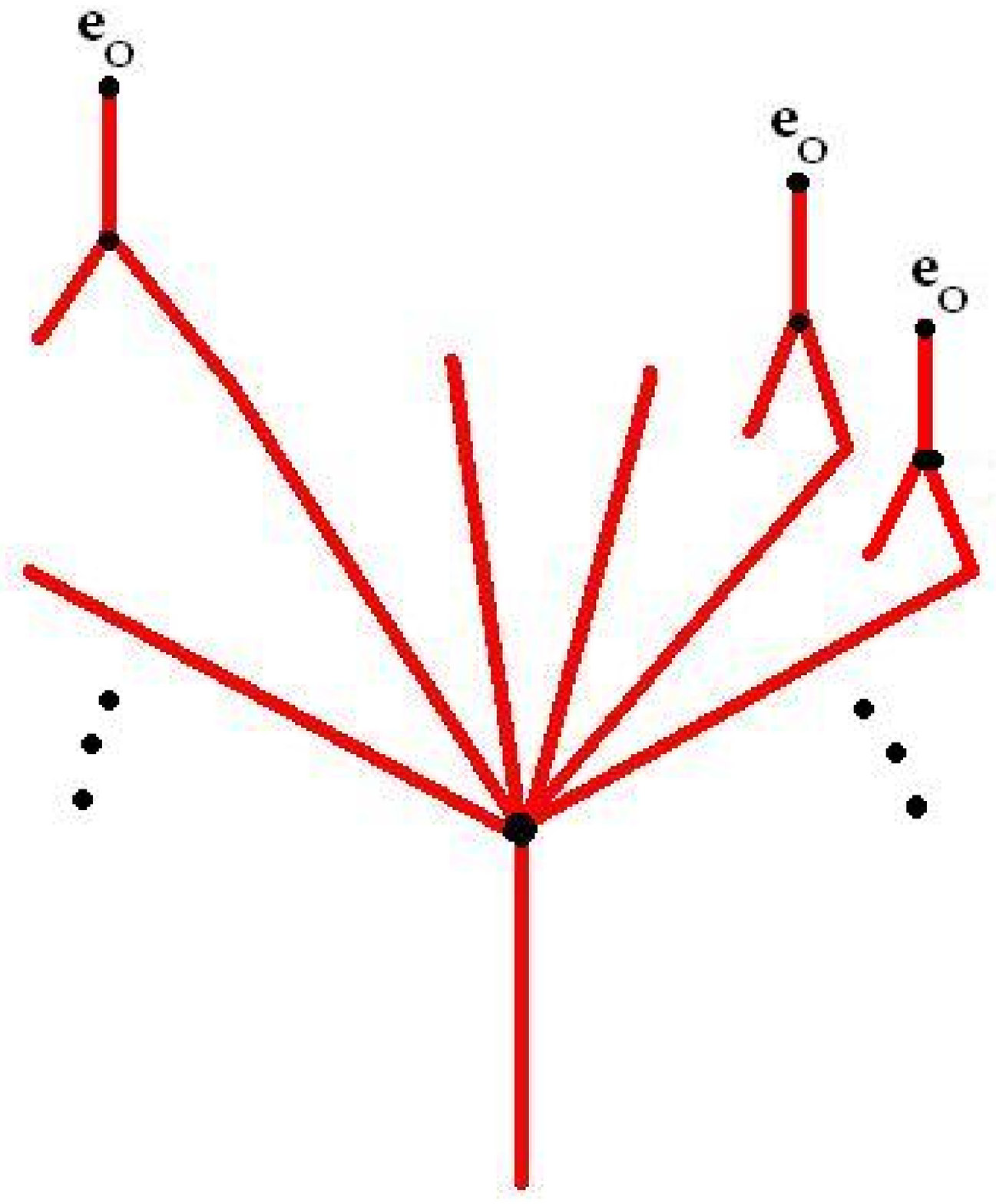}} \hspace{2in} & \begin{array}[b]{l} \mbox{\underline{explanation}:}\\ \mbox{All possible labelings are allowed.} \\ \mbox{The first and last stems can be} \\ \mbox{either output or input stems.}\end{array} \end{array}$$
\newpage

\noindent The following example shows what the image of a tree of this form is in $H_0(OC)$:
\begin{center}
\includegraphics[width=5in]{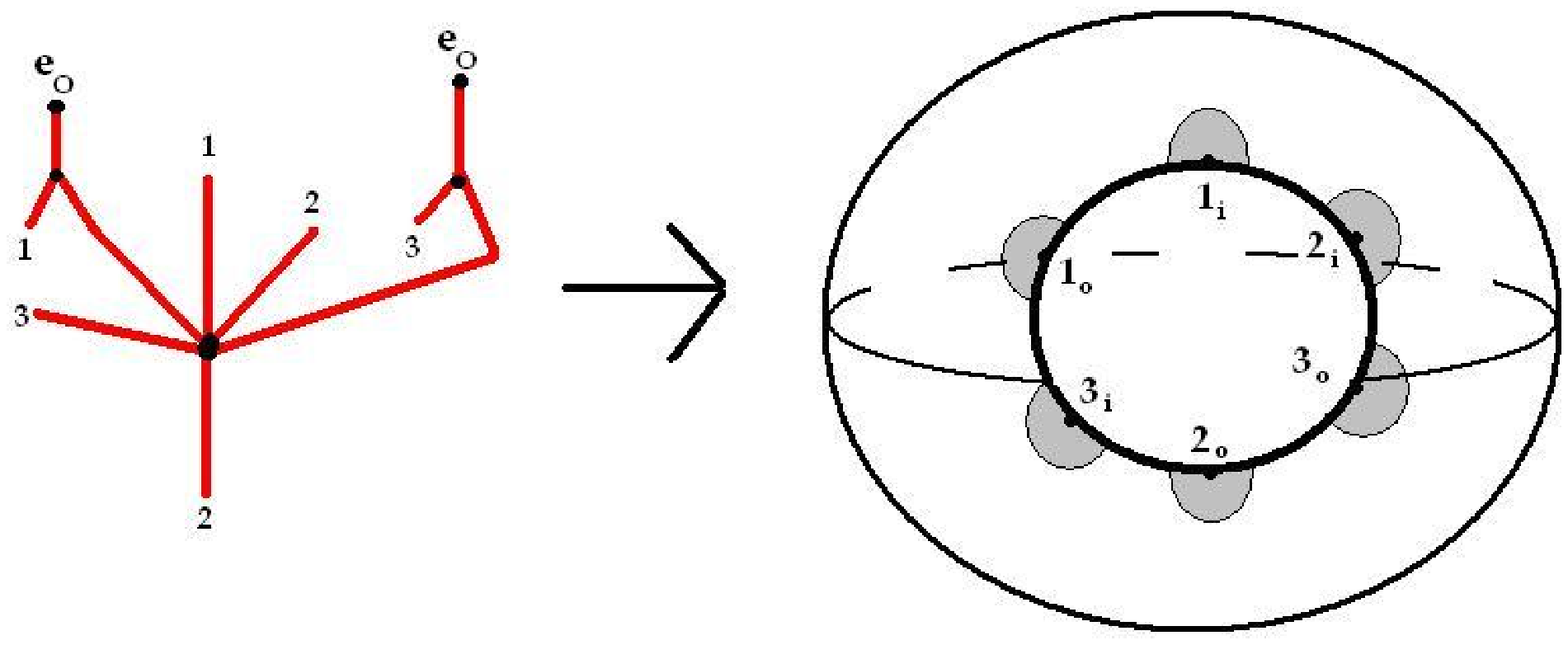}
\end{center}
\vspace{1in}

\noindent \underline{\bf {\Large Form 2}}\\
Trees with at least one closed output are in normal form if they look like:\\
$$\begin{array}{ll}
\put(0,-100){\includegraphics[width=2in]{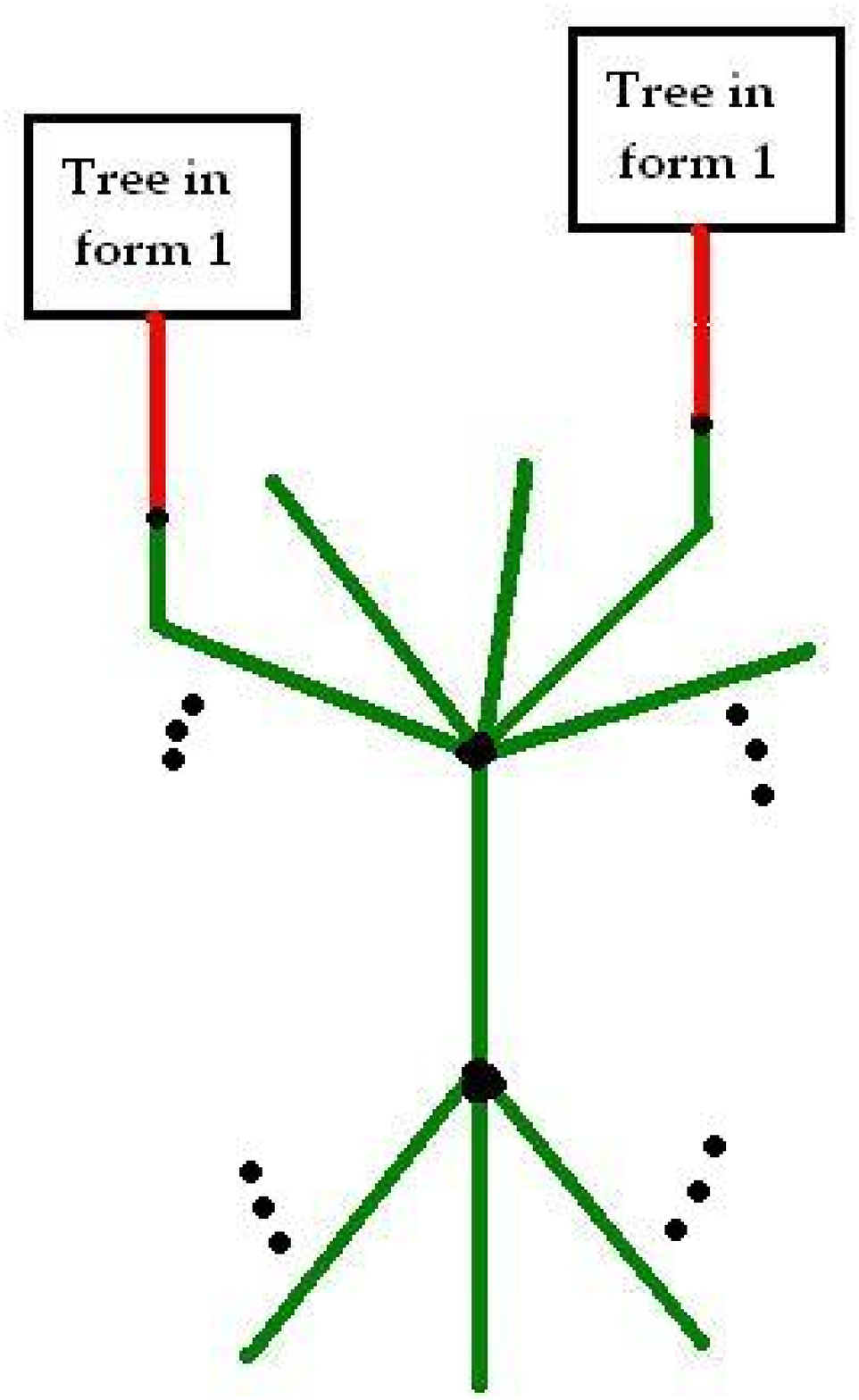}} \hspace{2in} & \begin{array}[b]{l} \mbox{\underline{explanation}:}\\ \mbox{All, none, or some of the green input} \\ \mbox{stems can have a tree in form 1 connected to it via} \\ \phi_{o \mapsto c} \end{array} \end{array}$$
\newpage

\noindent The following example shows what path component a tree of this form corresponds to:\\
\vspace{.3in}

\begin{center}
\includegraphics[width=5in]{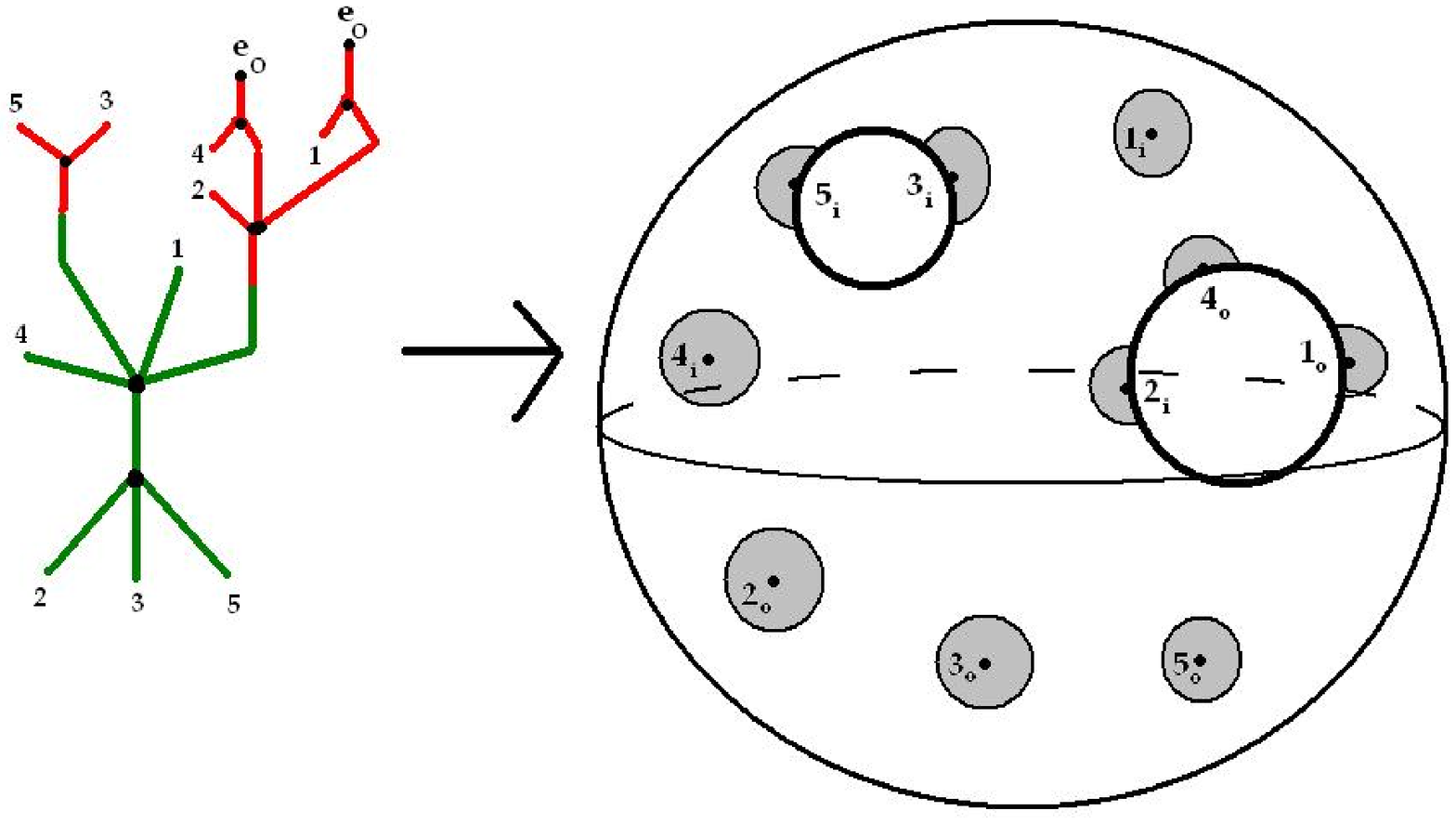}
\end{center}
\vspace{.3in}

\noindent \underline{\bf {\Large Form 3}}\\
A tree which is not equivalent to an entirely red tree but doesn't have a closed output must look like:\\
$$\begin{array}{ll}
\put(0,-100){\includegraphics[width=2in]{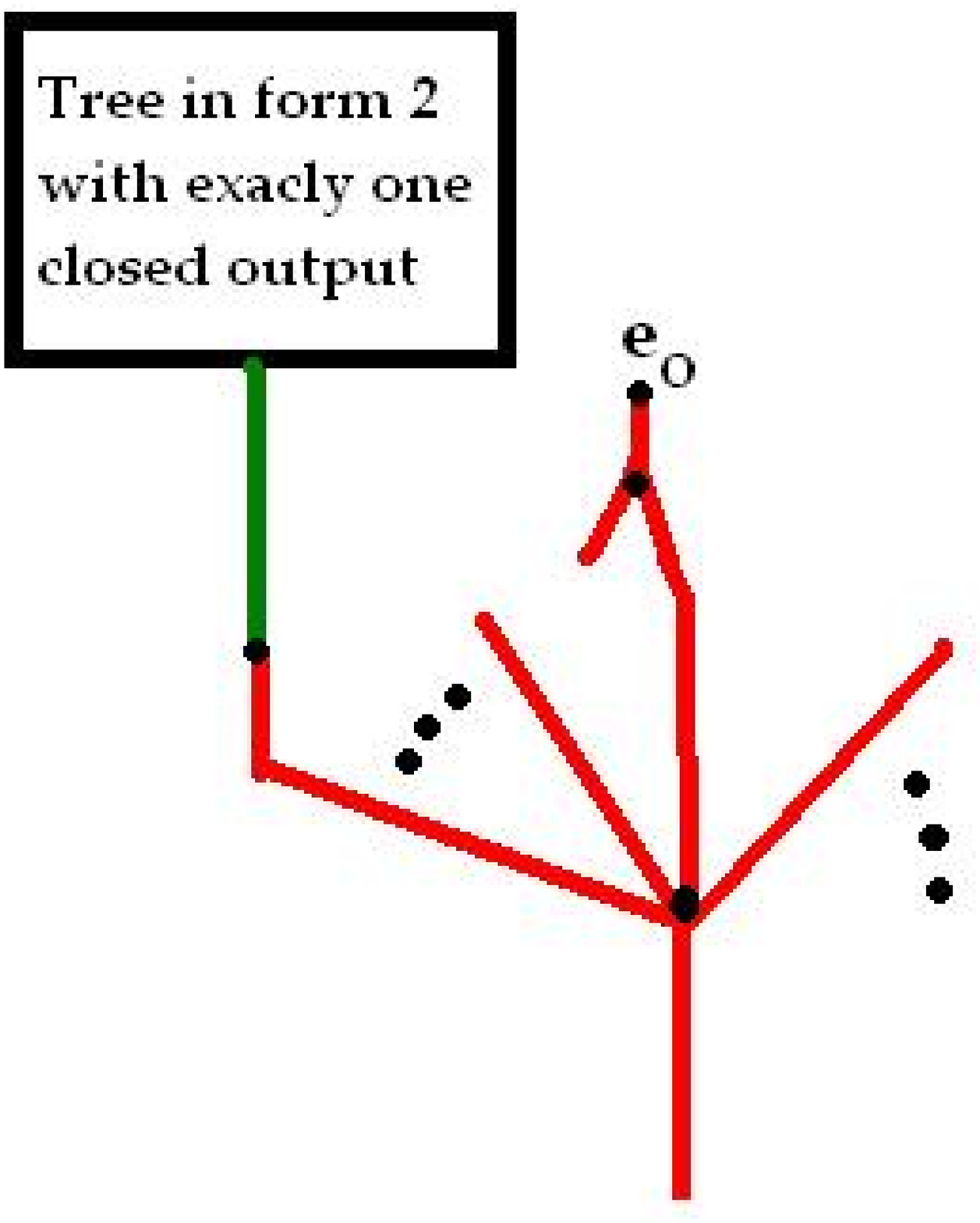}} \hspace{2in} & \begin{array}[b]{l} \mbox{\underline{explanation}:}\\ \mbox{The red tree at the bottom should be a tree of}\\ \mbox{form 1 with its left most stem an input stem.}\\ \mbox{Note: Requiring that the tree not be equivalent to}\\ \mbox{a red tree is the same thing as requiring that the}\\ \mbox{form 2 tree not be }e_c \end{array} \end{array}$$
\newpage
\noindent The following example shows what path component a tree of this form corresponds to:\\

\begin{center}
\includegraphics[width=5in]{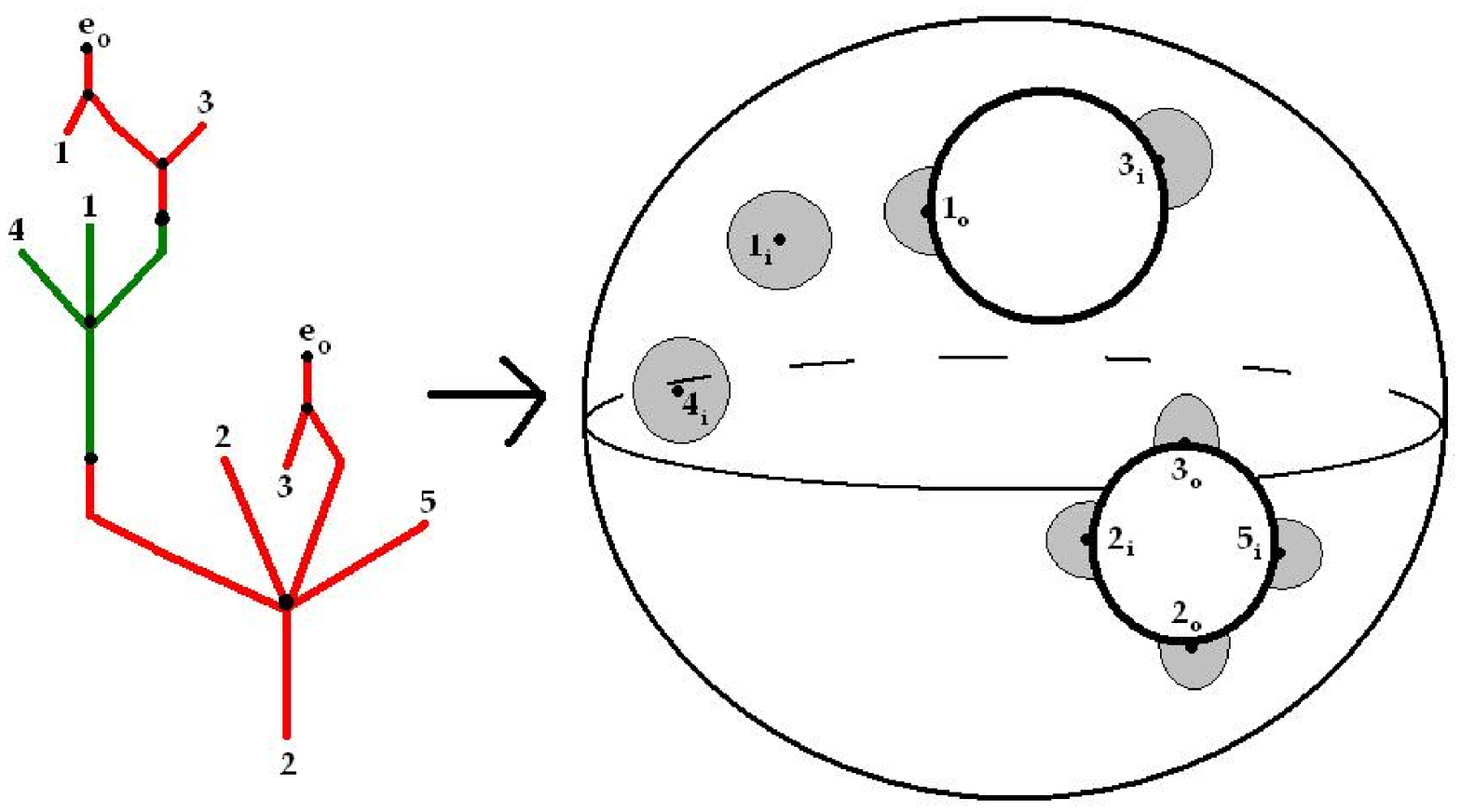}
\end{center}
\vspace{.3in}
\begin{claim}
There is a bijection between the equivalence classes of trees in normal form and the path components of the moduli space.
\end{claim}
\begin{proof}

As noted earlier in Proposition 2.1, it is clear that every path component corresponds to some tree in normal form.\\
Now, two trees in form 1 give the same path component iff they give the same circular permutation of the labels.  To show that two red trees giving the same circular permutation are equivalent, it suffices to show the following equivalence:\\
$$\begin{array}{ll}
\put(0,-75){\includegraphics[width=2.5in]{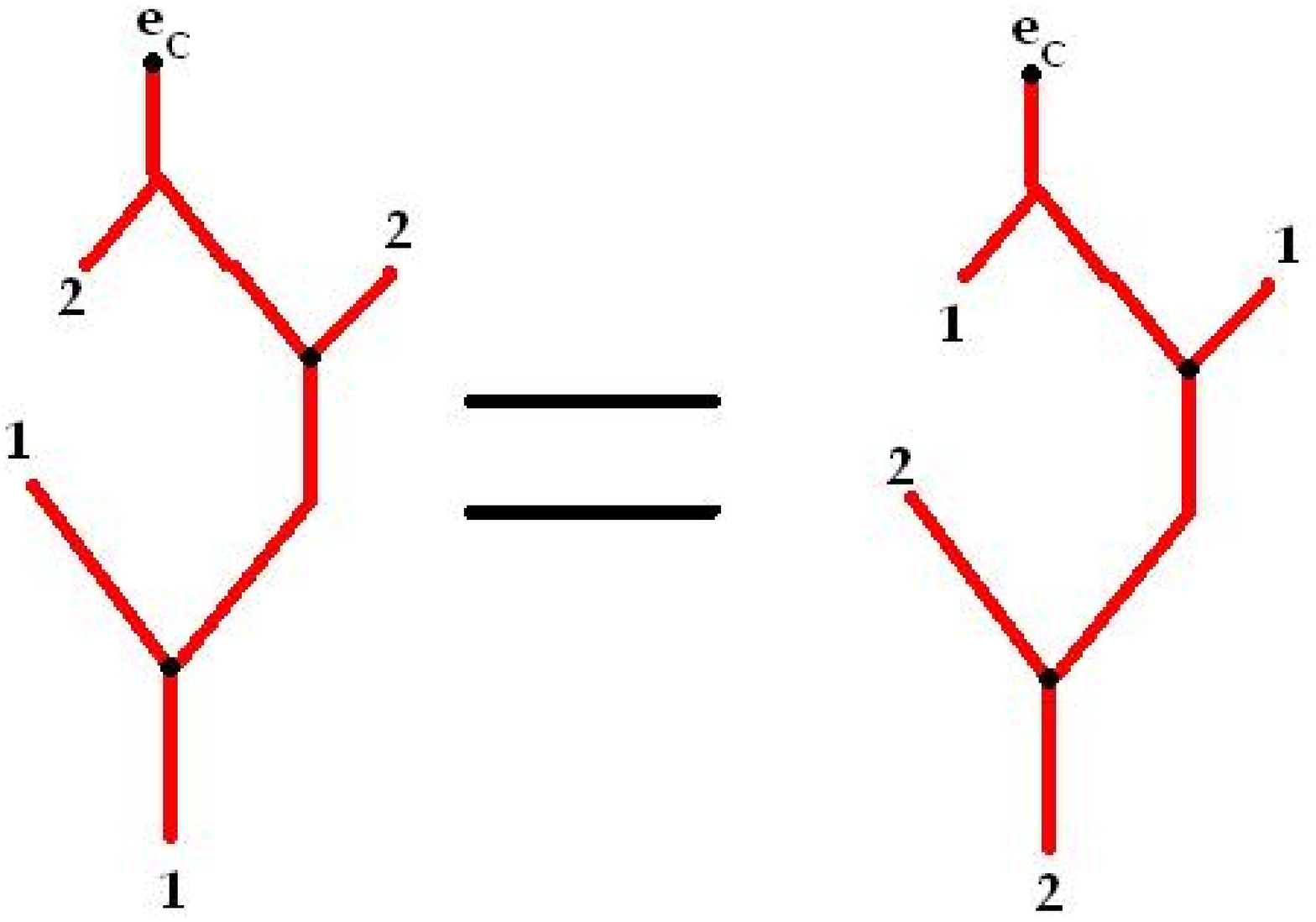}} \hspace{2.5in} & \begin{array}{l}
\mbox{\underline{note} that this covers the cases where}\\ \mbox{the output that we are "rotating" to the}\\ \mbox{root is the first or last stem.  To see this }\\ \mbox{just plug $e_o$ into the first or second input.} \end{array} \end{array}$$

\newpage
The equivalence can be shown through the following sequence of equivalences (the numbers on top of the equal signs say which relation is being used):
\begin{center}
\includegraphics[width=4in]{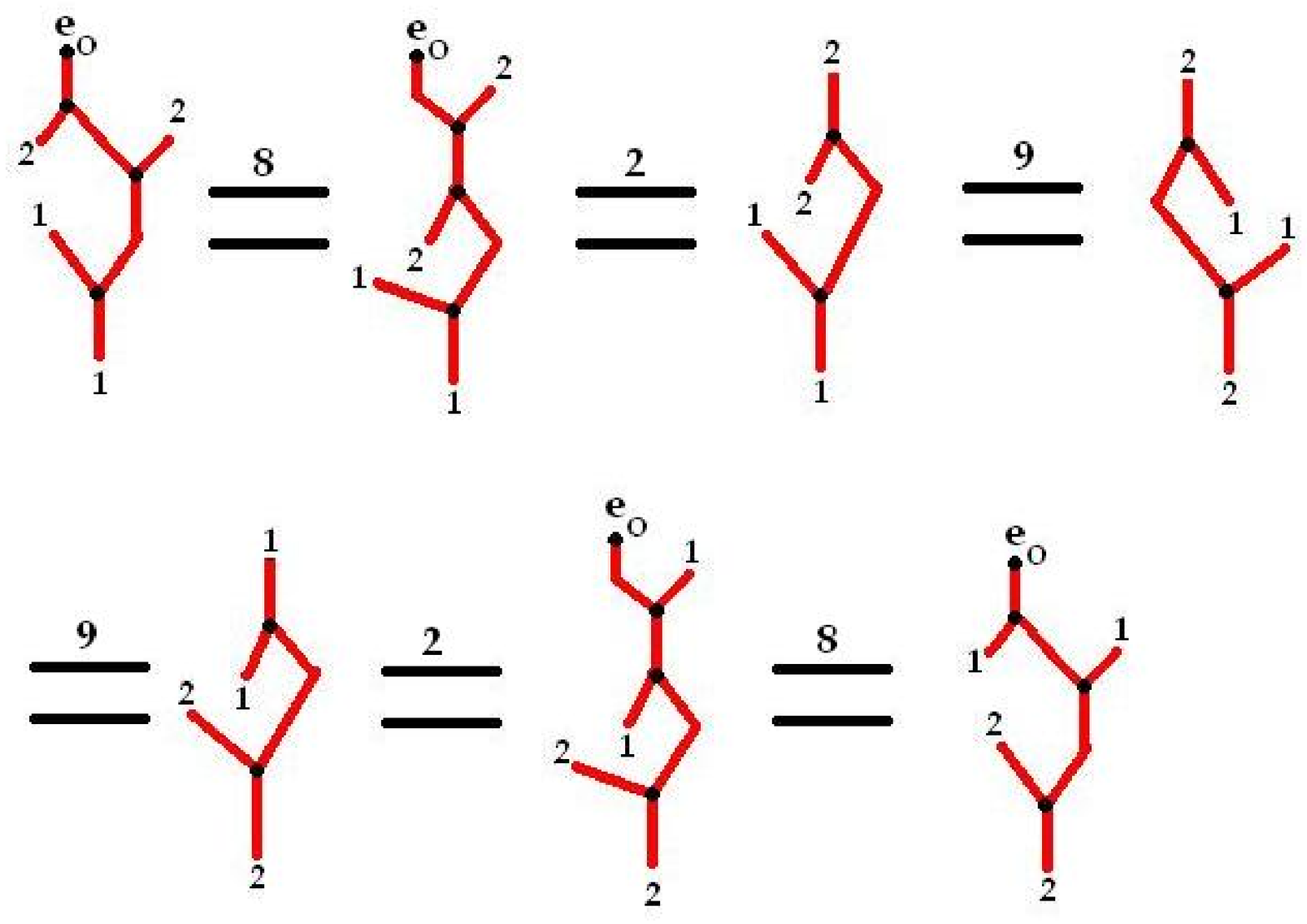}
\end{center}

Now consider trees of form 2.  Given that two trees of the form of figure A are equivalent as long as they have the same number of inputs and outputs, we see that it suffices to check that two trees of the form of figure B giving the same circular permutation of the OPEN labels are equivalent:
$$\begin{array}{cc}
\includegraphics[height=1.5in]{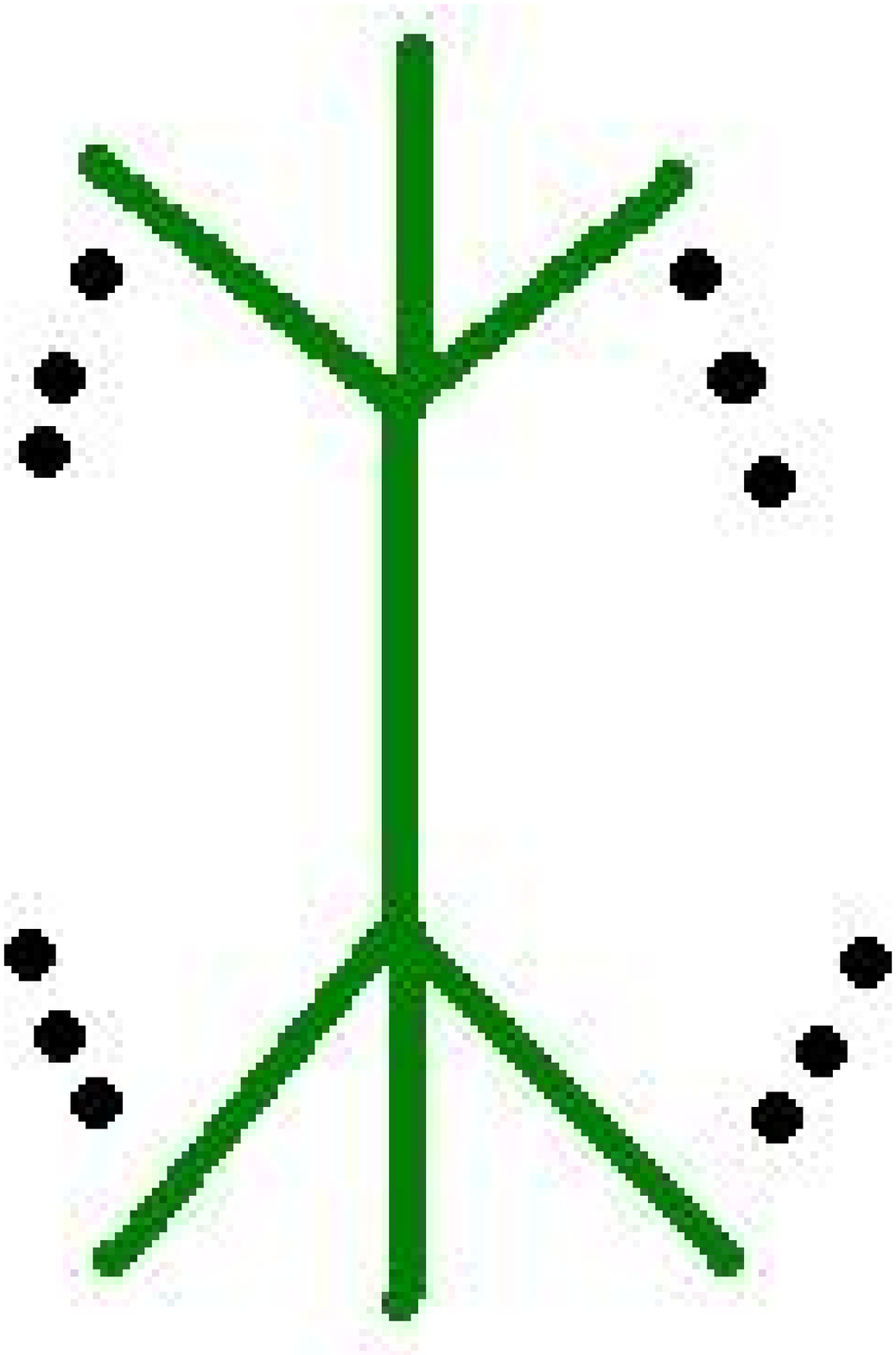} \hspace{1.5in} & \includegraphics[height=1.5in]{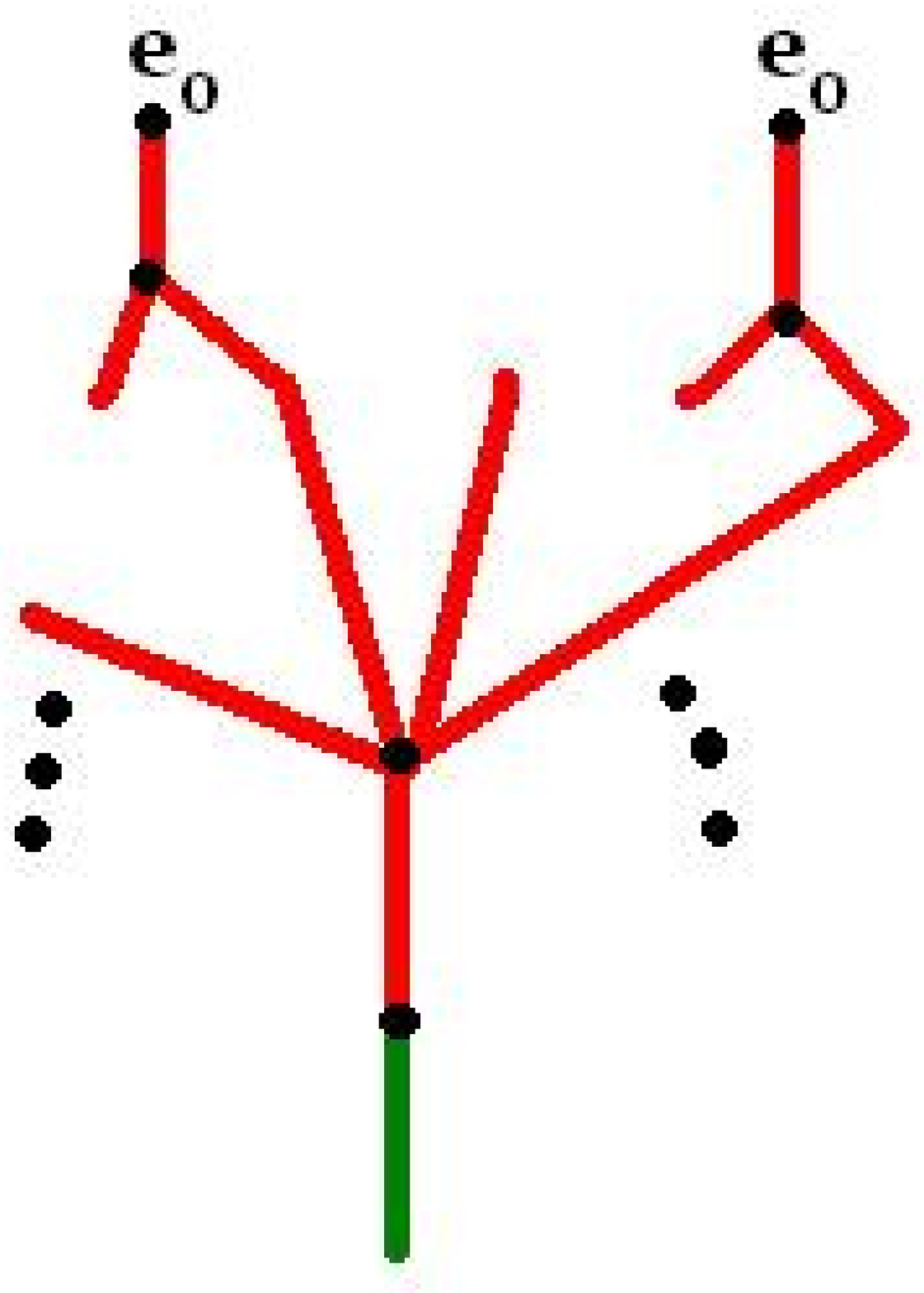} \\ \mbox{{\large Figure A}} \hspace{1.5in} & \mbox{{\large Figure B}}
\end{array}$$

\noindent But this follows directly from relation 5.
\vspace{.2in}

Finally, consider form 3 trees.   What we first need to see is that for each output of a form 3 tree, the tree is equivalent to another form 3 tree which has the output as the main root output of the bottom red tree.  If the output is already on the bottom red tree, then use the form 1 result to make this output the root output, then use relation 3b to move the form 2 stem back onto the left most stem of the bottom red tree.

Otherwise, the output belongs to one of the top red trees.  Using relation 5, we can assume this output stem is the right most stem of this red tree.  Then the following equivalence suffices:
\newpage

\begin{center}
\includegraphics[width=4.5in]{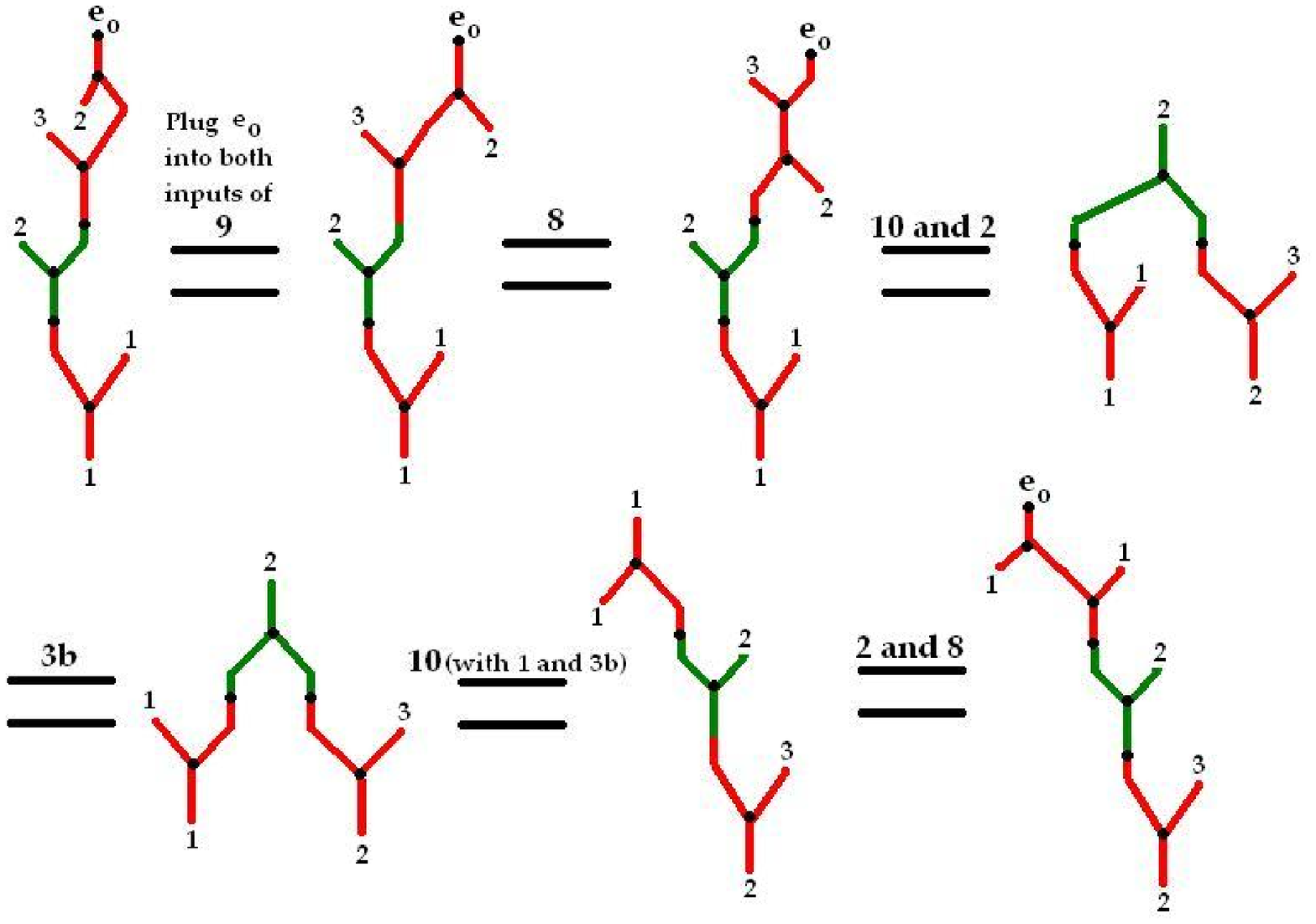}
\end{center}

Now take two trees of form 3 which correspond to the same path component.  By what was just shown, we can assume that their main root outputs are labeled the same.  Then it must be that their bottom red trees are exactly the same since otherwise the trees would go to different path components.  It must also be that their trees of form 2 correspond to the same path component so that they are equivalent by the last case.  Thus the two form 3 trees must be equivalent.

This concludes the proof of claim 2.5
\end{proof}

So now, to complete Theorem 2.3, all that needs to be shown is that any tree in $F(G)/S$ is equivalent to some tree in normal form.  To see this, first note that every generator is equivalent to a tree in normal form:\\

\begin{center}
\includegraphics[width=4.5in]{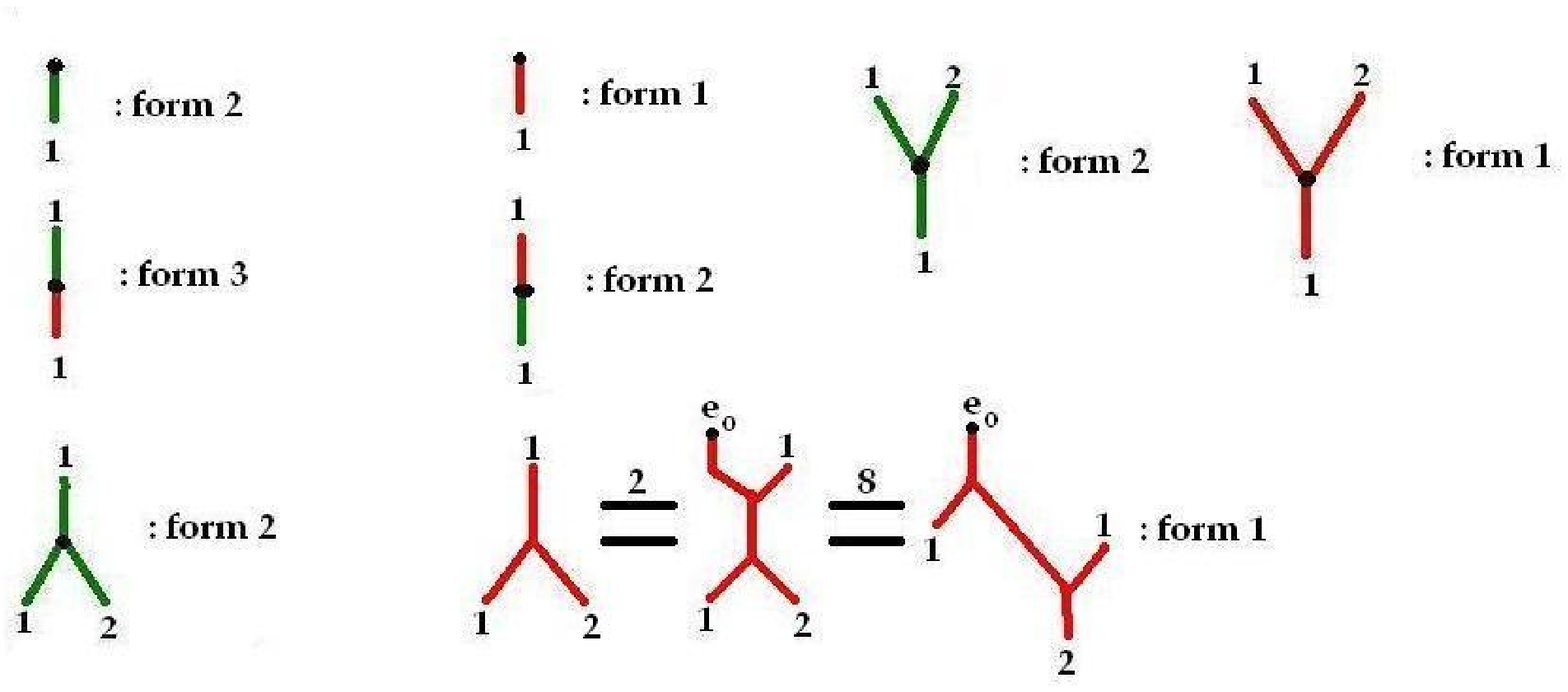}
\end{center}

Next we verify case by case that composing a tree in normal form with a generator gives a tree equivalent to a tree in normal form:
\newpage

\noindent \underline{\Large Form 1}\\
\begin{itemize}
\item Plugging $e_o$ into an input of a form 1 tree results in another tree of form 1
\item Plugging $\phi_{c \mapsto o}$ into an input is equivalent to a tree of form 3 by relation 3b
\item Consider plugging $\phi_{o \mapsto c}$ into an output.  By Claim 2.5 we can assume that the output is the main root output.  Thus we get a tree of form 2.
\item Composing with $m_o$, open multiplication, clearly gives another tree of form 1 after we assume again that when composing with an output of the form 1 tree, this output is the main root output.
\item Finally, consider composing with $\bigtriangledown_o$.  If we assume again that when we compose with an output of the tree it is the main root output, then using one of its normal forms below gives a tree of form 1:
\end{itemize}

\begin{center}
\includegraphics[width=2.5in]{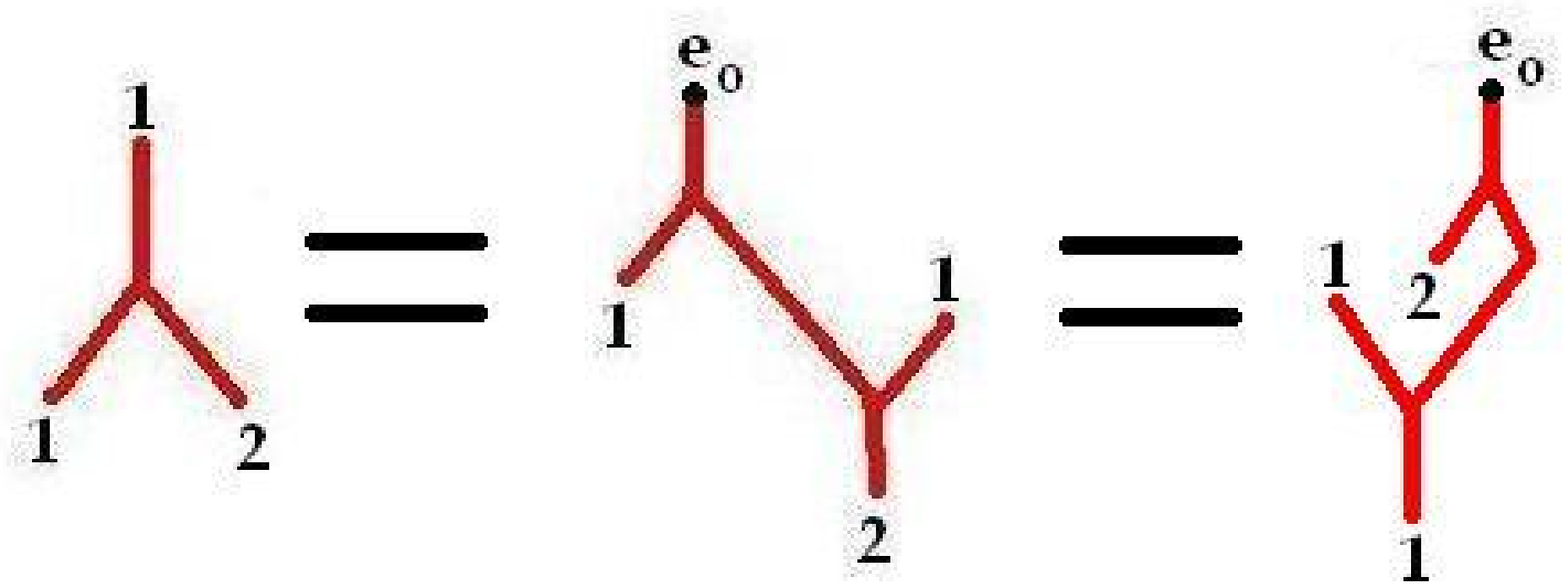}
\end{center}

\noindent \underline{\Large Form 2}\\
\begin{itemize}
\item Plugging in units give trees of form 2
\item If we plug the output of $m_c$ into an input we get another tree of form 2.  If we plug an input into an output, then we can assume the output is the rightmost output stem.  Then relation 7 suffices.
\item Same argument for $\bigtriangledown_c$
\item Plugging $m_o$ into an input gives another tree of form 2.  Now say we plug an input of $m_o$ into an output $o_1$ of the tree.  Then using the form 1 result, we first take an equivalent tree which has $o_1$  as the base output of the red tree:\\
\includegraphics[width=2.5in]{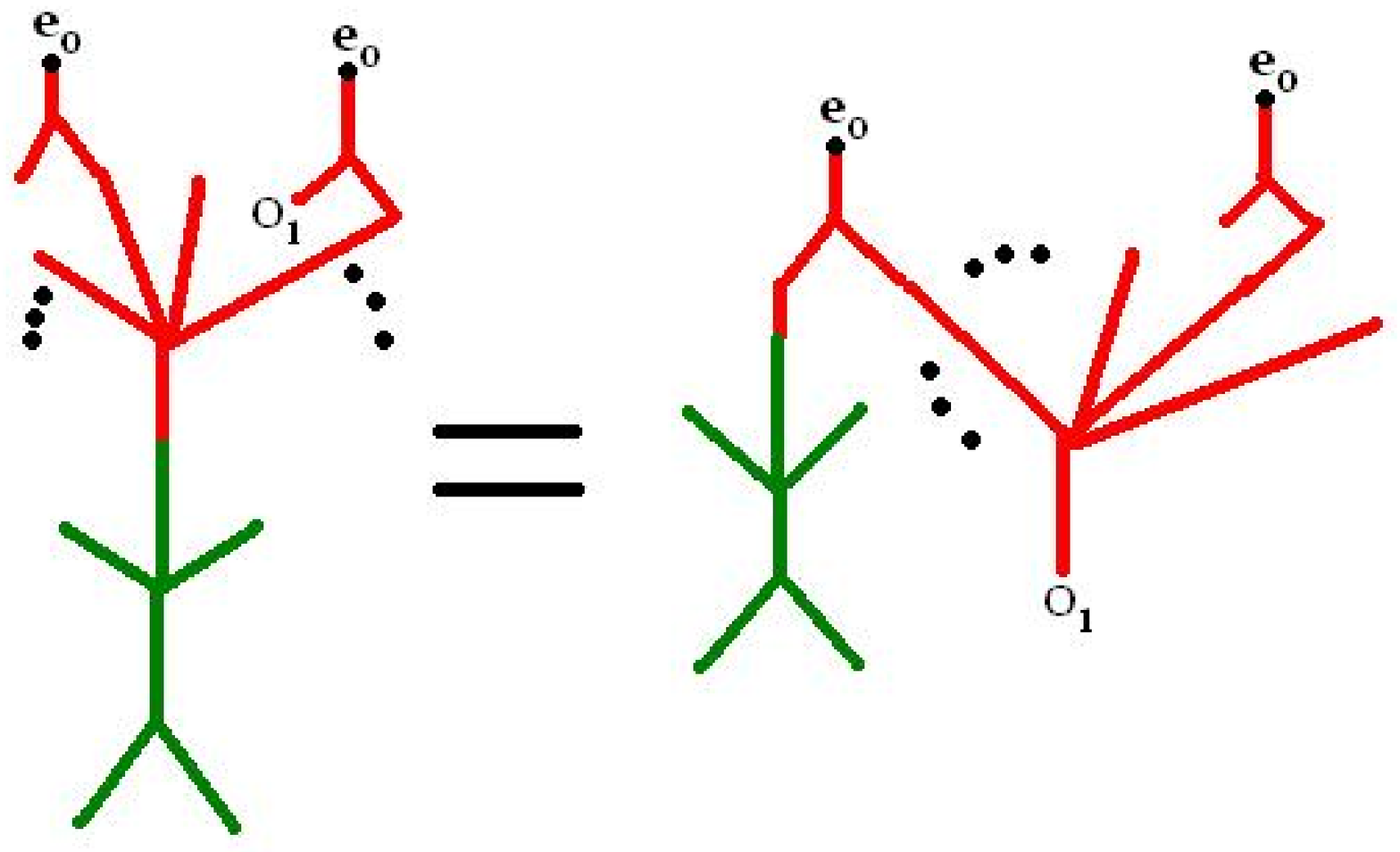}\\
Then plug the input into $o_1$ .  Then "rotate back" so that the green stem is again attached to the main red output stem, resulting in a tree of type 2.
\item The argument for $\bigtriangledown_o$ is the same after we use its appropriate normal form representation as in the form 1 case.
\item Plugging $\phi_{o \mapsto c}$  into an input gives another tree of type 2.  The following relation shows that plugging into an output gives a tree equiv to a tree of form 2:\\
\end{itemize}

\begin{center}
\includegraphics[width=4in]{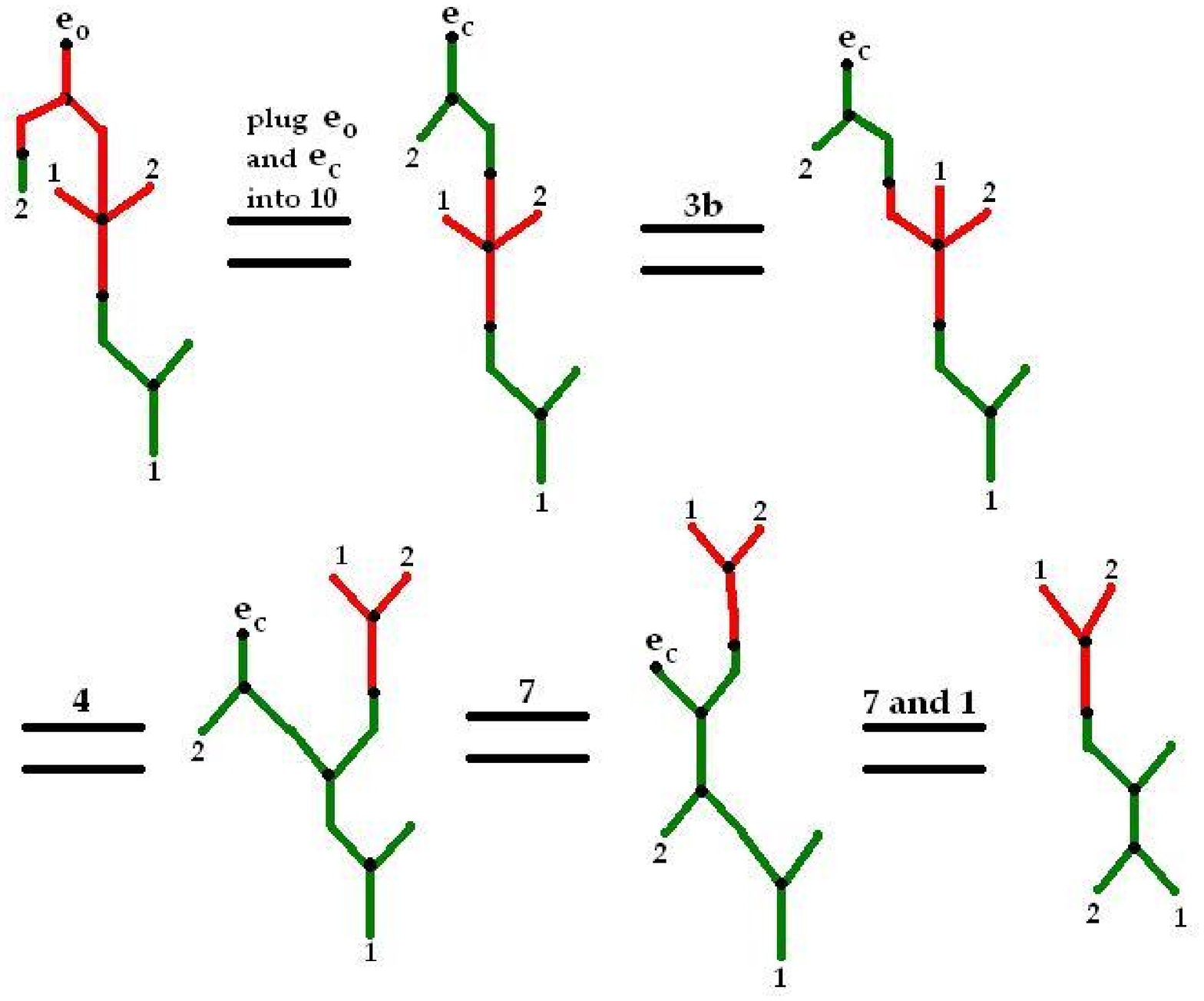}
\end{center}
\begin{itemize}
\item If we plug $\phi_{c \mapsto o}$ into the output of a form 2 tree with one closed output then we get a form 3 tree.  If there is more than one output, then we can assume the output being composed with is the right most output.  Then the following suffices:\\
\begin{center}
\includegraphics[width=3in]{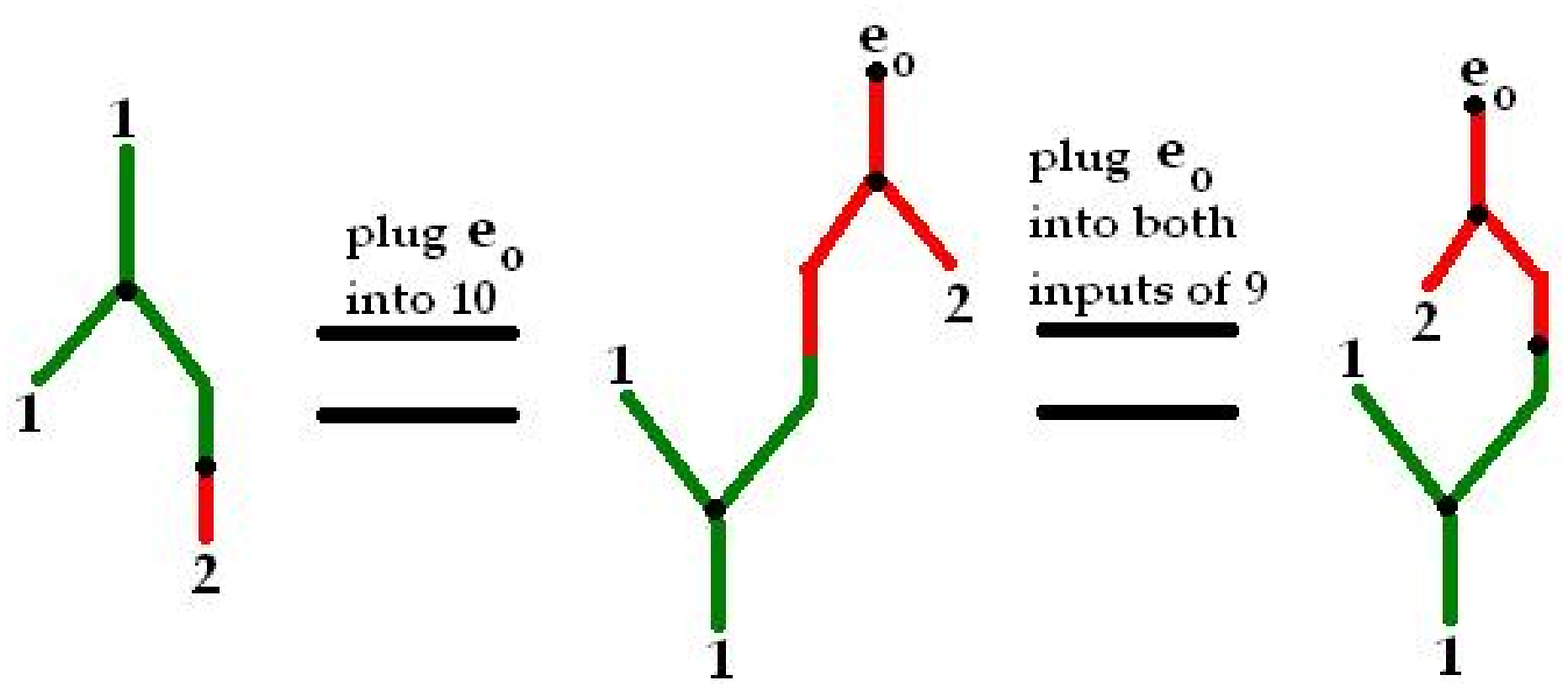}
\end{center}
For plugging into an input, use the following equivalence:\\
\begin{center}
\includegraphics[width=3in]{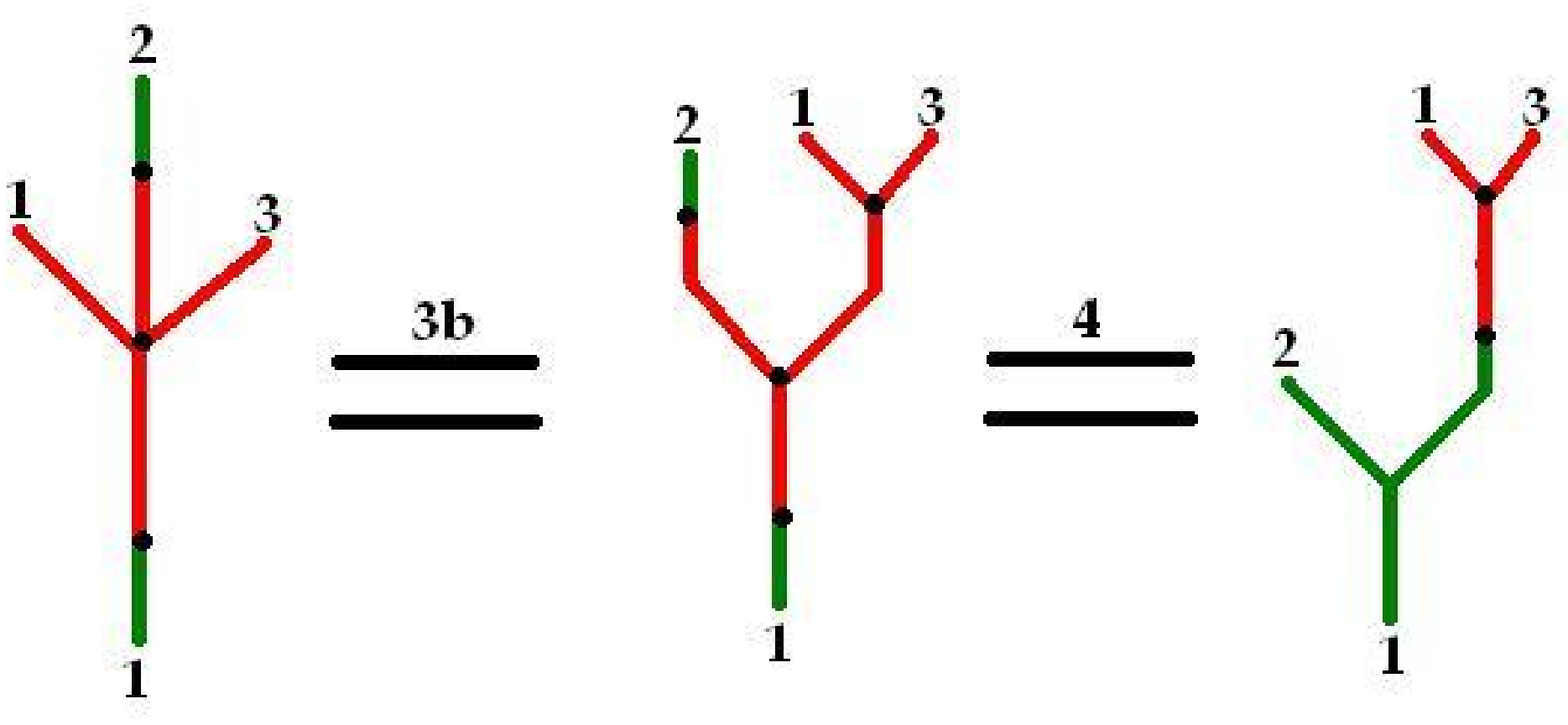}
\end{center}

\end{itemize}
\newpage
\noindent \underline{\Large Form 3}\\
\begin{itemize}
\item Plugging in units gives trees of form 3 (unless the form 2 tree plugged into the left stem is the identity. Then plugging in $e_c$ gives a form 1 tree.)
\item If we plug $\phi_{c \mapsto o}$ into an open input which is not on the bottom red tree, then it is plugged into the form 2 tree and this case has been covered.  If it is attached to an open input on the bottom red tree, then use this equivalence:\\
\begin{center}
\includegraphics[width=3in]{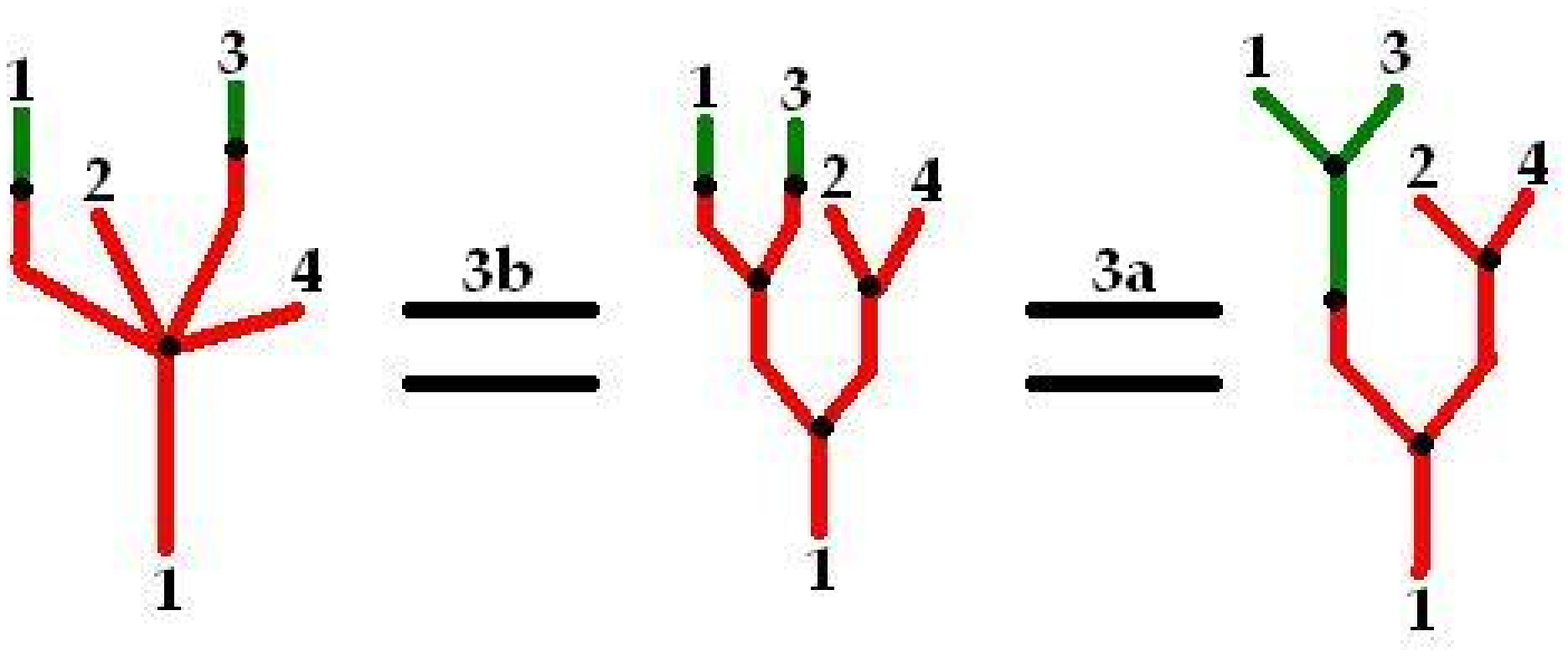}
\end{center}
\item Plugging $\phi_{o \mapsto c}$ into the input of a form 3 tree gives a form 3 tree.  If we plug into an output, then we can assume it is the main root output of the bottom red tree.  Then using relation 4 we immediately see it is equivalent to a tree of form 2.
\item Attaching $m_o$ and $\bigtriangledown_o$ works as before since, again, we can assume the output being plugged into is the main root output.
\item Plugging $m_c$ into an input gives another form 3 tree (and there are no closed outputs to compose with).
\item For plugging $\bigtriangledown_c$ into an input, we should get a tree equivalent to form 2.  To see that this is so, the following suffices:\\
\begin{center}
\includegraphics[width=4in]{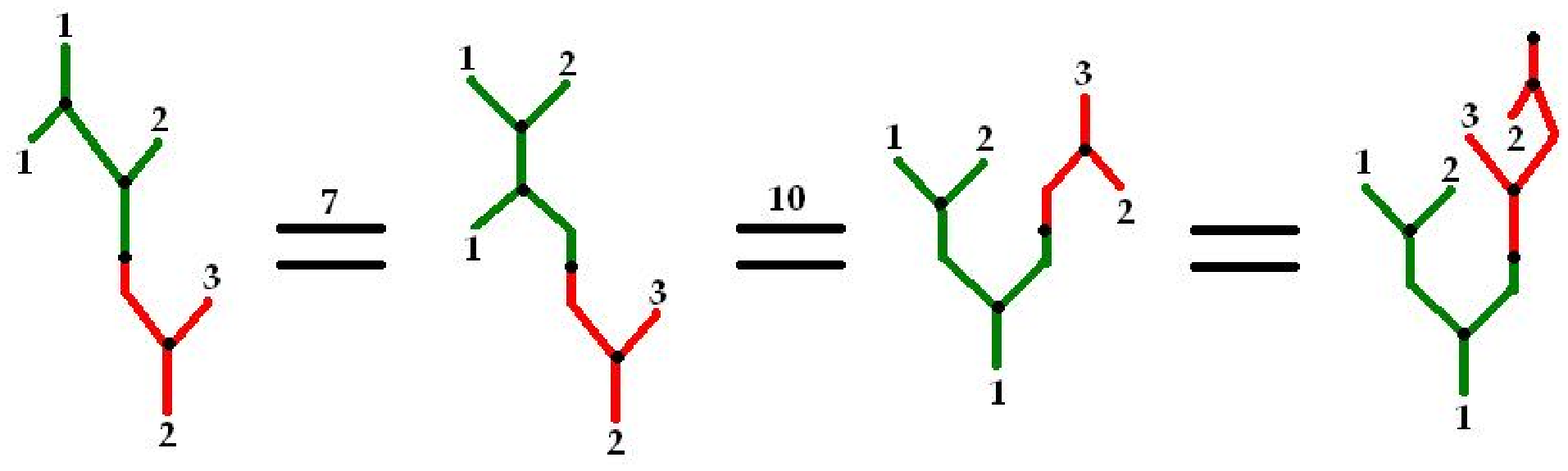}
\end{center}
\end{itemize}
\noindent This completes the form 3 case.

So inductively we now have that any tree formed by composing generators is equivalent to a tree in normal form.  And since the action by the symmetric groups is invariant on the set of trees in normal form, we see that all trees are equivalent to a tree in normal form.
\vspace{.2in}

\noindent This concludes the proof of Theorem 2.3.

\end{proof}

\section{Description of $H_*(C)$}

Before looking at the full $H_*(OC)$ let's first see what we need to add to our green degree 0 generators and relations in order to give a complete description of $H_*(C)$, where C is the moduli space of Riemann spheres with closed inputs/outputs and no boundary.
We add to the list of generators the degree one BV operator coming from $H_1(C)$ given by rotating the input parametrization 360 degrees:\\

\begin{description}
\item[g9]{\Large\bf BV Operator}

$$\begin{array}{ll}
\put(0,-50){\includegraphics[width=2.5in]{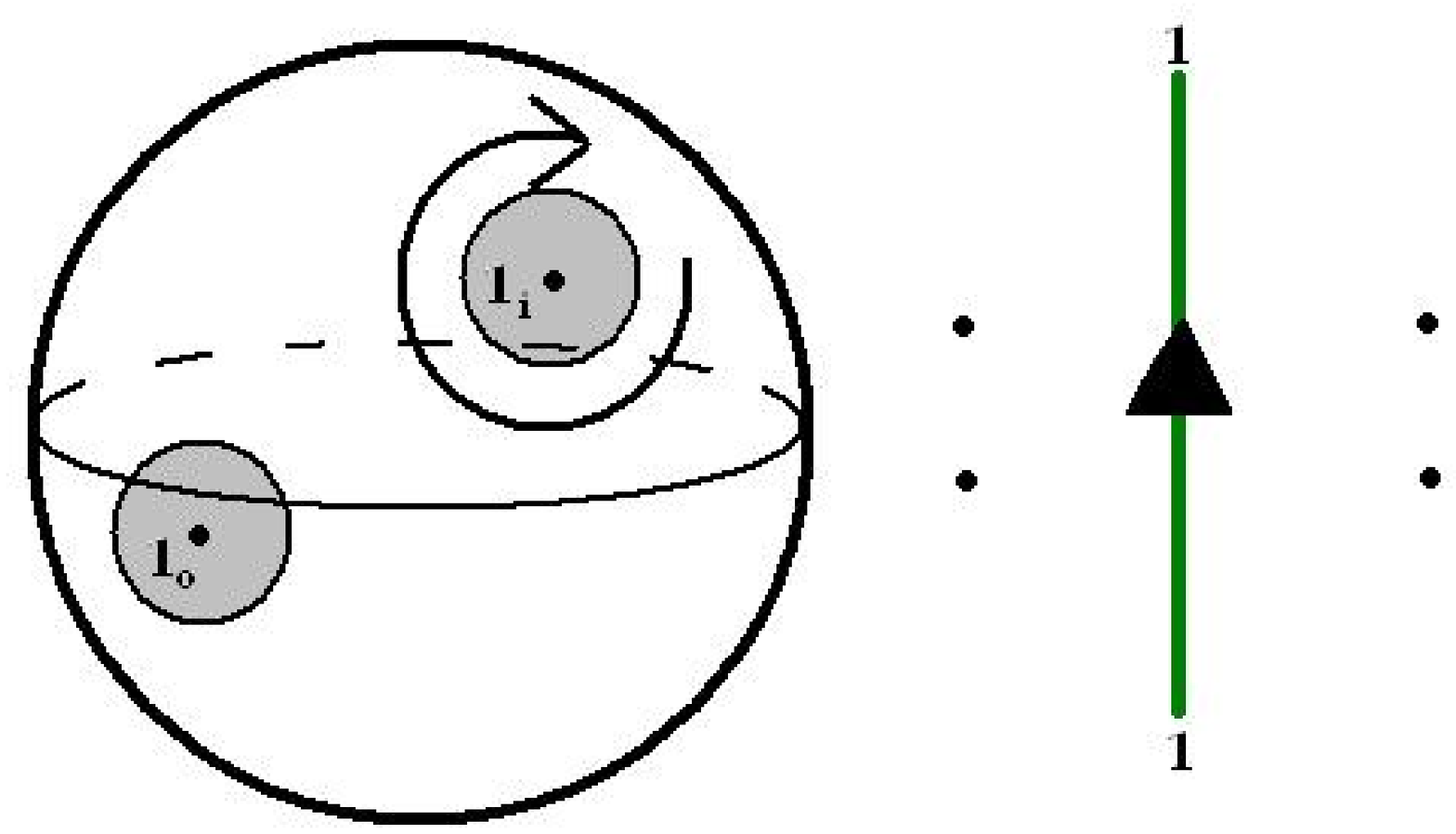}} \hspace{2.5in} & \mbox{{\LARGE $\Delta : V_c \mapsto V_c$}} \end{array}$$

\end{description}

\noindent {\large Then we have the following new relations:}\\

\begin{description}
\item[r11] The usual BV relations giving that $\Delta$ is a second order derivation of closed multiplication and $\Delta^2=0$:
\begin{center}
\includegraphics[width=5in]{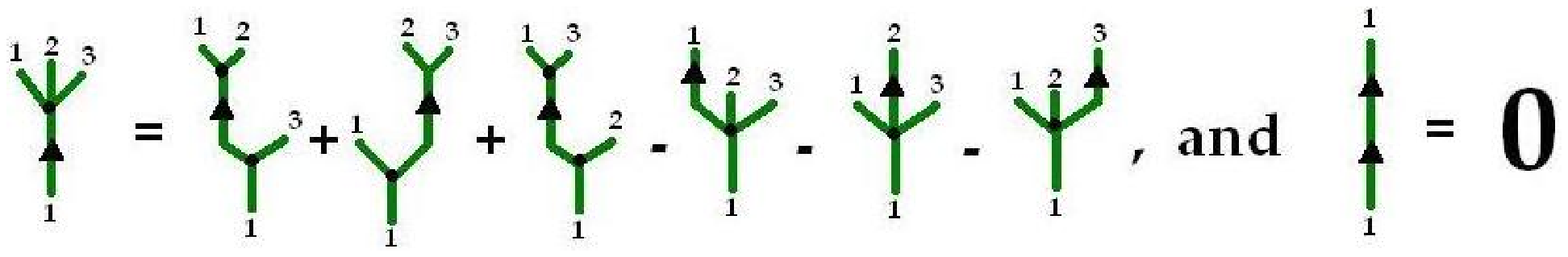}
\end{center}
\item[r12] The following relation holds since rotating the input on a sphere with an input at the north pole and an output at the south pole is conformally equivalent to rotating the output:\\
{\Large $(\bigtriangledown_c(a) \bullet \Delta ) \bullet b=a \bullet (\Delta \bullet \bigtriangledown_c(b))$}
\end{description}

\begin{center}
\includegraphics[width=2.5in]{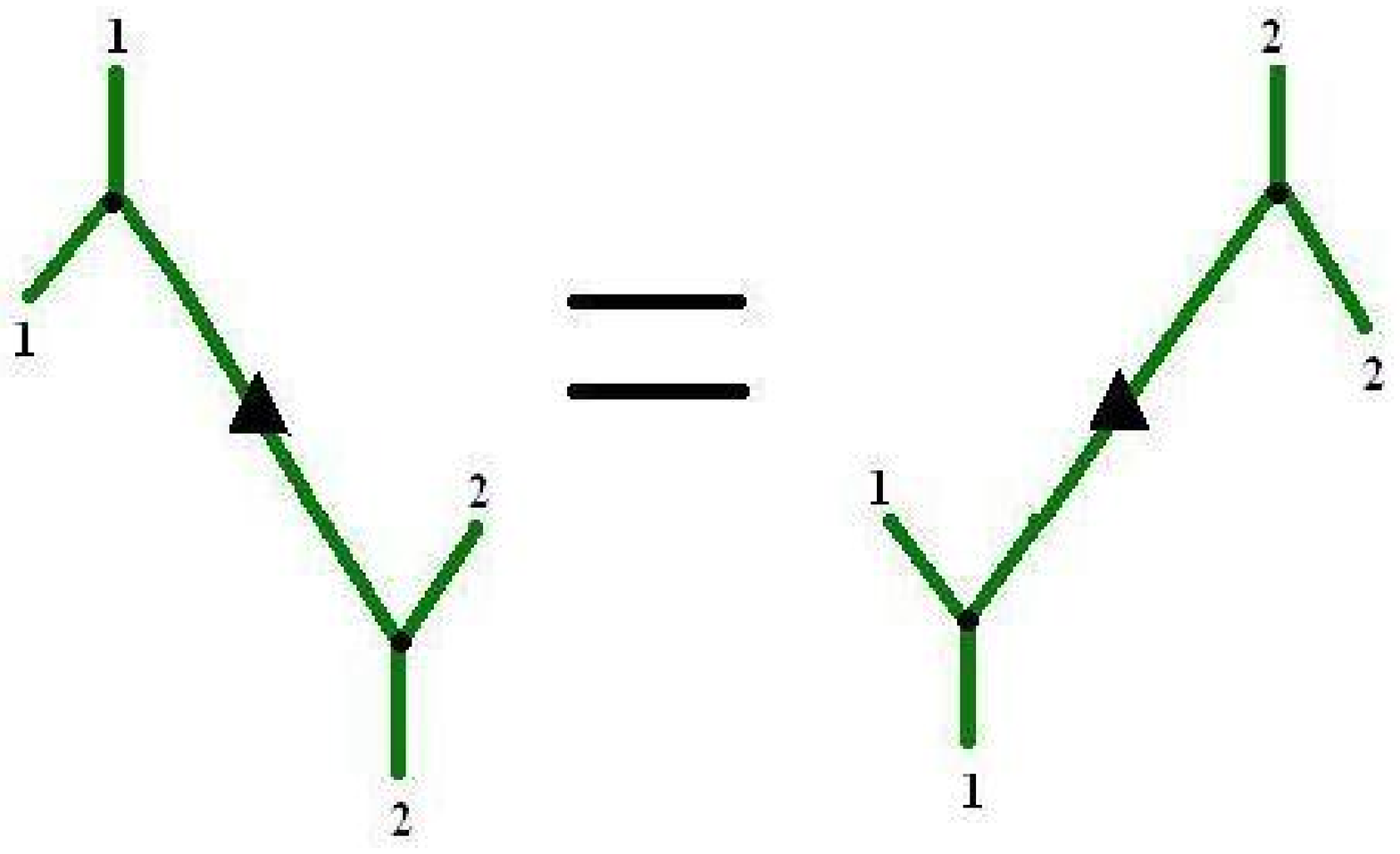}
\end{center}

\begin{thm}
The green (closed) generators and green relations completely describe $H_*(C)$.
\end{thm}
\begin{proof}
Again we consider the dioperad morphism from the free dioperad generated by the green trees representing our green generators  modulo the green relations to $H_*(C)$.  By Getzler's result \cite{getz1} describing an algebras over the homology of the operad formed by spheres with closed inputs and one closed output as BV algebras, we know that there is a vector space isomorphism between $H_*(C)$ and the span of equivalence classes of trees of the following form:\\

$$\begin{array}{ll}
\put(0,-70){\includegraphics[width=2in]{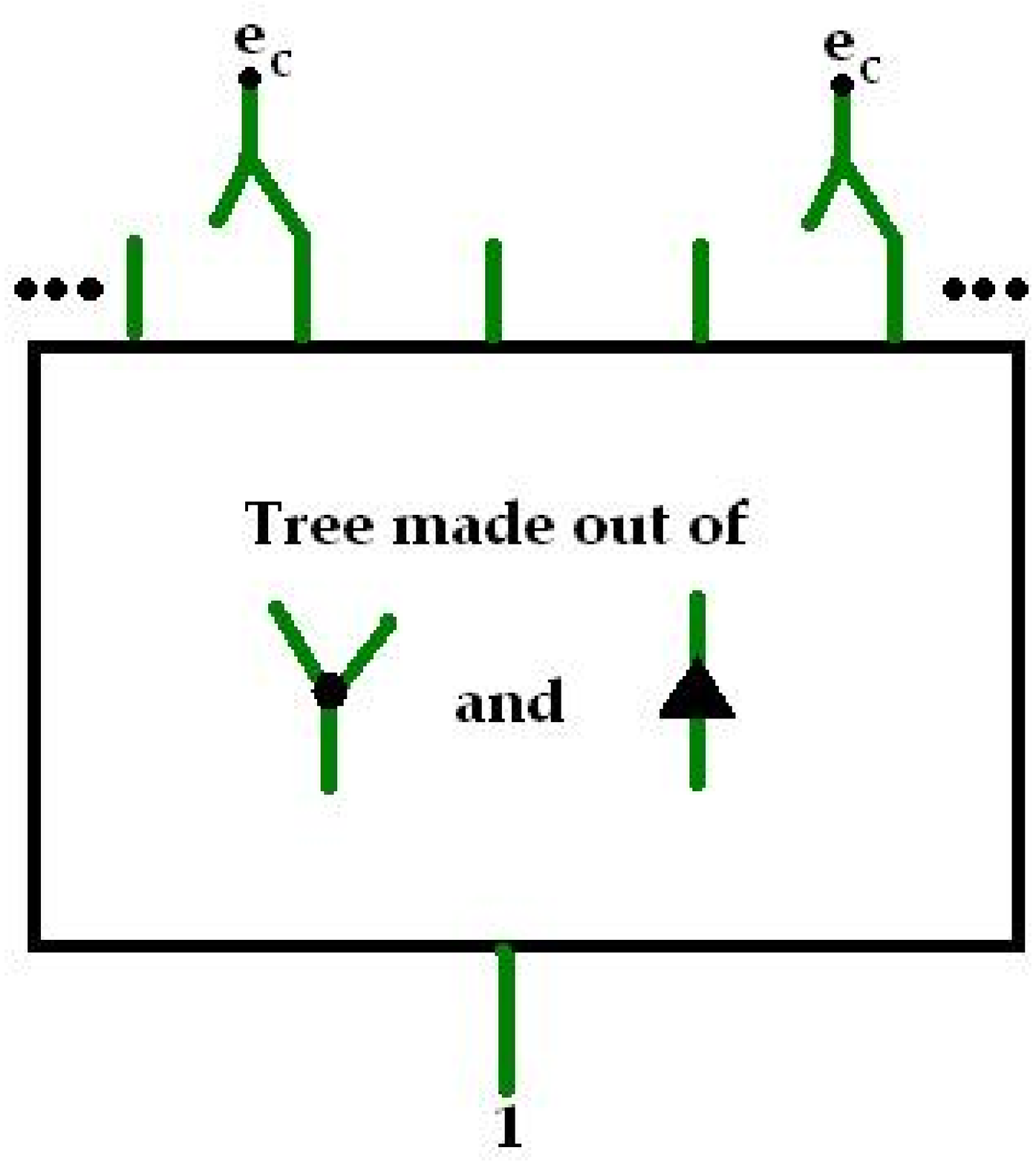}} \hspace{2in} & \begin{array}[b]{l} \mbox{Here it is required that the main}\\ \mbox{root output be labeled 1.}\end{array} \end{array}$$

\begin{df}
Now say a green tree is in {\em normal form} if it is of the above form except that the main root output is not required to be labeled 1.
\end{df}

\begin{claim}
Restricting the morphism to the span of equivalence classes of trees in normal form gives a vector space isomorphism onto $H_*(C)$.
\end{claim}

\begin{proof}
By the above fact, it suffices to show that any tree in normal form is equivalent to another tree in normal form which has its main root output labeled 1.  We can assume that the output labeled 1 is connected to the leftmost input of the tree made out of the multiplication and BV generators.  Then it suffices to check the claim for trees like the following:\\

\includegraphics[height=2in]{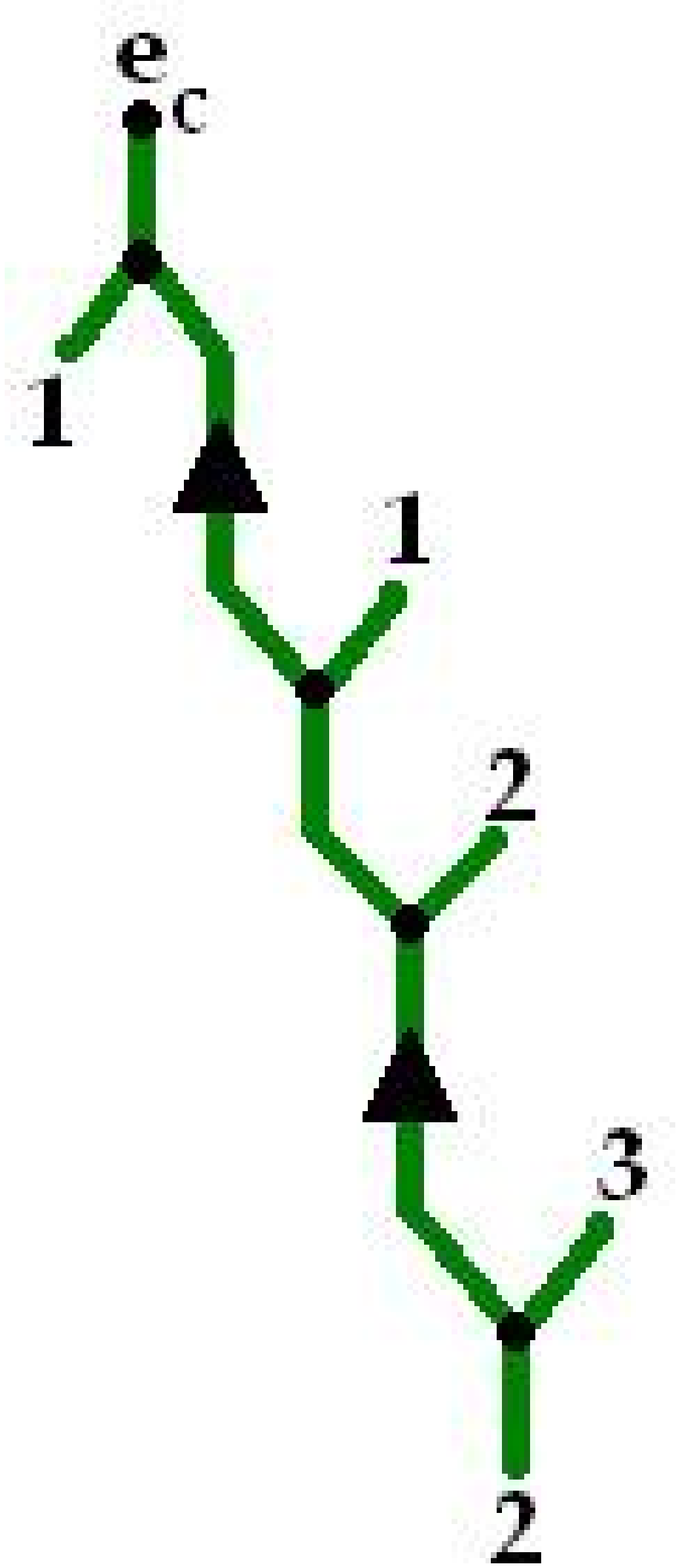}

I.e., it suffices to check it for trees which are formed by taking a tree made out of closed multiplication and the BV operator, where the right input of closed mult. never has anything composed into it, and plugging $\bigtriangledown_c$ and $e_c$ into its top left input.\\

But it is not too difficult to see how to achieve this using relations r7 and r12.\\
For example, for this tree we have:
\begin{center}
\includegraphics[width=4in]{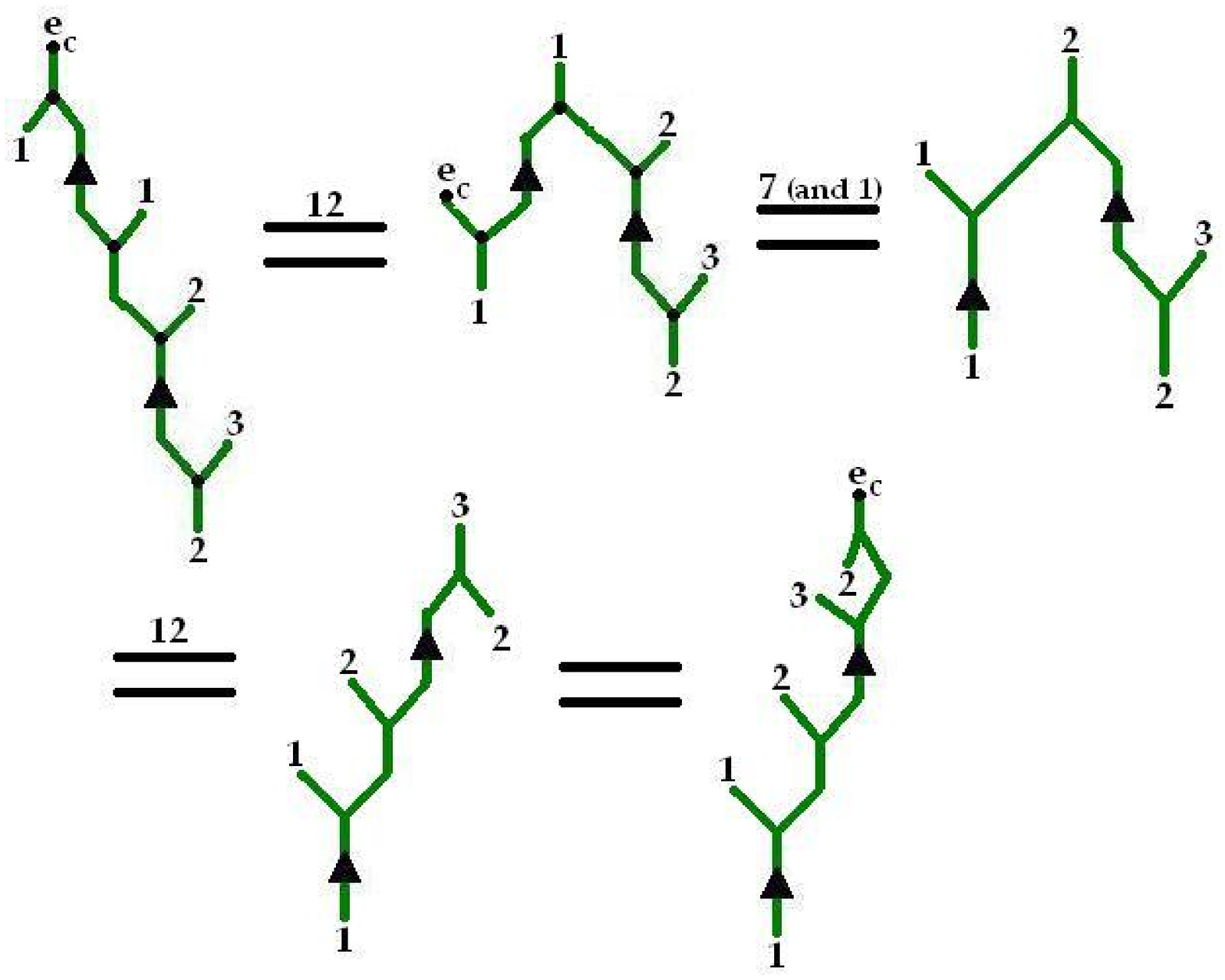}
\end{center}

\noindent \underline{\large\bf Warning}:  All equivalencies before this last one sufficed in showing the desired property, i.e. the property was directly implied by the equivalence.  This last equivalence is just an example for this particular tree and does not imply the property for all the trees that the property needs to be shown for.  However, it is easy to see that the method above of using a sequence of applications of relations 7 and 12, not necessarily alternating, will work for all of the trees that we need to check the property for.

\noindent Thus claim 3.3 is true.
\end{proof}

So all that we need to do is show that every tree in the dioperad is equivalent to a tree in normal form.  We proceed as in the last section by first noting that all the green generators are equivalent to normal form.  Next, we need to check that composing a normal form tree with a green generator gives a tree equivalent to a tree in normal form.  But using one of $\bigtriangledown_c$'s normal form representations,  we automatically get trees in normal form when composing with any green generator since by above we can assume that the when we compose with an output, that output is the main root output.

Since, again, the symmetric groups act invariantly on the normal form trees, all trees are equiv. to a tree in normal form.\\
Thus Theorem 3.1 is proved.

\end{proof}

\section{Description of $H_*(OC)$}

Taking care of $H_0(OC)$ and $H_*(C)$ first makes the arguing for go smoothly since now we can replace any degree 0 tree with another tree which goes to the same path component and we can put any green tree in the normal form of Sec 3 with the main output labeled as we desire.
We do not need to add anymore generators to our list of 9, and there is only one more relation that we need to add:\\
\begin{description}
\item[r13] $\Delta(\phi_{o \mapsto c}(e_c))=0$.  This relation comes from the fact that rotating the boundary component in a sphere with a closed output over the south pole and an empty boundary component over the north pole 360 degrees gives a constant map into the moduli space since boundaries are not parameterized.
\end{description}
\begin{thm}
The list of  9 generators (pg. 5-6,19) and 13 relations (pg. 6-9,19,22) completely describe $H_*(OC)$.
\end{thm}

\begin{proof}
Let's first discuss the homotopy type of an arbitrary path component of $OC$.  First note that $OC$ is h.e. to the moduli space of Riemann spheres with boundary, with labeled punctures in the interior and on the boundaries (no parameterizations in either case), and each puncture in the interior coming with a tangent direction.
Now consider a path component with $n$ interior punctures with directions, $k$ boundary components each having exactly one puncture, and $l$ empty boundary components.  Then up to homotopy we can replace a marked boundary component with a puncture and a direction and an empty boundary component with a puncture.  Thus it is h.e. to the moduli space of spheres with $n+k$ labeled punctures with directions and $l$ unlabeled punctures without directions.
For an arbitrary path component $P$ in $OC$, consider the forgetful fibre bundle $P \mapsto \bar{P}$ where  $\bar{P}$ is the moduli space obtained from $P$ by dropping all but one fixed open puncture on each non-empty boundary component.  Then it is clear that the fiber is contractible.

Using the above facts and mixing the descriptions of the homologies of the framed and non-framed little disks operads, we can get a vector space isomorphism between $H_*(P)$ and the span of the equivalence classes of a set of trees in $F(G)/S$, the free dioperad on all our generators modulo all our 13 relations.  We can see how the general case works while avoiding indexing messiness by assuming $P$ is of type, say, $\{1_i,2_i,1_o,2_o\},(3_i,3_o,4_o,4_i),(5_i,6_i),(5_o,6_o),(),()$.  If we choose a closed output, say $1_o$, then a set for which this is true is the set of all trees of  the form:\\
\begin{center}
\includegraphics[width=4in]{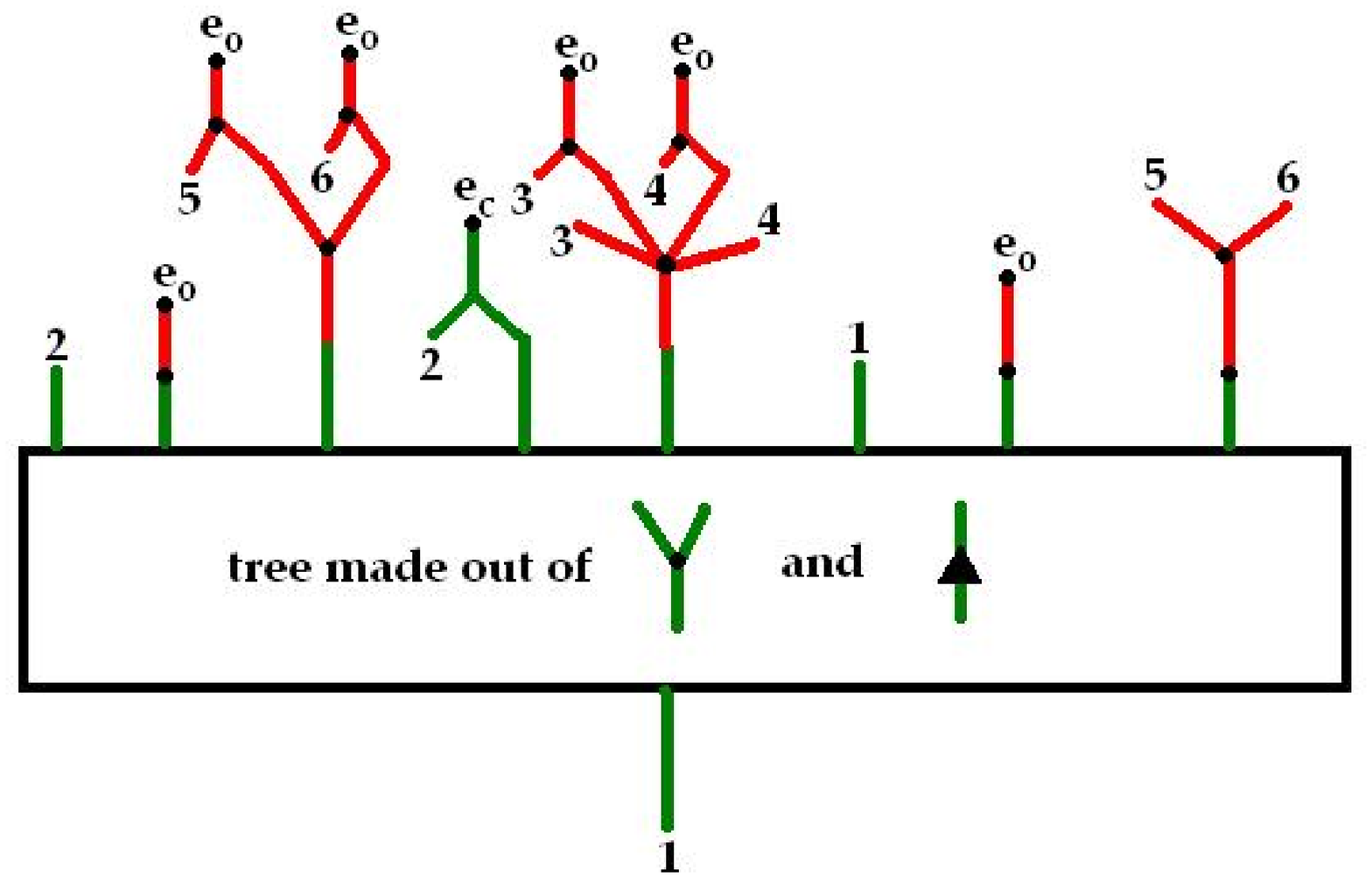}
\end{center}
\begin{center}
{\Large The main output stem must be labeled $1$}
\end{center}
\vspace{.2in}

If we choose an open output, say $3_o$, then another set whose equivalence classes span a vector subspace isomorphic to $H_*(P)$ is the set of trees of the following form:\\
\begin{center}
\includegraphics[width=4in]{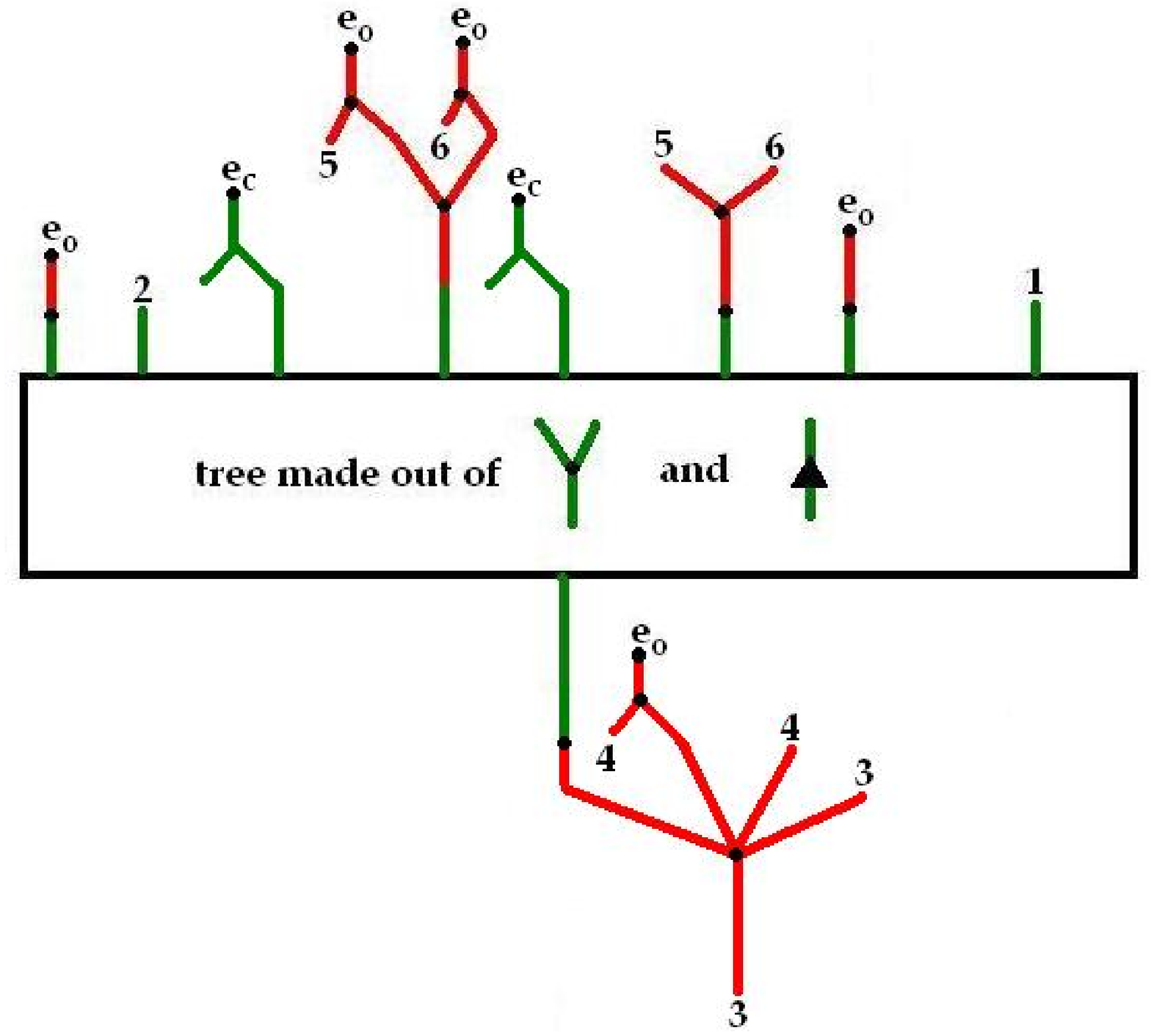}
\end{center}
\begin{df}
With this in mind we'll now say that a tree is in normal form if it is of either of the 2 forms above (with any labeling).
\end{df}

Note that since we allow the green tree in the second form to be $e_c$ , and $\phi_{c \mapsto o}(e_c)=e_o$, the second form contains completely red trees.  Thus the above two forms give all path components.
\begin{claim}
Restricting the morphism $f(G)/S \mapsto H_*(OC)$ to the span of equivalence classes of trees in normal form gives an onto vector space isomorphism.
\end{claim}
\begin{proof}
By the above discussion, along with the $H_0(OC)$ result, all that needs to be shown is that for any tree in normal form, and for any of its outputs, it is equivalent to another tree in normal form which has this output as the main root output.  But this follows from the same result for the green normal form trees of the last section and the $H_0(OC)$ result. For example, this shows the property for a tree of form two and an open output:\\
(we can assume the red tree containing the output is on the left most green input of the green tree)\\

\begin{center}
\includegraphics[width=4in]{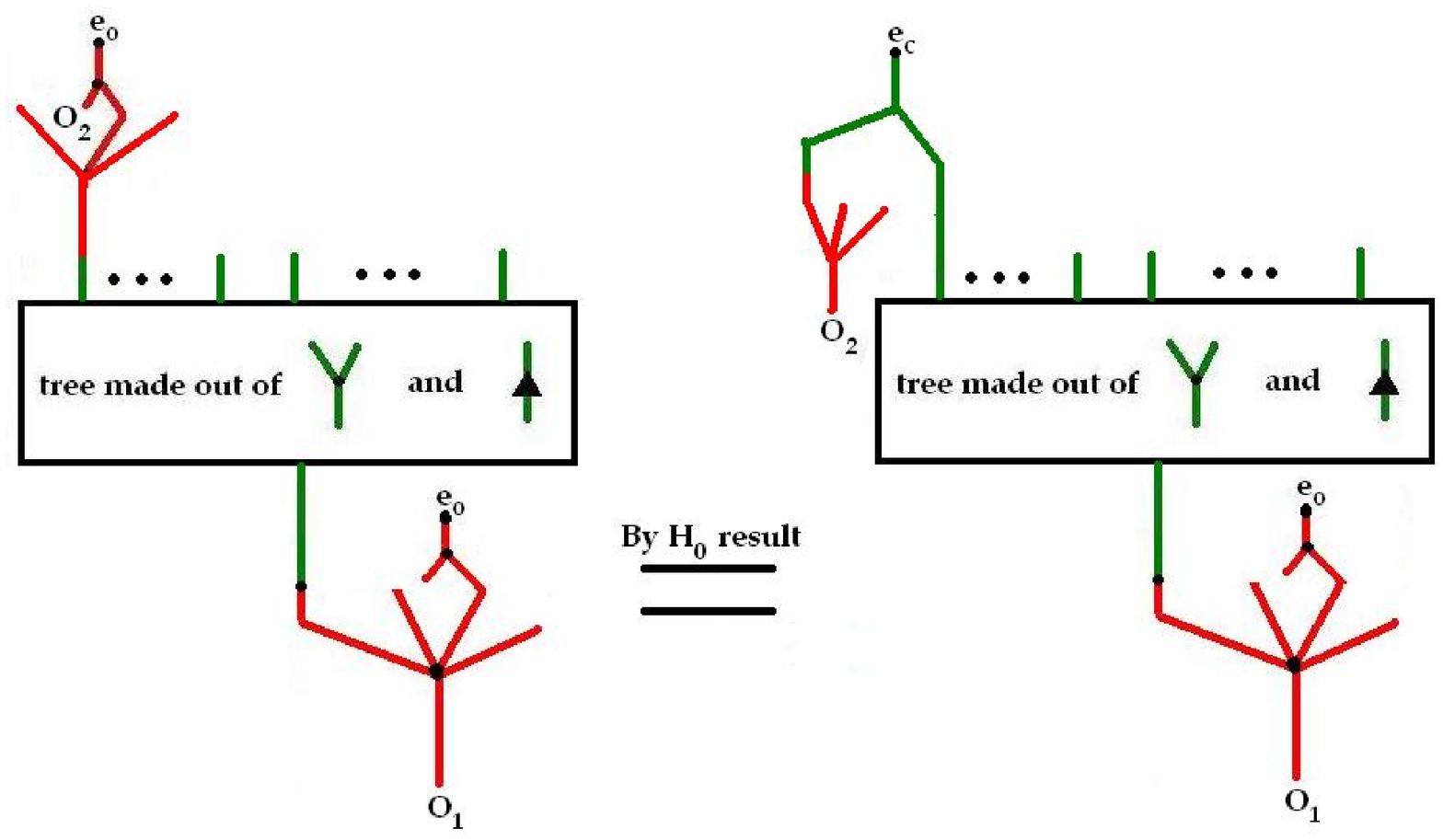}
\end{center}

\includegraphics[width=2.3in]{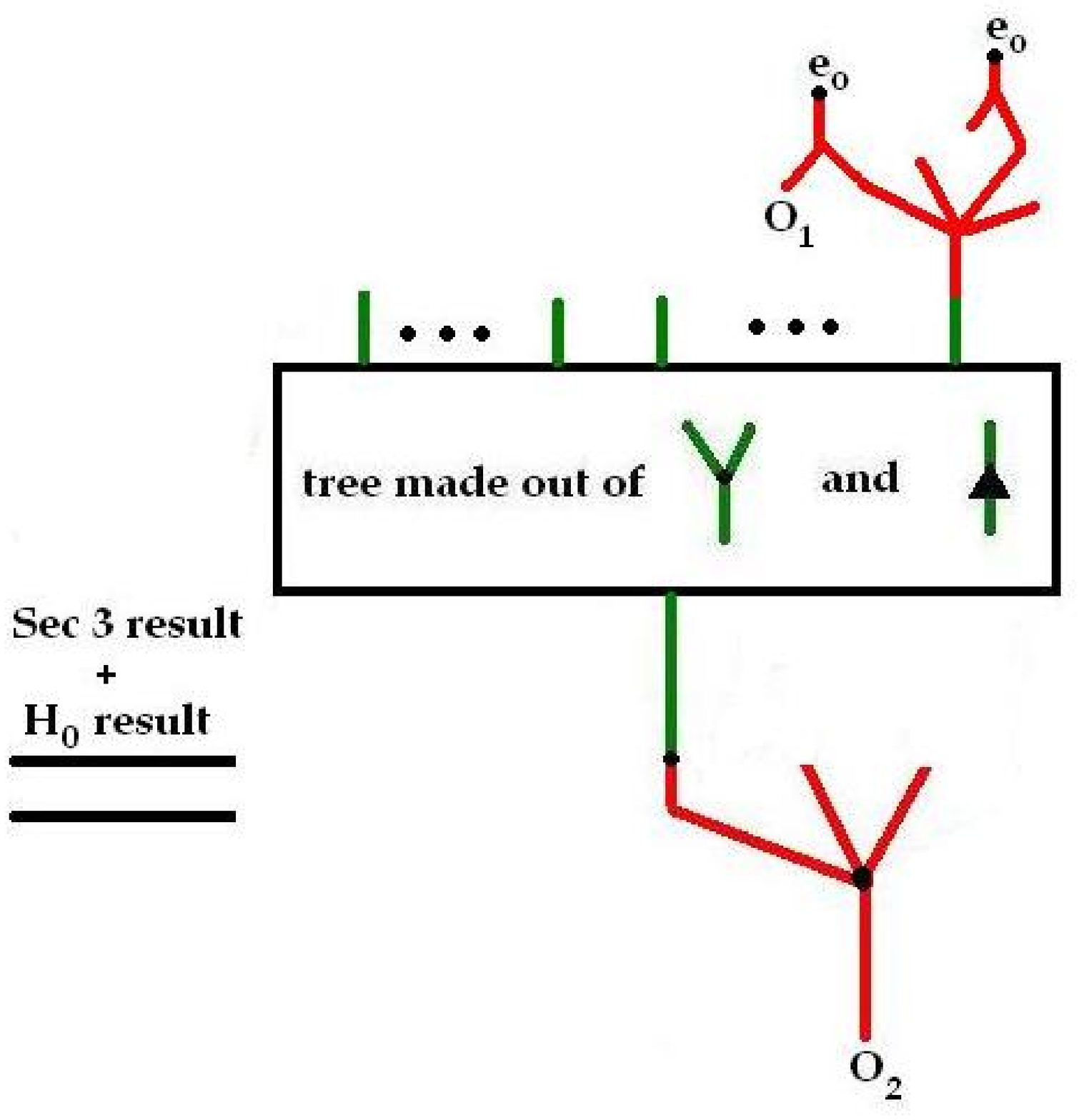}

\noindent This proves Claim 4.3.
\end{proof}

So in order to finish the proof Theorem 4, all that needs to be shown is that any tree is equivalent to a tree in normal form.  This again can be done inductively as in the last two sections.  But this is now straight forward to check given what has been developed so far and can be left to the reader.\\

\noindent This concludes the proof of the main Theorem 4.1

\end{proof}

\section{Semi-modular and cyclic structure of operad}

Restricting to the generators with one output and the relations only involving these generators, we get a complete description of the homology of the operad formed by the components of $OC$ with only one output (we'll abuse notation and call this operad $OC$ also).  This can be seen by restricting the proof above to trees with only one output. This operad is cyclic (see \cite{getz2}) in the sense that there is no natural output, requiring us to label it. The action which permutes all the labels extends the action which only permutes the input labels and it does it in a composition respecting way. Thus the homology forms a 2-colored cyclic operad in the category of graded vector spaces.

The definition given in \cite{getz2} of the cyclic endomorphism operad for a graded vector space $V$, finite dimensional in each degree, and a non-degenerate inner product $B$ on $V$ can be naturally extended to the colored case.  In our 2-colored case, the two vector spaces $V_c$ and $V_o$ come with inner products $B_c$ and$B_o$ which are used to identify $HOM(V_{i_1} \otimes V_{i_2} \otimes ... \otimes V_{i_n},V_j) \mbox{ with } HOM(V_j \otimes V_{i_1} \otimes V_{i_2} \otimes ... \otimes V_{i_n},F) \quad (F the base field) \mbox{ where } i_k,j \in \{C,O\}$.

Then since the permutation $(1, 2, 3)$ sends $m_c$ and $m_o$ to themselves,$(1, 2)$ sends $\Delta$ to itself, and $(1, 2)$ interchanges $\phi_{c \mapsto o}$ and $\phi_{o \mapsto c}$, we get that an algebra over the cyclic $H_*(OC)$ satisfies the following four additional relations:\\

\noindent {\large\bf Additional relations in an algebra over the cyclic operad $H_*(OC)$:}
$$\begin{array}{ll}
1) B_c(a,bc) = B_c(ab,c) \hspace{1in} & 2) B_o(a,bc) = B_o(ab,c) \\
3) B_c(\Delta(a),b)=(-1)^{|a|}B_c(a, \Delta(bc)) & 4)B_o(\phi_{c \mapsto o}(a),b)=B_c(a,\phi_{o \mapsto c}(b))
\end{array}$$
\vspace{.2in}

The operad $OC$ does not form a modular operad (see \cite{getz3}). This is because self-sewing in general results in a surface of genus $> 0$. However, sewing two open inputs on the same boundary component gives another sphere with one more boundary component:\\

\noindent $$ \begin{array}{lll} \put(-175,-75){\includegraphics[width=2.5in]{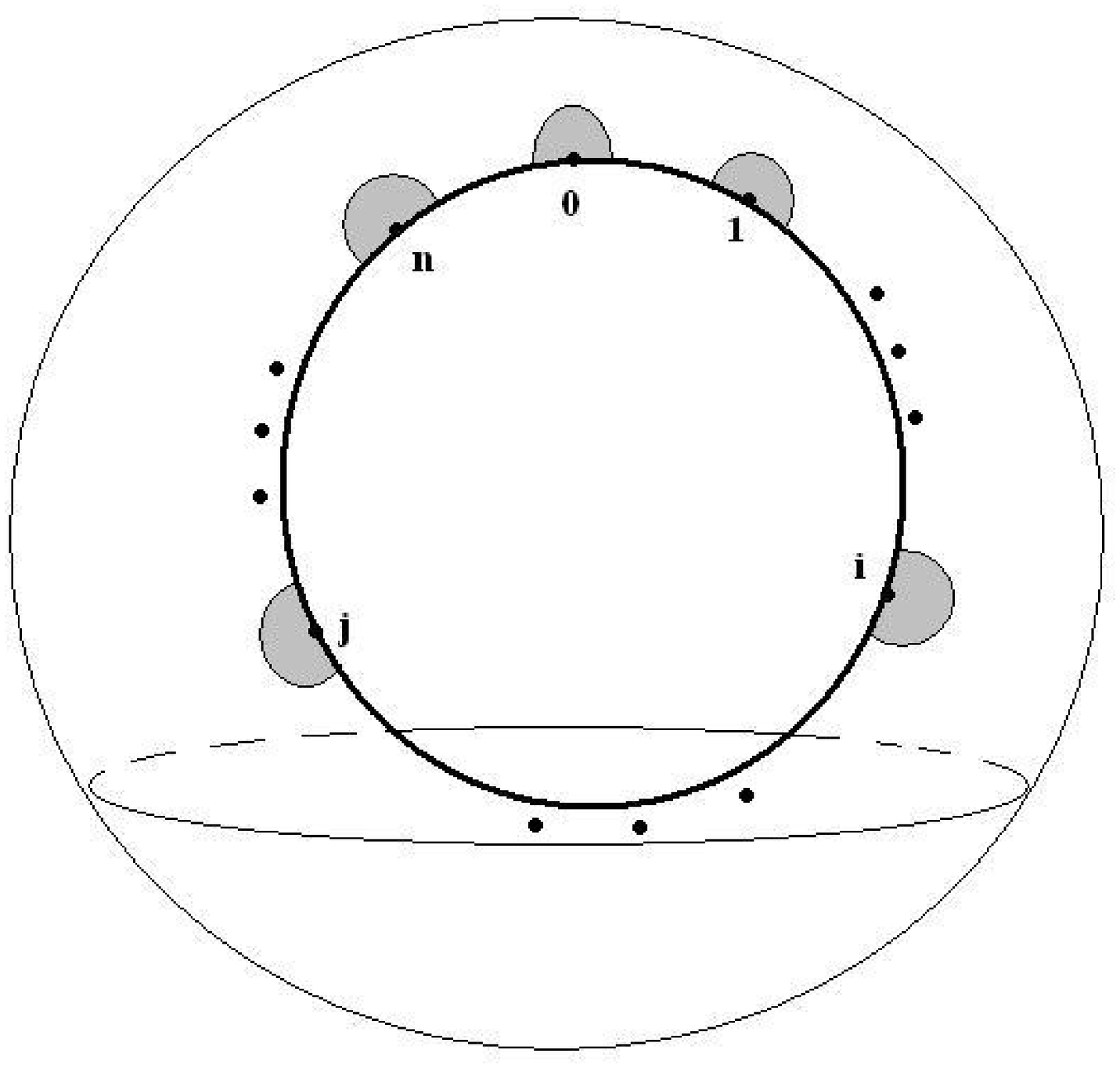}}& \begin{array}{c} \mbox{{\LARGE\bf $\Psi_{ij}$}}\\
\mbox{{\LARGE\bf $\longmapsto$}} \end{array} & \put(0,-75){\includegraphics[width=2.5in]{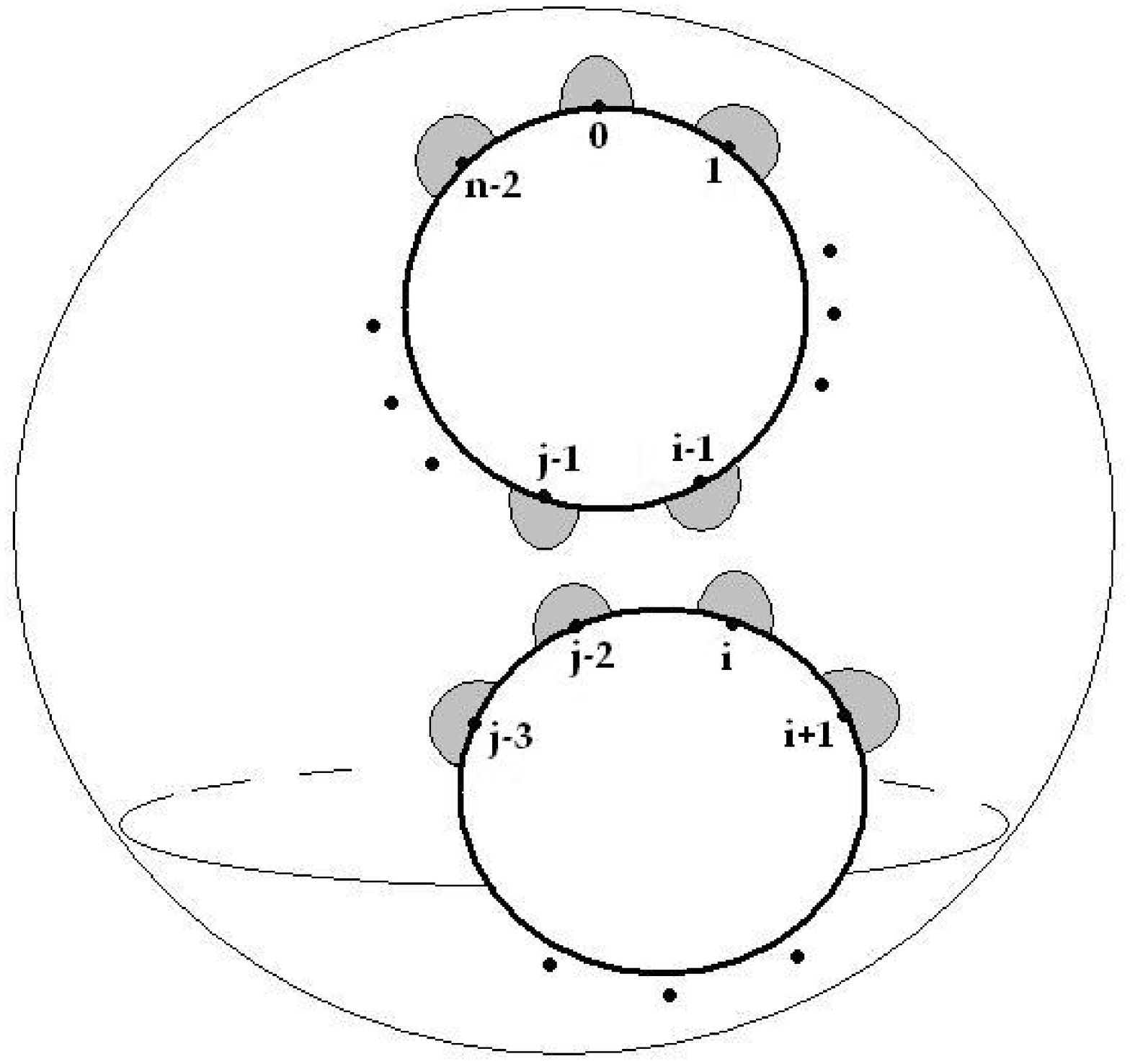}}
\end{array}$$

To see what extra algebra structure this adds, we need an appropriate endomorphism definition. In
order for a map from $H_*(OC)$ to this endomorphism operad to respect the contractions $\Psi_{ij}:H_*(OC(n)) \mapsto H_*(OC(n-2))$, we need the endomorphism contractions $\Psi_{ij}:End(n) \mapsto End(n-2)$ to be zero when applied to a homomorphism which is in the image of the homology of a path component in which the $i$ and $j$ are not open inputs on the same boundary component.  With this in mind we make the following:

\begin{df} {\bf\em OC cyclic semi-modular endomorphism operad}\\
For $(V_c,B_c),(V_o,B_o)$ as above, let $End(n)=\bigoplus_{type} End(type,n)$ where $type$ runs over all path component types of $OC$ with $n$ inputs and $End(type,n)=Hom(V_{i_1} \otimes ... \otimes V_{i_n},V_j)$, the $i_k's$ and $j$ either $O$ or $C$ depending on the type.  We have the cyclic structure as above, and we define the compositions so that $f \circ_i g \in End(type(f) \circ_i type(g),n+m-1)$.  We use $B_c$  and $B_o$ to identify $Hom(V_{i_1} \otimes ... \otimes V_{i_n},V_j)$ with $V_j \otimes V_{i_1} \otimes ... \otimes V_{i_n}$ and use $B_o$ to define the contractions $\Psi_{ij}$,just as in the definition of the modular endomorphism operad given by Getler in \cite{getz3}, {\em provided} the type has $i$ and $j$ as open inputs on the same boundary component.  Otherwise $\Psi_{ij}$ is defined to be zero.  $\Psi_{ij}$ should take $f \in End(type,n)$ to $End(type',n-2)$, $type'$ being the path comp of $OC$ containing the image of points in $type$ under self sewing input $i$ to $j$.
\end{df}
\vspace{.2in}
Considering the normal form representation of an element in $H_*(OC)$ from Section 4 and the fact that contractions are defined only for inputs on the same boundary component, we see that the extra structure is completely determined by the $\Psi_{ij}$'s restricted to completely red trees.  So we get one more relation:
\begin{center}
\includegraphics[width=4in]{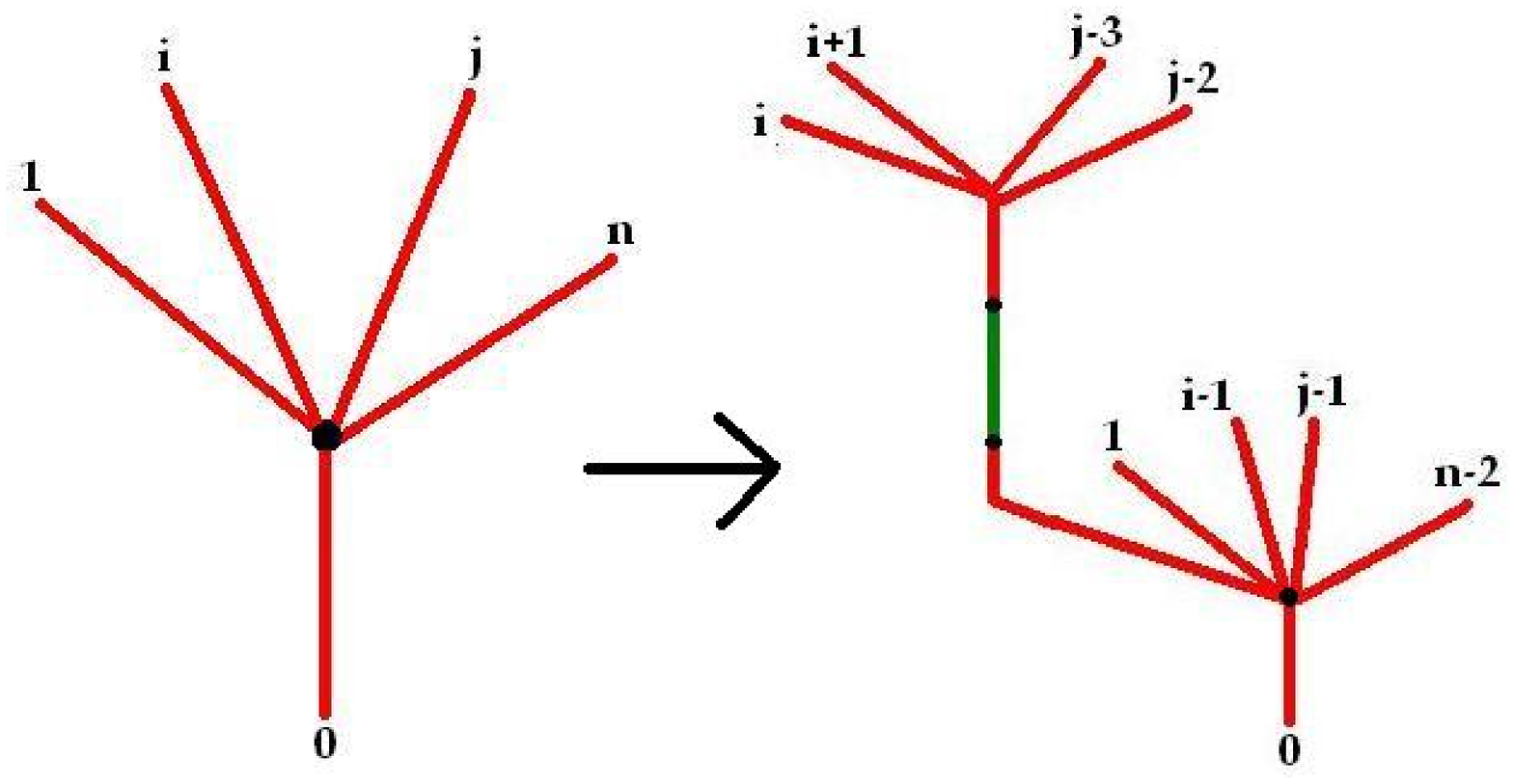}
\end{center}

(If $i=j-1$  then the top red tree is just $e_o$)\\
\\
Algebraically, this relation says if $m_o(n):V_o^{\otimes n} \mapsto V_o$ is the operation given by associative multiplication, $m_o(n)(a_1,...,a_n)=a_1 \cdots a_n$, then \\ $\Psi_{ij}(m_o(n))(a_1,a_2,...,a_{i-1},a_{i+1},...,a_{j-1},a_{j+1},...,a_n)=(a_1a_2 \cdots a_{i-1}a_{j+1} \cdots a_n)\phi(a_{i+1} \cdots a_{j-1})$ where $\phi=\phi_{c \mapsto o} \circ \phi_{o \mapsto c}$ and an empty parentheses means $(e_o)$  .

\section{Obtaining operad structure in open-closed string topology}
As in the last section, we'll abuse notation and refer to the operad inside of the dioperad $OC$ as $OC$ also.  This section is an extension to an open-closed setting of the construction by Voronov \cite{vor1} in which he defined the Cacti operad, h.e. to the operad formed by components of $C$ in sect. 3 having only one output, and showed how it produces the BV-structure on the homology of a free loop space LM of a compact oriented manifold $M$ of dimension $m$ given by Chas-Sullivan \cite{chas1}.  The goal is to define a 2-colored version of Cacti h.e. to the operad $OC$ and use it to realize the pair $H_*(LM),H_*(PM_K))$ as an algebra over $H_*(OC)$ where $K$ is a fixed oriented closed submanifold of dimension k and $PM_K$ is the space of paths starting and ending in $K$.

Let's first recall what the Cacti operad is and how we get $H_*(LM)$ as an algebra over its homology.  Basically, Cacti is what results when you get rid of everything in $C$ except the boundaries of the closed inputs and outputs.  A typical point in Cacti(n) is shown in the following picture:
\begin{center}
\includegraphics[width=2.5in]{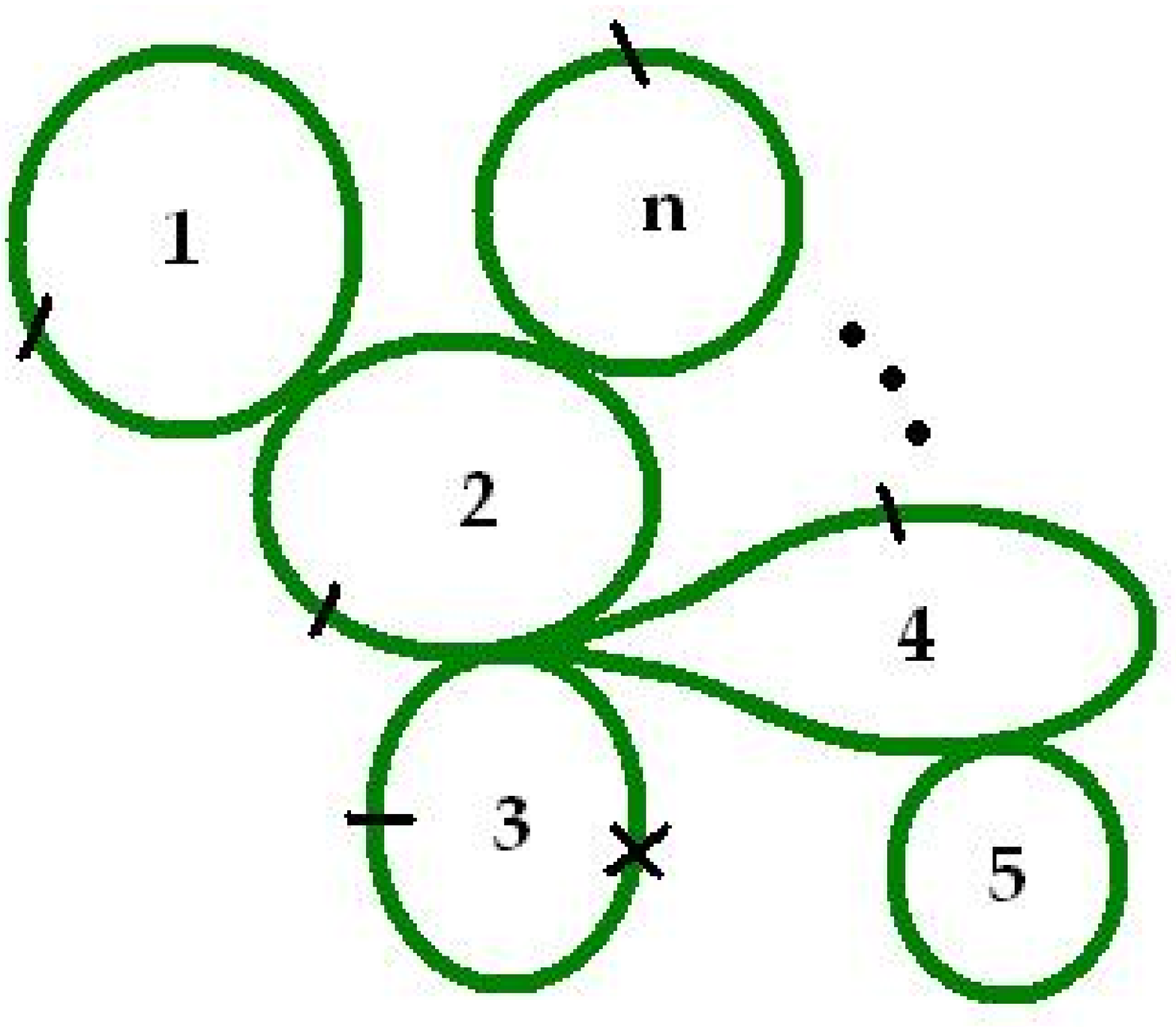}
\end{center}
\underline{\bf Explanation}:  A point is given by a tree like configuration of n labeled circles (lobes) of varying radii.  Each circle is parameterized by marking a point on it.  When 3 or more circles intersect at one point, they are given a cyclic order (can draw using counter clock-wise orientation of the plane).  Finally, if you start at any point on some circle and trace the picture in the counter clock-wise direction, using the cyclic orderings to jump from circle to circle, then the entire picture will be traversed before returning to the starting point.  Thus putting one more marked point on the picture gives it an $S^1$ parametrization and we consider the whole boundary as the output.  To compose two pictures, just replace the input circle of one picture with the entire second picture by rescaling the length of the second picture to math the length of the input and then identifying them via their parameterizations.\\

If we consider the space of maps from a point in $Cacti(n)$ into our manifold $M$, then restricting to the inputs gives us an embedding from this space of maps to $LM^n$ of finite codimension $(n-1)m$.  Restricting to the output gives a map into $LM$.  Thus we get the following diagram:\\
$Cacti(n) \times LM^n \stackrel{in}{\longleftarrow} L^{Cacti(n)}M \stackrel{out}{\longrightarrow} LM$ where $L^{Cacti(n)}M$ is the space of pairs $\{(c,f)|c \in Cacti(n), f \in Maps(c,M)\}$ and $in$ has finite codimension $(n-1)m$.

Applying the Pontryagin-Thom construction to the map $\in$ to get the push-forward map in homology, and then composing, we get the action $H_*(Cacti(n)) \otimes H_*(LM)^{\otimes n} \mapsto H_*(LM$.  This map has degree $-(n-1)m$ and gives an operad morphism, i.e. it commutes with composition and is equivariant.  The operations corresponding to the generators of $H_0(Cacti(2))$ and $H_1(Cacti(1))$ are exactly the BV operations of Chas-Sullivan.

In our situation, if we consider maps from a point in $OC$ to the manifold $M$ such that all boundary components map into the submanifold $K$, then restricting to the $S^1$ boundary of a closed input/output gives a point in $LM$ and restricting to the arc of the boundary of an open input/output gives a point in $PM^K$.  Thus we get a similar diagram as above with $OC,LM,PM_K$ and we just need to replace $OC$ by a skeletal model which will make the left arrow a finite codimensional embedding.

Right away we see, however, that we can not 'contract away' enough of a surface to get a finite codimensional embedding and get a space h.e. to $OC$ since we have these boundary components with no open inputs or outputs on them.  It is possible to define an operad which ignores these empty boundary components and gives us finite codim. embeddings, but the resulting action does not commute with composition.  This can be seen by considering the degrees of the induced operations.  Thus we are forced to keep the empty boundary components in our picture.  This in turn forces us to keep the part of the boundaries in between two open inputs so that composition gives an empty boundary when it should.

To handle the fact that this prevents us from directly getting finite codim. embeddings, consider the following ($PK$ is the space of paths in $K$):

$$ \begin{array}{r} PM_K \times PM_K \longleftarrow \\ \vspace{.42in} \\ PM_k \times PM_K \times PK \times PK \longleftarrow \end{array}
\put(0,-52){\includegraphics[height=1.5in]{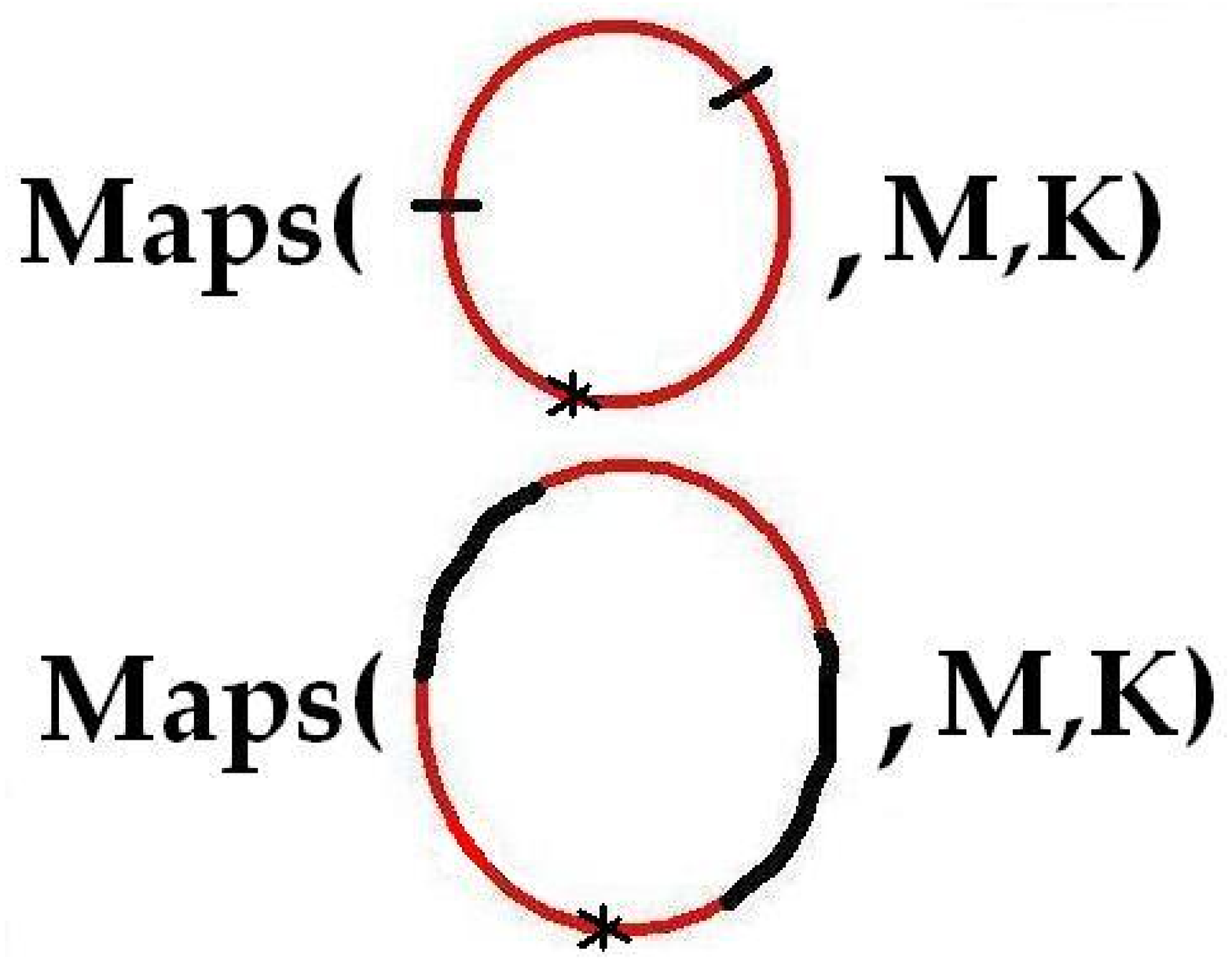}} \hspace{2in}
\begin{array}{l} \longrightarrow LM \\ \vspace{.42in} \\ \longrightarrow LM \end{array}$$

(The maps must send everything in black (except the main marked point) into the submanifold $K$)

Both of these maps are of finite codimension, the first of codim $2k$ and the second of codim $4k$.   So we can get the push forward maps in homology and compose to get operations.  The key observation is that if we plug the fundamental class of $H_k(PK)$ into the $H_*(PK)$'s then the two resulting operations are the same degree $-2k$ operation. This is the operation which is induced at the chain level by the function which takes two cells in $PM_K$ and transversally intersects the endpoints of the intervals of the first cell with the beginning points of the second cell and transversally intersects the beginning points of the first cell with the endpoints of the last cell.  This results in a chain in $LM$ of dimension $2k$ less than the sum of the dimensions of the two cells.

For the next observation, consider the operations given by the following two pictures:
\begin{center}
\includegraphics[width=4in]{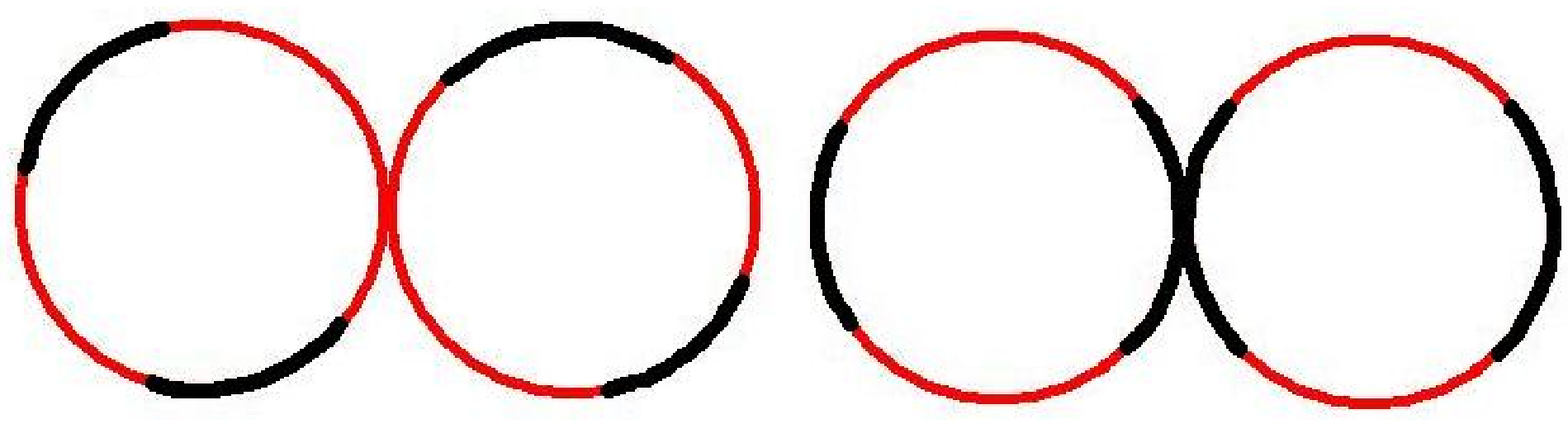}
\end{center}

The left picture results in an operation of degree $-(4k+m)$ while the right one gives a $-5k$ degree operation.  This is a problem since these two pictures would be in the same path component of our potential colored Cacti.  To remedy this, "ghost edges" are introduced:
\begin{center}
\includegraphics[width=4in]{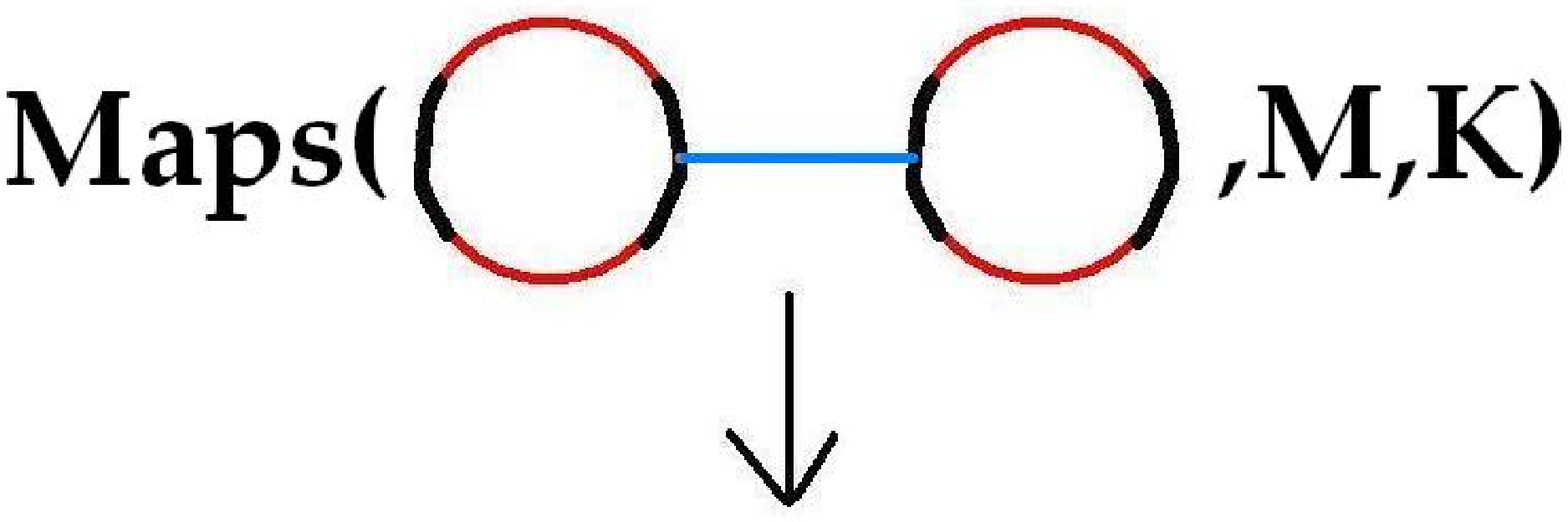}
\end{center}
\begin{center}
$PM_K \times PM_K \times PM_K \times PM_K \times PK \times PK \times PK \times PK \times PM$
\end{center}

This map has codimension $8k+2m$ no matter how the two circles are connected by the ghost edge.  If  we plug in the fundamental classes in $H_k(PK)$ and $H_m(PM)$ after getting the push forward, then we get a degree $-(4k+m)$ operation which is the same as the operation given by the left picture above.  This is the operation which takes 4 cells in $PM_K$, applies the operation described above to the first two and the last two resulting in two chains in $LM$, and then takes the loop product of these two chains as in Chas-Sullivan.\\
These two observations motivate the following:
\begin{df} {\bf\em Open-Closed Cacti}\\
The definition will be given by considering the following pictures which show show typical points in the configuration space.  First consider the case where the output is closed:
\begin{center}
\includegraphics[width=5in]{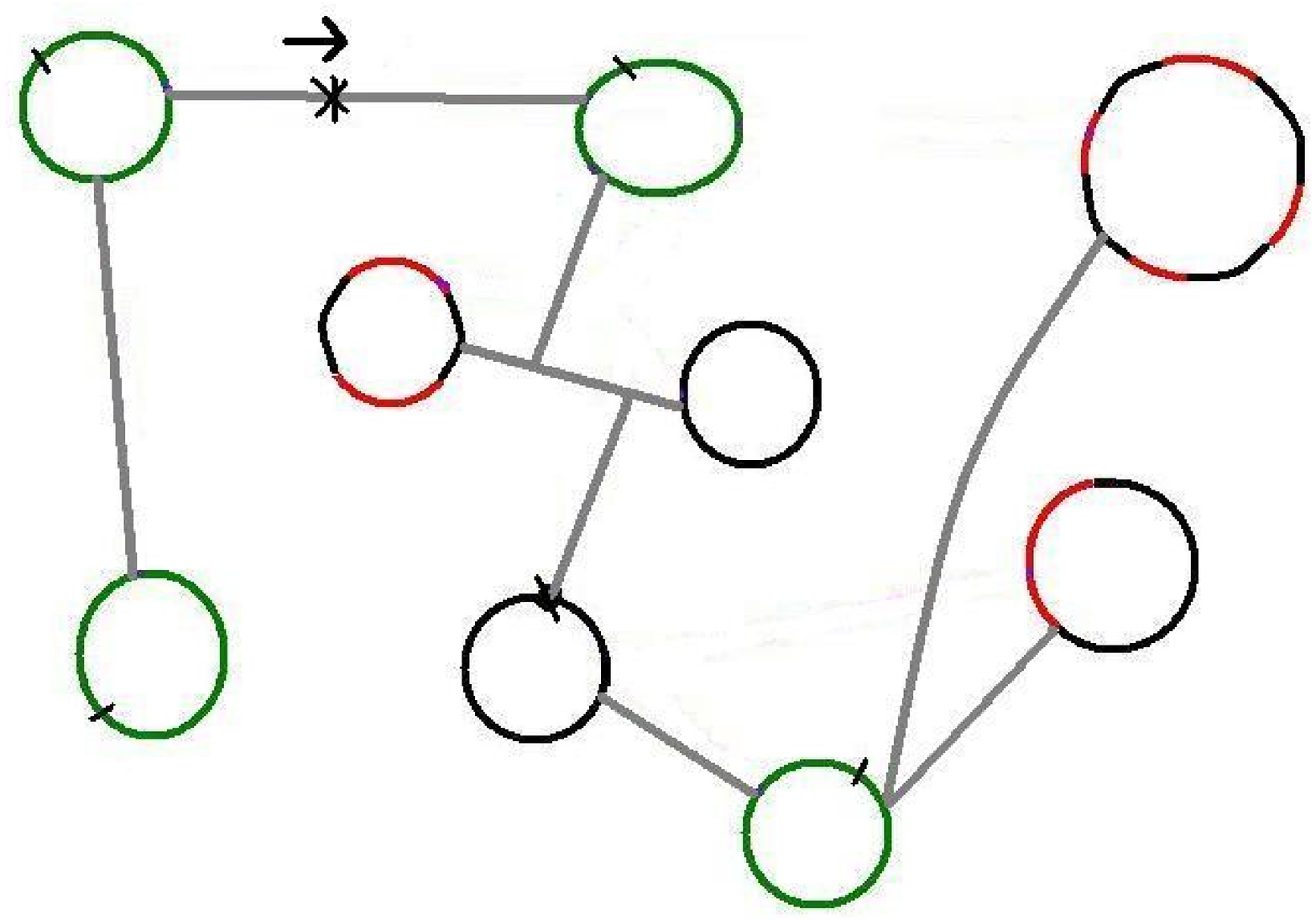}
\end{center}

\underline{Explanation}:\\
--The inputs are labeled.  Green circles are closed inputs and red intervals are open inputs.\\
--There is a marked point somewhere on the picture giving the starting point of the output.  If the marked point is on a circle, go in the counter clockwise direction.  If it's on a ghost edge, there needs to be an arrow pointing in the starting direction.  Then if we take the cyclic ordering of the edges meeting at a vertex to be given by counter clockwise orientation of the plane, there is a $S^1$   parametrization of the picture as in cacti.\\
--The closed inputs should have a mark as in cacti.\\
--For a black circle (empty boundary component) with more than one vertex, mark one of the vertices.  This marks where a boundary edge (black interval) "sews up" into a circle when composing.\\
--We can make our definition so that there is always $n-1$ ghost edges when there are $n$ circles.  We need this so that there are the same number of $PM$'s  to map into as above for any two pictures in the same path component.  This is done by choosing $t-2$ of the rays emanating from an intersection vertex which is not on a circle and has $t$  ghost rays emanating from it.  For the chosen $t-2$ rays, the rays are considered as ghost edges for which this vertex is an endpoint.  The other two rays are considered as one ghost edge and this vertex is just in the middle of it.  For example, in the above picture, there are three ghost edges, not 5, connecting the four circles in the middle.  The following picture shows the three different ways to connect three circles with 2 ghost edges which intersect off of the circles and the path in the configuration space which connects them.  After seeing this example (and after seeing the picture for the case of an open output), it is not hard to see why the path components of Open-Closed Cacti are in correspondence with the path components of $OC$:

\begin{center}
\includegraphics[width=3.5in]{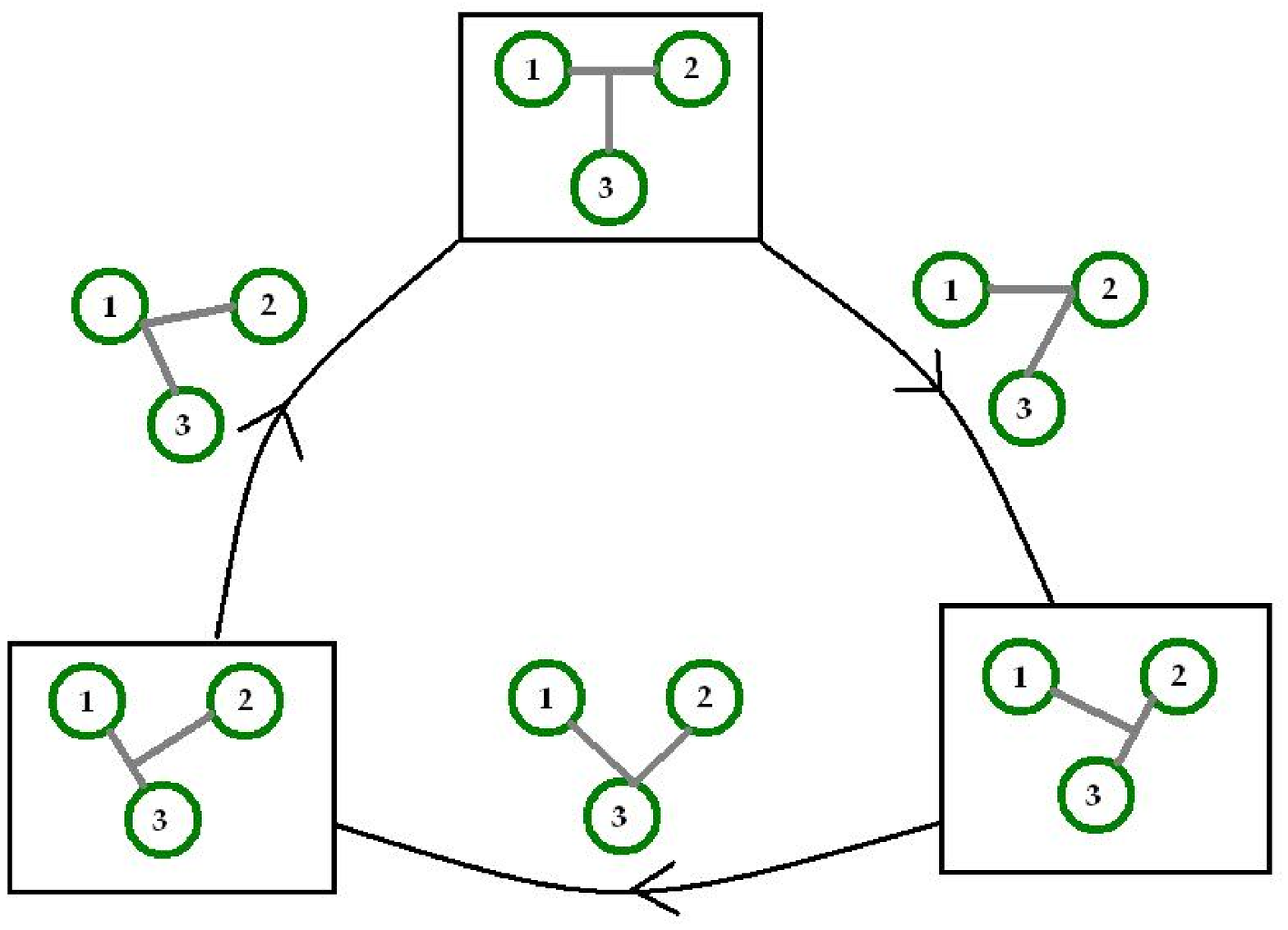}
\end{center}

Now let's look at the picture showing a typical point in the config space of O-C Cacti which has an open output:

\begin{center}
\includegraphics[width=5in]{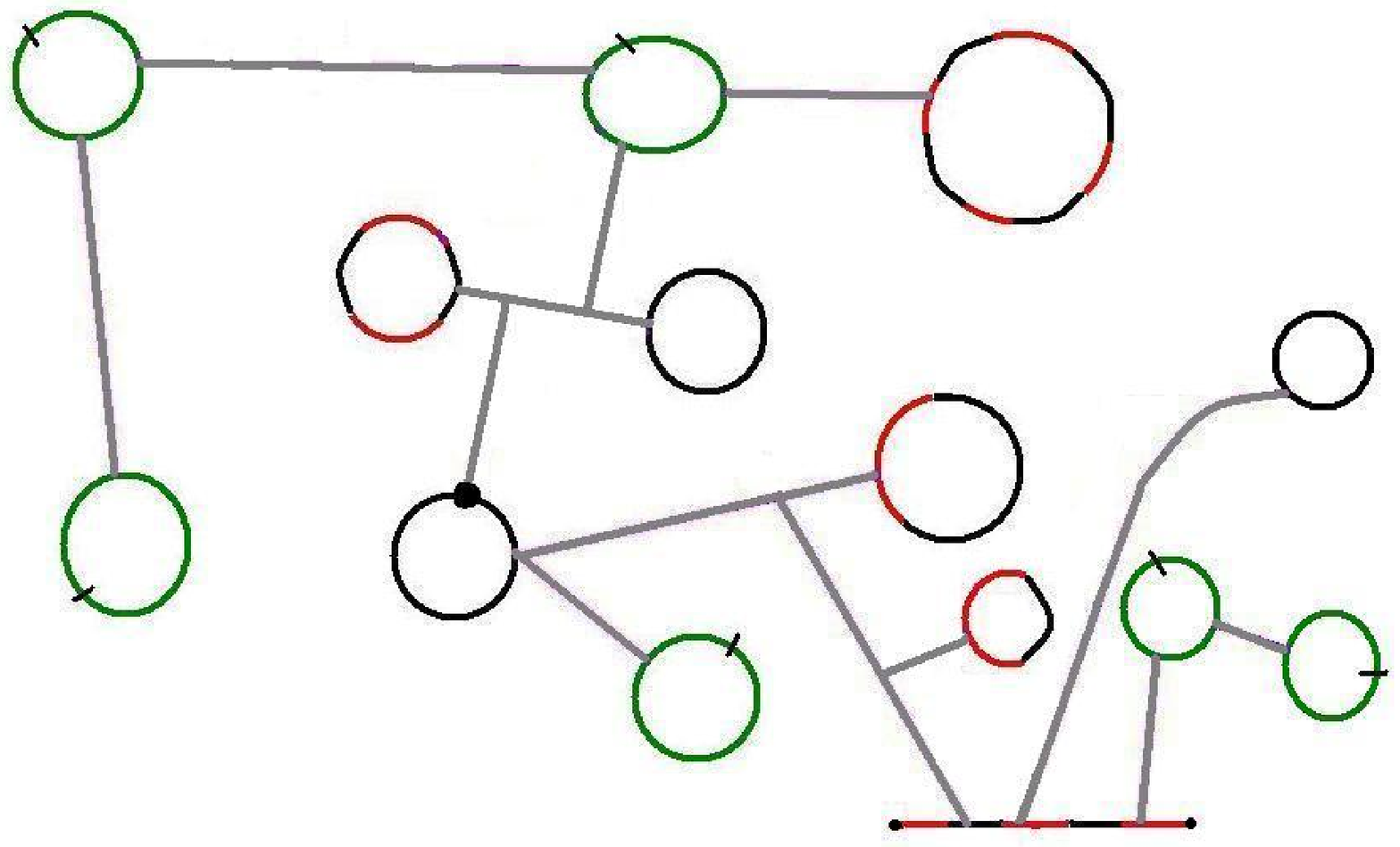}
\end{center}

\underline{Explanation}:\\
--Here, the number of ghost edges is equal to the number of circles.\\
--The interval at the bottom corresponds to boundary component containing the open output.\\
--If we start at the right of the interval and use the counter clock-wise orientation of the plane, we traverse the entire picture and end up at the left endpoint of the interval.  This gives us a parametrization of the picture by the unit interval and will serve as the output.\\
--If there are no open inputs on the same boundary comp as the output, then this interval should just be a black point (which should be sent to the submanifold $K$ when we consider maps into our manifold $M$ to get an action).

We can compose pictures by identifying the parametrization of a input circle or interval with the rescaled output parametrization of an entire picture, as in Cacti.  If the endpoint of a ghost edge lies on the input circle or interval, and if when we compose it happens that this endpoint is connected to the interior of a 2nd ghost edge in the picture we are replacing the input with, then this should not break the 2nd  ghost edge into two edges.  As described above, this endpoint just lies in the middle of this 2nd ghost edge.   Also, when connecting the endpoint to a 2nd ghost edge (interior or endpoint), there is natural way to cyclically order the ghost rays emanating from this point of attachment in the result.  Just "stick" the cyclic order of all the ghost edges in the input picture for which this point on the input circle (or interval) is an endpoint into the cyclic order of  the ghost rays emanating from this point of attachment in the output picture in between the ray which comes right before this point in the output parametrization and the ray which comes right after this point in the output parametrization.  For example, the following picture show the 3 different ways the ghost edges ending at a point on the input circle can be attached to the intersection point of two ghost edges:

\begin{center}
\includegraphics[width=5in]{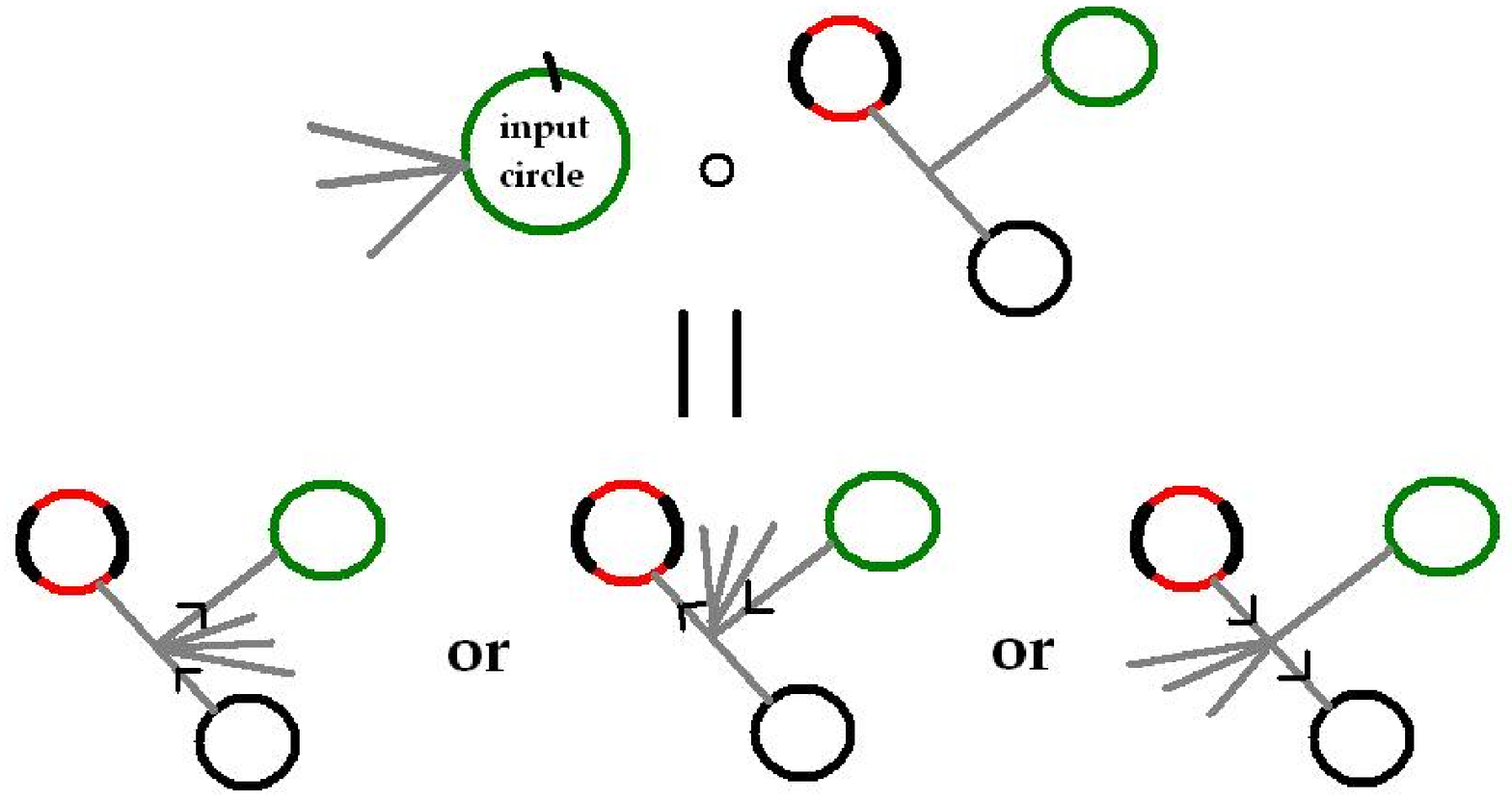}
\end{center}

Thus we see that after giving this configuration space its natural topology, as with Cacti or the Sullivan Chord Diagrams of \cite{cohen}, it forms a well-defined topological 2-colored operad which we'll call {\em Open-Closed Cacti} or {\em OC Cacti}.
\end{df}
\begin{prop}
$H_*(OC\,Cacti)$ and $H_*(OC)$ are isomorphic operads
\end{prop}
\begin{proof}
It is clear that the path components are in bijective correspondence.  Fix a path component of $OC\, Cacti$.  In a similar manner as in the proof of the homotopy type of a path component of $OC$ (dioperad) in Sec. 4, assume first that this path component has $n$ green circles (closed inputs), $k$  circles with exactly one open input, and $l$ completely black circles (empty boundary components).   If the path component has an open output then assume that there are no open inputs on the interval (so the interval is really a point), i.e. the boundary component possessing the open output does not have an open input for the points in the corresponding path component of $OC$.

Then this path component is h.e. to $Cacti(n+k+l)$ except that $l$ of the circles are not labeled and do not have a marked point.  And we know \cite{vor1} that Cacti is h.e. to $fD$, the framed little disks operad.  Thus it follows that this path component is h.e. to the configuration of $n+k+l$ disks inside the unit disk such that $n+k$ of the disks have labels and directions and l of the disks have neither directions nor labels.  But this is the homotopy type of the corresponding path component in $OC$, as described in Sec. 4.
Finally, just as in $OC$, the forgetful map from an arbitrary path component to the space obtained from this path comp. by dropping all but one fixed open input on each of the circles with 2 or more open inputs and dropping all inputs on the interval is a fibre bundle with contractible fibre.

Thus we have a vector space isomorphism from $H_*(OC)$ to $H_*(OC\,Cacti)$ given by these homotopy equivalences of the path components.   But it is clear that this map respects compositions of generators.  So they are isomorphic as colored operads.
\end{proof}

\newpage
\noindent \underline{\bf Obtaining action of $H_*(OC\,Cacti)$ on $(H_*(LM),H_*(PM_K))$}\\

First, let's consider the diagram associated to a fixed point $x \in OC\, Cacti$.  If $c_x$ is the number of closed inputs, $o_x$  the number of open inputs, $b_x$ the number of black intervals + the number of black circles, and $g_x$ the number of ghost edges, then we have:\\
$LM^{c_x} \times PM^{o_x}_K \times PK^{b_x} \times PM^{g_x} \stackrel{in}{\longleftarrow} Maps(x,M,K) \stackrel{out}{\longrightarrow} LM$ or $PM_K$

\noindent The maps should send all black intervals and circles into the submanifold $K$.  If the output is open, then the endpoints of the interval should also be sent into $K$ (if this interval is really a point, then this point should go into $K$).  Remember that one of the vertices of a black circle is marked so that a map restricted to this circle gives us a point in $LK \subset PK$.

The map $in$ is an embedding of finite codimension.  So we can get the push forward map in homology and then plug the fundamental classes into the $H_*(PK)$'s and $H_*(PM)$'s to obtain an operation \\
$H_*(LM)^{\otimes c_x} \otimes H_*(PM_K)^{\otimes o_x} \longmapsto H_*(LM)$ or $H_*(PM_K)$.  The codimension can be computed and thus the degree of this operation.  In the case of a closed output, the degree is $-[(\#\,of\,circles-1)m+(\#\,of\,open\,inputs)k]$.  In the case of an open output the degree is  $-[(\#\,of\,circles)m+(\#\,of\,open\,inputs-1)k]$.  These degrees are exactly what they should be.  That is, if we build the path component of $x$ out of the generators of $H_0(OC\,Cacti)$  via a tree in normal form, then the degree of the operation given by $x$ agrees with the degree of the operation given by composing in the same way the operations in $End(H_*(LM),H_*(PM_K))$ corresponding to points in the generating path components.

Next, consider the path component $P$ where $x \in P$.   We "almost" have the following diagram:\\
$P \times LM^{c_x} \times PM^{o_x}_K \times PK^{b_x} \times PM^{g_x} \stackrel{in}{\longleftarrow} L^P(M,K) \stackrel{out}{\longrightarrow} (LM$ or $PM_K)$, where $L^P(M,K)=\{(x,f)|x \in P,f \in Maps(x,M,K)\}$\\
The issue is that there are choices that would have to be made for each $x \in P$ in order to get a map $in$ and it is not clear that this can be done in a continuous way.  For the black intervals of the picture, there is a canonical way to choose the starting point of the interval and which $PK$ it corresponds to.  But for the ghost edges there is no canonical way to choose which endpoint is its starting point and which $PM$ it corresponds to.  Also, there is no canonical way to choose which $PK$  a black circle corresponds to.

There is no way to add these labels to our operad such that it doesn't matter the order in which we compose $n$ pictures into the $n$ inputs of a picture.  However, it is not necessary to change the definition of $OC$ Cacti.  Let $\bar{P}$ be the space obtained from $P$ by labeling the ghost edges and black circles, and marking one endpoint on each ghost edge.  Then $\bar{P} \mapsto P$ is a quotient map by a finite and free action so that the induced map $H_*(\bar{P}) \mapsto H_*(P)$  is onto (remember we are over a field of characteristic 0).  In fact, it should give an isomorphism when restricted to any path component in $\bar{p}$  (to see $\bar{P}$ might not be path connected, just consider the case of two green circles with one ghost edge connecting them).

Now, we do have the above diagram for $\bar{P}$.  So we can get the push forward and plug in fundamental classes to obtain\\
$H_*(\bar{P}) \otimes H_*(LM)^{\otimes c_x} \otimes H_*(PM_K)^{\otimes o_x} \longrightarrow H_*(LM)$ or $H_*(PM_K)$.\\
But it is clear that all the elements in $H_*(\bar{P})$ which are in the preimage of one element in $H_*(P)$ give the same operation since we are plugging in fundamental classes.  Thus we obtain our desired action.

\begin{thm}
This action of $H_*(OC\,Cacti)$ on $(H_*(LM),H_*(PM_K))$ is an operad action, i.e. it is equivariant and respects composition.
\end{thm}
\begin{proof}
We'll extend the proof given by Voronov in \cite{cohen2} for closed Cacti.  We consider the {\em category Corr of correspondences} where the topological actions described above are morphisms.
\begin{df} {\em Corr}\\
The objects are topological spaces, and a morphism (correspondence) between two spaces $X$ and $Y$ is a diagram $X \leftarrow X^{\prime} \rightarrow Y$ of continuous maps for some space $X^{\prime}$.  Two correspondences $X \leftarrow X^{\prime} \rightarrow Y$ and $Y \leftarrow Y^{\prime} \rightarrow Z$ are composed by taking a pullback:
\begin{equation*}
\begin{CD}
X' \times_Y Y' @>>> Y'\\
@VVV @VVV\\
X' @>>> Y
\end{CD}
\end{equation*}
which defines a new correspondence $X \leftarrow X' \times_Y Y' \to Z$.\\
Note that a map $X \rightarrow Y$ can be considered as a correspondence $X \stackrel{id}{\leftarrow} X \rightarrow Y$
\end{df}

Now, let $P_1$ be a path component in $OC\, Cacti(n)$, $P_2$ a path component in $OC\, Cacti(m)$ and let $P_3$ be the path component of $OC\, Cacti(n+m-1)$ such that $P_1 \circ_i P_2 \subset P_3$. Then we have the following diagram of correspondences:
\begin{equation*}
\begin{CD}
\bar{P_1} \times \bar{P_2} \times LM^{c_{\bar{P_3}}} \times PM_K^{o_{\bar{P_3}}} \times PK^{b_{\bar{P_3}}} \times PM^{g_{\bar{P_3}}} @>{\circ_i \times \id}>> \bar{P_3} \times LM^{c_{\bar{P_3}}} \times PM_K^{o_{\bar{P_3}}} \times PK^{b_{\bar{P_3}}} \times PM^{g_{\bar{P_3}}}\\
@VVV @VVV\\
\bar{P_1} \times LM^{c_{\bar{P_1}}} \times PM_K^{o_{\bar{P_1}}} \times PK^{b_{\bar{P_1}}} \times PM^{g_{\bar{P_1}}} @>>> (LM\,or\,PM_K)
\end{CD}
\end{equation*}
The top map is given by operad composition and is just a regular map considered as a correspondence.\\
Then, just as in the proof of the closed case, we can see that the diagram commutes by composing correspondences and seeing that both compositions are equal to the following correspondence:
$\bar{P_1} \times \bar{P_2} \times LM^{c_{\bar{P_3}}} \times PM_K^{o_{\bar{P_3}}} \times PK^{b_{\bar{P_3}}} \times PM^{g_{\bar{P_3}}} \leftarrow L^{\bar{P_1} \circ_i \bar{P_2}}(M,K) \rightarrow (LM$ or $PM_K)$, where $L^{\bar{P_1} \circ_i \bar{P_2}}(M,K)=\{(x,y,f)|x \in P_1,y \in P_2, f \in Maps(x \circ_i y,M,K)\}$.
The left arrow of this correspondence is a finite codimensional embedding whose codimension is the sum of the codimensions of two correspondences being composed.  We can then apply the Pontryagin-Thom construction to get push forward maps in homology and use the functoriality of the homology with respect to this construction along with its naturality on pull back diagrams to get a commutative diagram in homology.  Finally, we plug in the fundamental classes $e_o \in H_*(PK)$ and $e_c \in H_*(PM)$ and replace $\bar{P_i}$ with $P_i$ to see that our action of $H_*(OC\,Cacti)$ on $(H*(LM),H_*(PM_K))$ commutes with composition.
\end{proof}

Let's end by looking at this action when the two pictures being composed result in a black circle, (corresponding to a surface with a boundary component with no open inputs or outputs on it) to get an understanding of why this idea of adding a $PK$ for every black interval and circle, and then plugging in the fundamental class in homology, is the right thing to do.  We know that in the moduli space these boundary components correspond to $\phi_{o \mapsto c}(e_o)$.  For example, the normal form tree on pg. 23 corresponds to an element in the homology of a path component of $OC$ whose points have two empty boundary components.

The simplest case where this situation arises is when compose the following degree 0 generators (the pictures represent the path components that they belong to):\\
\includegraphics[height=4in]{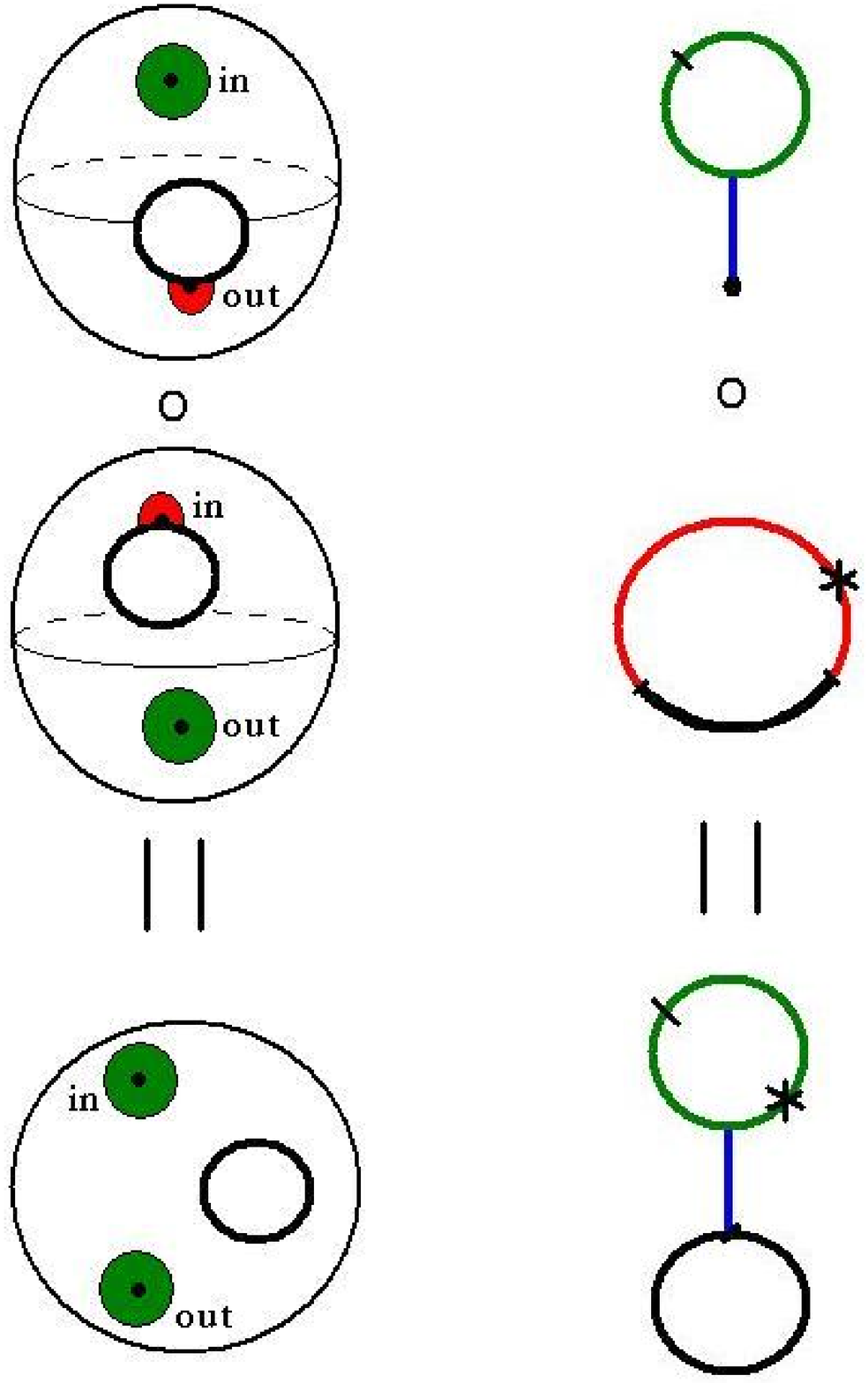} \hspace{.5in} $\begin{array}[b]{l} {\Large\phi_{c \mapsto o}:H_*(LM) \mapsto H_*(PM_K)} \\ \vspace{1in} \\ {\Large \phi_{o \mapsto c}:H_*(PM_K) \mapsto H_*(LM)} \\ \vspace{1in} \\ {\Large \phi_{o \mapsto c}(\phi_{c \mapsto o}):H_*(LM) \mapsto H_*(LM)} \\ \vspace{.5in} \end{array}$

The top operation comes from the correspondence whose left arrow goes from (the space of maps of the cacti picture into $M$ such that the black vertex at the bottom of the ghost edge goes into $K$) to $LM \times PM$.  This map has codimension $m + (m-k)$.  So after getting the push forward map and plugging in $e_c \in H_m(PM)$ we get a degree $-(m-k)$ operation.  This operation is induced from the chain operation which takes a cell in $LM$ and transversally intersects the marked points of the loops with $K$

The left arrow of the correspondence for the middle operation goes into $PM_K \times PK$ and has codim $2k$.  So after plugging in $e_o \in H_k(PK)$ we get a degree $-k$ operation.  The picture tells us that this is the operation which transversally intersects the two endpoints of each path in a cell of $PM_K$.

The left arrow for the correspondence giving the operation of the composition picture goes into $LM \times PM \times PK$.  It has codimension $m+m+k$ since $(x,y,z) \in LM \times PM \times PK$ is in the image iff $x(t_0)=y(0)$ (codim $m$), $y(1)=z(0)$ (codim $m$) and $z(0)=z(1)$ (codim $k$).  The operation $H_*(LM) \otimes H_*(PK) \mapsto H_*(LM)$, obtained by plugging in $e_c \in H_m(PM)$, transversally intersects the two endpoints of the paths in $K$ then takes the loop product of the result with the loops in $M$.  Thus plugging in $e_o \in H_k(PK)$ we get the degree $-m$ operation that takes in $a \in H_*(LM)$ and outputs $a\phi_{o \mapsto c}(e_o)$.  And this is exactly what should result since, in $H_*(OC)$, the normal form tree representing the operation $\phi_{o \mapsto c}(\phi_{c \mapsto o}(a))$ is the tree giving the operation $a\phi_{o \mapsto c}(e_o)$ (these two operations being the same is implied by relation $r4$).

{\bf Last Remark:} The constructions of this section can be easily extended to the situation where we consider a set of submanifolds of $M$, $\{K_{\lambda}\}_{\lambda \in \Lambda}$ and the spaces $PM_{(K_{\lambda},K_{\lambda^{\prime}}})$ of paths starting in $K_{\lambda}$ and ending in $K_{\lambda^{\prime}}$.  We get an action on the homology of these spaces by labeling the black intervals and circles of $OC \,Cacti$ by elements in the set $\Lambda$, called D-branes in the literature.  An open input can only be composed with an open output if their pair of D-branes match.  Thus we get a colored operad, one color for each pair of D-branes (repeats allowed), and the homology of this operad acts on the homology of the loop space and all of the path spaces.

Also, I have come up with a graph model extending $OC\, Cacti$ to a full PROP which acts in string topology.  Both of these extensions will be included in my thesis.  As mentioned in the introduction, I have learned recently that Antonio Ramirez has independently done similar work involving an open-closed graph PROP and string topology in his thesis.

\end{document}